\documentclass[preprint,11pt,3p,times]{article}

\usepackage[amssymb,thinqspace,thinspace]{SIunits}
\usepackage{amsmath}
\usepackage{amsfonts}
\usepackage[mathscr]{eucal}
\usepackage{graphicx}
\usepackage{subfigure}
\usepackage{color}
\usepackage[ruled]{algorithm2e}

\usepackage{authblk} 

\usepackage[left=2.2cm,right=2.2cm,top=2.0cm,bottom=2.0cm]{geometry}

\usepackage{hyphenat}

\usepackage{multirow}

\usepackage{tikz}

\usepackage{pgfplots}
\pgfplotsset{compat = newest}

\usepackage[colorlinks=true,
            linkcolor=black,
            citecolor=black,
            urlcolor=black,
            bookmarksnumbered=true,
            pdfstartpage=1,
            pdfpagemode=UseOutlines]{hyperref}

\newcommand{\mycite}[1]{} 
\newcommand{\myciteDetail}[2]{} 

\newcommand{\secref}[1]{Section~\ref{#1}} 
\newcommand{\secsref}[2]{Sections~\ref{#1} and~\ref{#2}} 
\newcommand{\secsrangeref}[2]{Sections~\ref{#1} -- \ref{#2}} 
\newcommand{\figref}[1]{Figure~\ref{#1}} 
\newcommand{\figssref}[3]{Figures~\ref{#1}, \ref{#2}, and~\ref{#3}} 
\newcommand{\Figref}[1]{Figure~\ref{#1}} 
\newcommand{\figsrangeref}[2]{Figures~\ref{#1} -- \ref{#2}} 
\newcommand{\tabref}[1]{Table~\ref{#1}} 
\newcommand{\algref}[1]{Algorithm~\ref{#1}} 
\newcommand{\Algref}[1]{Algorithm~\ref{#1}} 

\newcommand{\teo}[1]{\ensuremath{\underline{\mathbf{#1}}}} 
\newcommand{\mao}[1]{\ensuremath{\mathbf{#1}}} 
\newcommand{\mat}[1]{\ensuremath{\mathbf{#1}}} 

\newcommand{\mr}[1]{\ensuremath{\mathrm{#1}}} 

\newcommand{\mbs}[1]{\ensuremath{\boldsymbol{#1}}} 
\newcommand{\define}[1]{\emph{#1}} 

\newcommand{\indexedDom}[2]{#1^{\left(#2\right)}} 
\newcommand{\indexedNode}[2]{#1_{#2}} 

\newcommand{\sub}[2]{#1_{#2}} 

\newcommand{\Bound}{\ensuremath{\Gamma}} 
\newcommand{\DBound}{\ensuremath{\Bound_{\DCond}}} 
\newcommand{\NBound}{\ensuremath{\Bound_{\NCond}}} 
\newcommand{\CoupBound}{\ensuremath{\Bound_{\CoupCond}}} 
\newcommand{\CoupBoundDisc}{\ensuremath{\Bound_{\CoupCond,\Disc}}} 
\newcommand{\CoupBoundSubdomain}[1]{\ensuremath{\Bound_{\CoupCond,#1}}} 

\newcommand{\CoupCond}{\ast} 
\newcommand{\DCond}{\mr{D}} 
\newcommand{\NCond}{\mr{N}} 

\newcommand{\indexedCoupling}[1]{{#1}_{\CoupCond}} 

\newcommand{\indexedSlave}[1]{#1^\MoSlave}
\newcommand{\indexedMaster}[1]{#1^\MoMaster}

\newcommand{\bulk}[1]{#1^\dom}

\newcommand{\trans}[1]{#1^{\mr{T}}} 

\newcommand{\Disc}{\mr{h}} 
\newcommand{\disc}[1]{#1_{\Disc}} 

\newcommand{\meshfac}{\kappa} 

\newcommand{\dx}{\,\mr{d}} 

\newcommand{\bxi}{\mbs{\xi}}

\newcommand{\blamb}{\mbs\lambda}

\newcommand{\MapMaToSl}{\chi} 
\newcommand{\MoM}{\mbs{\mathscr M}} 
\newcommand{\MoD}{\mbs{\mathscr D}} 
\newcommand{\MoMaster}{{\mr{ma}}} 
\newcommand{\MoSlave}{{\mr{sl}}} 

\newcommand{\ttime}{\ensuremath{t}} 
\newcommand{\ttimeend}{\ensuremath{T}} 
\newcommand{\Dt}{\Delta t} 
\newcommand{\SpectralRadiusGenAlpha}{\rho_\infty} 

\newcommand{\AnyQuantity}{\left(\bullet\right)} 
\newcommand{\HOT}[2]{\mathcal{O}\left(#1^{#2}\right)} 

\newcommand{\matDiv}{\mathrm{Div}} 

\newcommand{\ShapeFunc}{N} 
\newcommand{\ShapeFuncLM}{\varPhi} 
\newcommand{\hEle}{h} 
\newcommand{\ParamCoord}{\xi} 
\newcommand{\ParamCoordVec}{\mbs\bxi} 

\newcommand{\dom}{\ensuremath{\Omega}} 

\newcommand{\matConfig}[1]{#1_{0}} 

\newcommand{\trac}{h} 
\newcommand{\normal}{n} 
\newcommand{\lmb}{\blamb} 
\newcommand{\lm}{\lambda} 
\newcommand{\res}{f} 

\newcommand{\matCoord}{X} 

\newcommand{\weakForm}{\delta\mathscr{W}} 
\newcommand{\WeakForm}[1]{\weakForm_{#1}} 

\newcommand{\bodyforce}{b} 
\newcommand{\density}{\rho} 
\newcommand{\coord}{x} 
\newcommand{\sol}{u} 
\newcommand{\bc}[1]{\ensuremath{\bar{#1}}} 

\newcommand{\variation}{\delta} 

\newcommand{\disp}{u} 
\newcommand{\pressure}{p} 
\newcommand{\pki}{P} 
\newcommand{\youngs}{E} 
\newcommand{\poisson}{\nu} 

\newcommand{\gap}{g} 

\newcommand{\average}[1]{\overline{#1}} 

\newcommand{\jacobian}{J} 

\newcommand{\indTimeStep}{n} 
\newcommand{\indDomain}{i} 
\newcommand{\indNode}{j} 
\newcommand{\indNodeSlave}{k} 
\newcommand{\indNodeMaster}{\ell} 
\newcommand{\indProc}{p} 
\newcommand{\altIndProc}{q} 
\newcommand{\indSubdomain}{m} 
\newcommand{\altIndSubdomain}{n} 
\newcommand{\indBin}{b} 

\newcommand{\indContactEvalEvent}{c} 
\newcommand{\numContactEvalEvents}{C} 

\newcommand{\indexedNormal}[1]{#1_{\normal}}

\newcommand{\Ghosted}{\mathrm{gh}}

\newcommand{\procID}[2]{#1_#2}
\newcommand{\subdomainID}[2]{#1_#2}

\newcommand{\procName}[1]{'proc~$#1$'}

\newcommand{\nproc}{n^\mr{proc}} 
\newcommand{\nsubdomain}{M} 
\newcommand{\ndim}{d} 
\newcommand{\nnode}{n^\mr{nd}} 
\newcommand{\nEle}{n^\mr{el}} 
\newcommand{\nEleMin}{n^\mr{el}_{\mr{min}}} 
\newcommand{\nNodesSlave}{\indexedDom{n}{1}} 
\newcommand{\nNodesSlaveWithLM}{\indexedDom{m}{1}} 
\newcommand{\nNodesMaster}{\indexedDom{n}{2}} 
\newcommand{\nEleSlave}{n^\mr{el,\MoSlave}} 
\newcommand{\nEleMaster}{n^\mr{el,\MoMaster}} 

\newcommand{\tGhosting}{\ttime_{\Ghosted}} 
\newcommand{\tRedistribute}{\ttime_{\mr{LB}}} 
\newcommand{\tEvaluate}{\ttime_{\mr{eval}}} 
\newcommand{\tEvaluateOnProc}[1]{\ttime_{\mr{eval}, #1}} 
\newcommand{\tAssemble}{\ttime_{\mr{ass}}} 
\newcommand{\tContactTotal}{\ttime_{\mr{total}}} 
\newcommand{\tAccumulated}{\ttime_{\mr{acc}}} 

\newcommand{\costStorage}{s} 
\newcommand{\costCommunication}{\sigma} 
\newcommand{\costCommunicationMeasure}{\zeta} 

\newcommand{\imbalanceThresholdTime}{\hat{\eta}_{\mr{t}}}

\newcommand{\imbalanceRatioTime}{\eta_{\mr{t}}}
\newcommand{\imbalanceRatioEles}{\eta_{\mr{e}}}

\newcommand{\listGhostEles}{\{e_{\Ghosted}\}}
\newcommand{\listOfBinsEnclosingSubdomain}[1]{\{\indexedSlave{\mathfrak{B}}_{#1}\}} 
\newcommand{\listOfNeighboringBins}[1]{\{\mathfrak{B}^{\mathrm{nb}}_{#1}\}} 

\newcommand{\binSize}{\beta}
\newcommand{\minBinSize}{\binSize_{\mathrm{min}}}

\newcommand{\inhouse}{in\hyp{}house}
\newcommand{\meshtying}{meshtying}

\newcommand{\multilevel}{multi\hyp{}level}

\newcommand{\multiphysics}{multi\hyp{}physics}

\newcommand{\nonlinear}{nonlinear}

\newcommand{\nonmatching}{non\hyp{}matching}

\newcommand{\runtime}{run\hyp{}time}

\newcommand{\speedup}{speed\hyp{}up}

\newcommand{\timetosolution}{time-to-solution}

\newcommand{\master}{master}

\newcommand{\slave}{slave}
\newcommand{\Slave}{Slave}

\newcommand{\SoftwarePackage}[1]{\textsc{#1}} 
\newcommand{\zoltan}{\SoftwarePackage{Zoltan}}
\newcommand{\trilinos}{\SoftwarePackage{Trilinos}} 
\newcommand{\baci}{\SoftwarePackage{Baci}} 

\newcommand{\reqSlave}{\textbf{R2}} 
\newcommand{\reqMaster}{\textbf{R1}} 

\newcommand{\ie}{i.e.}

\newcommand{\eg}{e.g.} 
 
\newcommand{\cf}{cf.}
\newcommand{\wrt}{w.r.t.}

\usetikzlibrary{arrows.meta}
\usetikzlibrary{bending}
\usetikzlibrary{patterns}
\usetikzlibrary{positioning}
\usetikzlibrary{snakes}
\usetikzlibrary{shapes.misc}
\usetikzlibrary{shapes.geometric}

\tikzstyle{domainline}=[thick]
\tikzstyle{gridline}=[thin]
\tikzstyle{hiddenline}=[thick,dashed]
\tikzstyle{dbcline}=[thick]
\tikzstyle{tractionline}=[thick]
\tikzstyle{tractionarrow}=[thick,-{Latex}]
\tikzstyle{measline}=[thick, {Latex}-{Latex}]
\tikzstyle{measradius}=[thick, -{Latex}]
\tikzstyle{measauxline}=[thick]
\tikzstyle{coordsline}=[thick,{Latex[open,fill=white]}-{Latex[open,fill=white]}]
\tikzstyle{coordline}=[thick,-{Latex[open,fill=white]}]
\tikzstyle{symmetryline}=[very thick,dash pattern=on 10pt off 3pt on \the\pgflinewidth off 3pt]

\tikzstyle{darkgreen}=[green!50!black]

\tikzstyle{springlineal}=[snake=zigzag,thick,line before snake=0.5cm,line after snake=0.5cm,segment length=6,segment amplitude=5,join=round]

\tikzstyle{flowchartarrow}=[thick,-{Latex}]

\usetikzlibrary{calc}

\newcommand{\convexpath}[2]{
[   
    create hullnodes/.code={
        \global\edef\namelist{#1}
        \foreach [count=\counter] \nodename in \namelist {
            \global\edef\numberofnodes{\counter}
            \node at (\nodename) [draw=none,name=hullnode\counter] {};
        }
        \node at (hullnode\numberofnodes) [name=hullnode0,draw=none] {};
        \pgfmathtruncatemacro\lastnumber{\numberofnodes+1}
        \node at (hullnode1) [name=hullnode\lastnumber,draw=none] {};
    },
    create hullnodes
]
($(hullnode1)!#2!-90:(hullnode0)$)
\foreach [
    evaluate=\currentnode as \previousnode using \currentnode-1,
    evaluate=\currentnode as \nextnode using \currentnode+1
    ] \currentnode in {1,...,\numberofnodes} {
-- ($(hullnode\currentnode)!#2!-90:(hullnode\previousnode)$)
  let \p1 = ($(hullnode\currentnode)!#2!-90:(hullnode\previousnode) - (hullnode\currentnode)$),
    \n1 = {atan2(\y1,\x1)},  
    \p2 = ($(hullnode\currentnode)!#2!90:(hullnode\nextnode) - (hullnode\currentnode)$),
    \n2 = {atan2(\y2,\x2)},  
    \n{delta} = {-Mod(\n1-\n2,360)}
  in 
    {arc [start angle=\n1, delta angle=\n{delta}, radius=#2]}
}
-- cycle
} 

\usetikzlibrary{decorations.text}

\graphicspath{{fig/}}

\makeatletter
\def\ps@pprintTitle{%
  \let\@oddhead\@empty
  \let\@evenhead\@empty
  \def\@oddfoot{\reset@font\hfil\thepage\hfil}
  \let\@evenfoot\@oddfoot
}
\makeatother

\begin{document}

\title{Scalable computational kernels for mortar finite element methods}

\author[1,2]{Matthias Mayr}
\author[2]{Alexander Popp}

\affil[1]{\small Data Science \& Computing Lab, Universit\"{a}t der Bundeswehr M\"{u}nchen,\\Werner-Heisenberg-Weg 39, D-85577 Neubiberg, Germany}
\affil[2]{\small Institute for Mathematics and Computer-Based Simulation, Universit\"{a}t der Bundeswehr M\"{u}nchen,\\Werner-Heisenberg-Weg 39, D-85577 Neubiberg, Germany}

\date{}

\maketitle

\begin{abstract}
Targeting simulations on parallel hardware architectures,
this paper presents computational kernels for efficient computations in mortar finite element methods.
Mortar methods enable a variationally consistent imposition of coupling conditions at high accuracy,
but come with considerable numerical effort and cost for the evaluation of the mortar integrals to compute the coupling operators.
In this paper, we identify bottlenecks in parallel data layout and domain decomposition
that hinder an efficient evaluation of the mortar integrals.
We then propose a set of computational strategies to restore optimal parallel communication and scalability for the core kernels devoted to the evaluation of mortar terms.
We exemplarily study the proposed algorithmic components in the context of three-dimensional large-deformation contact mechanics,
both for cases with fixed and dynamically varying interface topology,
yet these concepts can naturally and easily be transferred to other mortar applications, {\eg} classical {\meshtying} problems.
To restore parallel scalability,
we employ overlapping domain decompositions of the interface discretization independent from the underlying volumes
and then tackle parallel communication for the mortar evaluation by a geometrically motivated reduction of ghosting data.
Using three-dimensional contact examples,
we demonstrate strong and weak scalability of the proposed algorithms
up to 480 parallel processes as well as study and discuss improvements in parallel communication related to mortar finite element methods.
For the first time,
dynamic load balancing is applied to mortar contact problems with evolving contact zones,
such that the computational work is well balanced among all parallel processors independent of the current state of the simulation.

\emph{Keywords:} Mortar methods, contact mechanics, interface problems, parallel algorithms, finite elements, domain decomposition
\end{abstract}

\section{Introduction}
\label{sec:Intro}

Mortar finite element methods (FEM) are nowadays well established in a variety of application areas in computational science and engineering
as discretization technique for the coupling of {\nonmatching} meshes.
Their general applicability in a vast range of problems as well as their mathematical properties, {\eg} variational consistency,
make them one of the most popular choices among interface discretization techniques.
They are undoubtedly the most preferred choice for robust finite element discretization
in computational contact mechanics undergoing large deformations~\cite{DeLorenzis2014a,Puso2004a,Puso2004b,Wohlmuth2011a,Wriggers2006a}.
However, the numerical effort and computational cost is high and can be considered a bottleneck in many scenarios.
This paper discusses several performance challenges of mortar methods in the context of parallel computing
and proposes remedies to reduce the overall runtime, obtain optimal scalability as well as reduce parallel communication and memory consumption.
As a demanding prototype application,
several test cases from computational contact mechanics showcase the proposed algorithms and their impact on runtime and parallel scalability.

Originally being developed in the context of domain decomposition
for the weak imposition of interfacial constraints~\cite{Belgacem1999a,Bernardi1993a},
mortar methods soon became popular in {\meshtying}~\cite{Puso2004c,Puso2003a}
and contact mechanics problems~\cite{Belgacem1998a,Gitterle2010a,McDevitt2000a,Popp2009a,Popp2010a,Popp2012b,Puso2004b,Puso2008a,Yang2005a,Yang2008a}.
Recently, mortar methods for {\meshtying} problems have regained attention due to the rise of isogeometric analysis
and the need for isogeometric patch coupling~\cite{Dittmann2019a,Dittmann2020a,Dornisch2015a,Dornisch2017a,Hesch2012a,Wunderlich2019a,Zou2018a}.
A variety of papers discusses mortar methods in the context of isogeometric analysis for contact problems,
among them~\cite{Dittmann2014a,Duong2019a,DeLorenzis2012a,DeLorenzis2014a,Seitz2016a}.
Moreover, mortar methods have spread to other single-field problems,
{\eg} contact mechanics including wear~\cite{Farah2016b} or fluid dynamics~\cite{Ehrl2014a},
as well as a variety of surface-coupled {\multiphysics} problems,
among them fluid-structure interaction~\cite{Hesch2014a,Kloeppel2011a,Mayr2015a}
or the simulation of lithium-ion cells in electrochemistry~\cite{Fang2018a}.
Lately, also volume-coupled problems have been addressed by mortar methods~\cite{Farah2016a}.
Despite their significant computational cost,
the popularity of mortar methods over classical node-to-segment, Gauss-point-to-segment, and other collocation-based approaches
is based on their mathematical properties such as their variational consistency and stability.
Compared to two-dimensional problems,
an efficient mortar evaluation is much more critical in three-dimensional problems,
which are at the same time of great practical relevance in real-world applications.

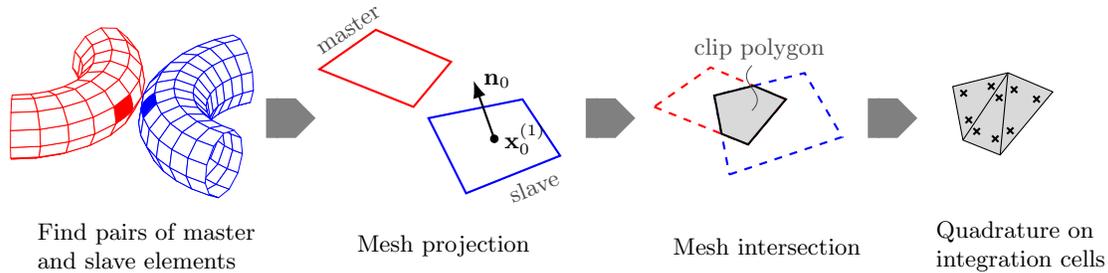
\begin{figure}
\begin{center}

\usetikzlibrary{shapes.misc}
\tikzset{cross/.style={cross out, draw=black, minimum size=2*(#1-\pgflinewidth), inner sep=0pt, outer sep=0pt, thick},
cross/.default={2pt}}

\begin{tikzpicture}[font=\footnotesize]

\def\xDistance{10cm} 
\def\yDistance{6cm} 

\def\normalStartX{4.25}
\def\normalStartY{1.95}
\def\normalEndX{3.7}
\def\normalEndY{3.5}

\def\fontColor{black!70}
\def\clipPolygonColor{gray!30}

\begin{scope}[scale=0.4,shift={(-2.3cm,-2.5cm)}] 

\begin{scope} 
\begin{axis}[
axis equal image,
hide axis,
z buffer = sort,
view = {122}{30},
scale = 1.5
]

\addplot3[ 
surf,
fill=white,
shader = faceted,
faceted color = red,
thick,
line join = round,
samples = 12,
samples y = 10,
domain = 0:2*pi,
domain y = 1/6*pi:5/6*pi
](
{(3+sin(deg(\x)))*cos(deg(\y))},
{(3+sin(deg(\x)))*sin(deg(\y))},
{cos(deg(\x))}
);

\addplot3[ 
surf,
fill = red,
shader = faceted,
faceted color = red,
thick,
samples = 12,
samples y = 10,
domain = 2/11*2*pi:3/11*2*pi,
domain y = 1/6*pi+4/9*4/6*pi:1/6*pi+5/9*4/6*pi,
](
{(3+sin(deg(\x)))*cos(deg(\y))},
{(3+sin(deg(\x)))*sin(deg(\y))},
{cos(deg(\x))}
);

\end{axis}
\end{scope} 

\begin{scope} 
\begin{axis}[
axis equal image,
hide axis,
z buffer = sort,
view = {122}{30},
scale = 1.5,
xshift = 3.4cm,
yshift = -0.7cm
]

\addplot3[ 
surf,
fill=white,
shader = faceted,
faceted color = blue,
thick,
line join = round,
samples = 12,
samples y = 14,
domain = 0:2*pi,
domain y = -0.7*pi:0
](
{(3+sin(deg(\x)))*cos(deg(\y))},
{(3+sin(deg(\x)))*sin(deg(\y))},
{cos(deg(\x))}
);

\addplot3[ 
surf,
fill = blue,
shader = faceted,
faceted color = blue,
thick,
samples = 12,
samples y = 14,
domain = 1/11*2*pi:2/11*2*pi,
domain y = -0.7*pi+7/13*0.7*pi:-0.7*pi+8/13*0.7*pi,
](
{(3+sin(deg(\x)))*cos(deg(\y))},
{(3+sin(deg(\x)))*sin(deg(\y))},
{cos(deg(\x))}
);
\end{axis}
\end{scope} 

\end{scope} 


\begin{scope}[scale=0.5, shift={(6.0,0)}]
\draw [-Triangle Cap,line width=15pt,draw=gray] (0,0) -- (1,0);
\end{scope}


\begin{scope}[scale=0.5, shift={(7.5,-2.5)}]

\draw [thick,color=blue] (3.5,0.5) -- (6,1.5) -- (5,3.0) -- (2.5,2.5) -- cycle;
\node [below left, rotate=22, text=\fontColor] at (6.0,1.5) {\slave};

\draw [measradius] (\normalStartX,\normalStartY) -- (\normalEndX,\normalEndY);  
\node [draw,circle,inner sep=1pt,fill] at (\normalStartX,\normalStartY) {};
\node [right] at (\normalStartX,\normalStartY) {$\indexedDom{\mao\coord_0}{1}$};
\node [right] at (\normalEndX,\normalEndY) {$\mao\normal_0$};

\begin{scope}[shift={(-1.4,1.5cm)}]
  \draw [thick,color=red] (1,2.3) -- (3.5,1.3) -- (4.5,2.5) -- (2.5,3.35) -- cycle;
  \node [above right, rotate=34, text=\fontColor] at (1,2.3) {\master};
\end{scope}
\end{scope}


\begin{scope}[scale=0.5, shift={(14.5,0)}]
\draw [-Triangle Cap,line width=15pt,draw=gray] (0,0) -- (1,0);
\end{scope}


\begin{scope}[scale=0.5, shift={(15.0,-2)}]

\draw [thick,dashed,color=blue] (3.0,0.5) -- (6,1.5) -- (5,3.2) -- (2.6,2.6) -- cycle;

\draw [thick,dashed,color=red] (1,2.3) -- (3.5,1.3) -- (4.5,2.5) -- (2.5,3.35) -- cycle;

\draw [thick,fill=\clipPolygonColor] (2.79381,1.58247) -- (3.5,1.3) -- (4.5,2.5) -- (3.64815,2.86204) -- (2.6,2.6) -- cycle;
\node [above, text=\fontColor] at (3.8,3.3) {clip polygon};
\draw [color=\fontColor] (3.6,2.2) .. controls (4.1,2.7) and (3.1,2.9) .. (3.6,3.4);

\end{scope}


\begin{scope}[scale=0.5, shift={(22.0,0)}]
\draw [-Triangle Cap,line width=15pt,draw=gray] (0,0) -- (1,0);
\end{scope}


\begin{scope}[scale=0.7, shift={(14.5,-2)}]

\draw [fill=\clipPolygonColor] (2.79381,1.58247) -- (3.5,1.3) -- (4.5,2.5) -- (3.64815,2.86204) -- (2.6,2.6) -- cycle;
\draw (2.79381,1.58247) -- (3.64815,2.86204);
\draw (3.64815,2.86204) -- (3.5,1.3);

\draw (2.903898,1.965320) node[cross] {};
\draw (3.331068,2.605105) node[cross] {};
\draw (2.806993,2.474085) node[cross] {};
\draw (3.053898,1.748653) node[cross] {};
\draw (3.406993,1.607418) node[cross] {};
\draw (3.481068,2.388438) node[cross] {};
\draw (3.691358,1.760340) node[cross] {};
\draw (4.191358,2.360340) node[cross] {};
\draw (3.765433,2.541360) node[cross] {};

\end{scope}


\begin{scope}[shift={(0,-1.7)}]

\node[text width=3cm] at (1.3,0) {Find pairs of {\master} and {\slave} elements};
\node at (5.2,0) {Mesh projection};
\node at (9.5,0) {Mesh intersection};
\node[text width=2.5cm] at (13,0) {Quadrature on integration cells};

\end{scope}

\end{tikzpicture}
\caption{Main steps of 3D mortar coupling (\emph{from left to right}):
Pairs of {\master} and {\slave} elements, that (i) are potentially in contact,
need to be (ii) projected onto each other along the normal vector~$\mao\normal_0$
to (iii) compute the mesh intersection
and to (iv) perform numerical quadrature of mortar contributions on integration cells.
}
\label{fig:MortarSegmentation}
\end{center}
\end{figure}

When using a Lagrange multiplier field~$\teo\lmb$ to impose constraints on the subdomain interfaces,
mortar methods discretize~$\teo\lmb$ on the so-called \emph{{\slave} side} of the interface.
The numerical effort of mortar methods is usually related to the search for nearest neighbors,
local projection of meshes and subsequent clipping and triangulation of intersected meshes,
as well as the resulting segment-based numerical integration, {\cf} \figref{fig:MortarSegmentation}.
While these operations themselves are already expensive,
implicit contact solvers need to perform them in \emph{every {\nonlinear} iteration},
rendering this a possible feasibility bottleneck or at least a performance impediment,
which becomes even more demanding through the necessity of consistent linearizations of all mortar terms.
The parallelization of contact search algorithms has been addressed in \cite{Hansen2016a} for example,
where standard domain-decomposition-based spatial search is enhanced with thread-level parallelism.
To {\speedup} the subsequent evaluation of contact terms,
various integration strategies are available,
among them \define{element-based} and \define{segment-based} integration,
{\cf}~\cite{Brivadis2015a,Farah2015a,Maday2002a,Wilking2017a}.
Segment-based integration subdivides each {\slave} element into segments having no discontinuities of the integrands within their domain.
This yields a highly accurate quadrature, though is computationally expensive.
Element-based integration on the other hand reduces the effort of clipping and triangulation the intersected meshes
by employing higher-order integration schemes to deal with weak discontinuities at element edges,
though brings along a less accurate evaluation of the mortar integrals.
While the segment-based integration strategy is unequivocally preferable due to its accuracy,
it comes at significantly higher computational cost.
Furthermore, systems of linear equations arising from mortar-based interface discretizations require
tailored preconditioning techniques for an efficient iterative solution procedure.
Depending on the specific details of the discretization, the resulting linear system might exhibit saddle-point structure.
Efficient preconditioners to be used in conjunction with Krylov solvers are available
in literature~\cite{Achdou1999a,Casarin1996a,Stefanica2001a,Wieners1999a,Wieners2003a,Wiesner2018a,Wiesner2021a,Wohlmuth2000b}
and, thus, are not in the scope of this paper.
We rather focus on the cost of evaluating all mortar-related terms.

As outlined previously,
many theoretical aspects of mortar methods have already been discussed and solved in the literature,
{\eg} the choice of discrete basis functions~\cite{Flemisch2007a,Lamichhane2005a,Lamichhane2007a,Popp2010a,Popp2012b,Wohlmuth2000a,Wohlmuth2012a},
numerical quadrature~\cite{Brivadis2015a,Farah2015a,Maday2002a,Wilking2017a},
conservation laws~\cite{Hesch2009a,Hesch2011a,Yang2005a},
or contact search algorithms~\cite{Benson1990a,Williams1995a,Williams1999a,Yang2008a,Yang2008b,Zhong1989a,Zhong1990a}.
However, computational aspects of mortar methods for contact problems
-- especially in the context of parallel computing --
have largely been neglected by the scientific community so far.
To fill this gap,
this work is motivated and guided by the quest for parallel scalability of \emph{all} algorithmic components of mortar methods
for arbitrarily evolving contact zones in three-dimensional problems.
Therefore, we analyze the computational kernels of mortar finite element methods and design their interplay to assure parallel scalability.
To the best of our knowledge,
most contributions in literature have focused on the serial case ({\ie} one processor) only
or have embedded mortar methods into existing parallel finite element codes without specific provisions.
An exception to this observation is the work of Krause and Zulian~\cite{Krause2016a},
where a parallel approach to the variational transfer of discrete fields between unstructured finite element meshes
as well as the associated proximity and intersection detections are described in detail
and examples for the evaluation of grid projection operators are given for various surface and volume projection problems.
Yet, Krause and Zulian~\cite{Krause2016a} spare dynamic contact problems with evolving contact zones,
which are of particular importance in engineering applications.
In the present contribution,
we analyze several schemes to subdivide mortar interface discretizations into subdomains suitable for parallel computing
and discuss their interplay with distributed memory architectures of computing clusters to achieve parallel scalability.
Thereby, we follow a message-passing parallel programming model
that utilizes the message passing interface (MPI) for communication between address spaces of different processes~\cite{mpi4_0}.
Finally, we develop and showcase a dynamic load balancing strategy
to address the particular needs of contact problems with evolving contact configurations and interface topologies for three-dimensional problems.

By starting from an analysis of the computational cost of the evaluation of mortar terms,
which is most commonly related to the {\slave} side of the contact interface,
we identify three main tasks, which will directly lead to the postulation of two essential requirements for parallel and scalable computational kernels for mortar finite element methods:
\begin{itemize}
\item For the geometrical task of identifying close {\master} and {\slave} nodes within the contact search,
each {\slave} node needs access to the position of every node of the {\master} side of the interface discretization.
While the distribution of the {\master} interface discretization to several compute nodes enables larger problem sizes,
it requires advanced ghosting ({\ie} sending data between different processors) of interface quantities to reduce the overall communication and memory footprint.
We will propose ghosting strategies that take a measure of geometric proximity between {\master} and {\slave} nodes into account
to pre-compute and reduce the list of {\master} nodes/elements to be communicated.
\item To efficiently parallelize the evaluation of mortar terms,
we will start from a baseline approach
where interfacial subdomains are aligned with the subdomains of the underlying bulk domain.
This method is straightforward to implement, preserves data locality,
and reduces communication between parallel processes.
However, it does not include all processes in the evaluation of the mortar terms and, thus, is not scalable.
We will then devise strategies for redistributing the interface domain decomposition
in order to increase parallel efficiency and scalability of the mortar evaluation.
\item As the contact configuration and area often changes over the course of a simulation,
we will propose a dynamic load balancing scheme.
Therefore, we will monitor characteristic quantities of the parallel evaluation of all mortar terms
and will trigger an adaptation of the interface domain decomposition
if the current state and computational behavior of the simulation indicates a deterioration of parallel performance.
\end{itemize}
We will discuss these approaches in detail and demonstrate their scaling behavior and applicability to large three-dimensional problems.
Although our current work studies scalable computational kernels for mortar methods in the context of classcial finite element analysis,
all findings are equally valid for isogeometric mortar methods ({\ie} NURBS-based interface discretizations).

The remainder of this paper is organized as follows:
After a brief description of the contact problem, its discretization, and suitable solution techniques in \secref{sec:Problem},
the implications of storing mortar discretizations on distributed memory machines will be discussed in \secref{sec:DistributedMemory}.
Domain decomposition approaches for an efficient evaluation of the mortar integrals will then be developed in \secref{sec:LoadBalancing}.
\secref{sec:NumEx} presents several numerical studies to assess communication patterns and demonstrate the parallel scalability of the proposed methods
in the context of computational contact mechanics,
before we conclude with some final remarks in \secref{sec:Conclusion}.

\section{Problem formulation and finite element discretization}
\label{sec:Problem}

While mortar methods are applicable to a broad spectrum of problems and partial differential equations (PDEs),
finite deformation contact problems are nowadays certainly one of the most appealing and challenging application areas for mortar methods in computational mechanics.
Hence, we focus on contact problems now, but keep the generality of mortar evaluations in mind.

\subsection{Governing equations}
\label{sec:GovEq}

In general, mortar methods allow for the coupling of several physical domains governed by PDEs
through enforcing coupling conditions at various coupling surfaces or interfaces.
Without loss of generality,
we focus our presentation on the two-body contact problem with bodies~$\indexedDom{\dom}{1}$ and~$\indexedDom{\dom}{2}$
which potentially come into frictionless contact along their contact boundaries~$\indexedDom{\CoupBound}{1}$ and~$\indexedDom{\CoupBound}{2}$, respectively.
Each subdomain~$\indexedDom{\dom}{\indDomain}, \indDomain\in\{1,2\}$
is governed by the initial boundary value problem of finite deformation elasto-dynamics,
reading
\begin{align*}
\matDiv \indexedDom{\teo\pki}{\indDomain} + \indexedDom{\matConfig{\bc{\teo\bodyforce}}}{\indDomain} & = \indexedDom{\matConfig{\density}}{\indDomain} \indexedDom{\ddot{\teo{\disp}}}{\indDomain}
& \text{ in } \indexedDom{\matConfig{\dom}}{\indDomain} \times [0, \ttimeend],\\
\indexedDom{\teo\disp}{\indDomain} & = \indexedDom{\bc{\teo\disp}}{\indDomain}
& \text{ on } \indexedDom{\DBound}{\indDomain} \times [0, \ttimeend],\\
\indexedDom{\teo\pki}{\indDomain} \indexedDom{\matConfig{\teo\normal}}{\indDomain} & = \indexedDom{\matConfig{\bc{\teo\trac}}}{\indDomain}
& \text{ on } \indexedDom{\NBound}{\indDomain} \times [0, \ttimeend],\\
\indexedDom{\teo\disp}{\indDomain} (\indexedDom{\teo\matCoord}{\indDomain}, 0) & = \indexedDom{\bc{\teo\disp}}{\indDomain} (\indexedDom{\teo\matCoord}{\indDomain})
& \text{ in } \indexedDom{\matConfig{\dom}}{\indDomain},\\
\indexedDom{\dot{\teo\disp}}{\indDomain} (\indexedDom{\teo\matCoord}{\indDomain}, 0) & = \indexedDom{\bc{\dot{\teo\disp}}}{\indDomain} (\indexedDom{\teo\matCoord}{\indDomain})
& \text{ in } \indexedDom{\matConfig{\dom}}{\indDomain}
\end{align*}
with the unknown displacement field~$\teo\disp$, the first Piola-Kirchhoff stress tensor~$\teo\pki$, the body force vector~$\matConfig{\bc{\teo\bodyforce}}$, density~$\matConfig{\density}$,
normal vector~$\matConfig{\teo\normal}$, and traction vector~$\matConfig{\bc{\teo\trac}}$ in the initial configuration~$\matConfig{\dom}$.
Furthermore, prescribed boundary and initial values are marked with~$\bc{\AnyQuantity}$.
First and second time derivatives are given as~$\dot{\AnyQuantity}$ and~$\ddot{\AnyQuantity}$, respectively.

For frictionless contact,
the contact constraints are typically given by the Hertz-Signorini-Moreau conditions, reading
\begin{align*}
\indexedNormal{\gap} \geq 0, \quad
\indexedNormal{\pressure} \leq 0, \quad
\indexedNormal{\pressure}\indexedNormal{\gap} = 0
\qquad
\text{ on }
\CoupBound \times [0, \ttimeend],
\end{align*}
with the contact pressure~$\indexedNormal{\pressure}$ along the contact interface~$\CoupBound$
and the gap function~$\indexedNormal{\gap}$ denoting the normal distance between the two bodies.
To later distinguish between the two sides of the contact interface,
we follow the traditional naming scheme and refer to~$\indexedDom{\CoupBound}{1}$
carrying the Lagrange multiplier as so-called ``{\slave}'' side~$\indexedSlave{\CoupBound}$,
while~$\indexedDom{\CoupBound}{2}$ denotes the ``{\master}'' side~$\indexedMaster{\CoupBound}$.

Since this paper is concerned with the efficient evaluation of the mortar terms on parallel computing clusters,
we will detail the discretization of all mortar-related terms in \secref{sec:Discretization}.
However, to keep the focus tight and concise,
we refer to the extensive literature for any further details on the finite element formulation and
discretization~\cite{Flemisch2007a,Lamichhane2005a,Lamichhane2007a,Popp2010a,Popp2012b,Wohlmuth2000a,Wohlmuth2012a},
the solution of the nonlinear problem via active set strategies~\cite{Hartmann2007a,Hintermueller2003a,Hueeber2005a,Hueeber2008a,Popp2009a},
as well as for details on the structure of the arising linear systems of equations and efficient
solvers thereof~\cite{Achdou1999a,Casarin1996a,Stefanica2001a,Wieners1999a,Wieners2003a,Wiesner2018a,Wiesner2021a,Wohlmuth2000b}.

\subsection{Discretization}
\label{sec:Discretization}

In order to perform the spatial discretization with FEM,
we assume the existence of a weak form of the contact mechanics problem summarized in \secref{sec:GovEq}.
For the additional terms arising in contact mechanics,
a Lagrange multiplier field~$\teo\lmb$ is introduced into the weak form to enforce the contact constraints,
leading to a mixed method with a variational inequality,
where both the primal field~$\teo\sol$ as well as the dual variable~$\teo\lmb$ need to be discretized in space.

For the sake of a concise presentation,
we skip the details of the FEM applied to the three-dimensional solid bodies~$\indexedDom{\matConfig{\dom}}{\indDomain}, \indDomain\in\{1,2\}$.
Considering the contact interface,
we adopt from the volume discretization the isoparametric concept with the parameter coordinate~$\mao\bxi = [\ParamCoord_1, \ParamCoord_2]$
and the shape functions~$\indexedNode{\ShapeFunc}{\indNodeSlave} \left(\mao\ParamCoordVec\right)$
defined at node~$\indNodeSlave$ of all $\nNodesSlave$ nodes on the discrete {\slave} surface~$\indexedSlave{\CoupBoundDisc}$
and~$\indexedNode{\ShapeFunc}{\indNodeMaster} \left(\mao\ParamCoordVec\right)$
defined at node~$\indNodeMaster$ of all ~$\nNodesMaster$ nodes on the discrete {\master} surface~$\indexedMaster{\CoupBoundDisc}$, respectively.
The interpolation of the displacement field on element level is then given as
\begin{align}
\label{eq:AnsatzDiscretization}
\indexedDom{\mao\sol}{1}\left(\mao\ParamCoordVec,\ttime\right)
= \sum_{\indNodeSlave=1}^{\nNodesSlave} \indexedNode{\ShapeFunc}{\indNodeSlave} \left(\mao\ParamCoordVec\right) \indexedNode{\mao\sol}{\indNodeSlave}\left(\ttime\right)
,\quad
\indexedDom{\mao\sol}{2}\left(\mao\ParamCoordVec,\ttime\right)
= \sum_{\indNodeMaster=1}^{\nNodesMaster} \indexedNode{\ShapeFunc}{\indNodeMaster} \left(\mao\ParamCoordVec\right) \indexedNode{\mao\sol}{\indNodeMaster}\left(\ttime\right).
\end{align}
As usual in mortar methods,
the Lagrange multiplier field~$\teo\lmb$ is discretized on $\nNodesSlaveWithLM$ nodes of the discrete slave surface~$\indexedSlave{\CoupBoundDisc}$,
reading
\begin{align}
\label{eq:AnsatzDiscretizationLagrangeMult}
\mao\lmb\left(\mao\ParamCoordVec,\ttime\right)
= \sum_{\indNode=1}^{\nNodesSlaveWithLM} \indexedNode{\ShapeFuncLM}{\indNode} \left(\mao\ParamCoordVec\right) \indexedNode{\mao\lmb}{\indNode}\left(\ttime\right),
\end{align}
where~$\indexedNode{\ShapeFuncLM}{\indNode} \left(\mao\ParamCoordVec\right)$ denotes
the Lagrange multiplier shape function at node~$\indNode$.
Thereby, either standard or dual shape functions can be used.

Inserting \eqref{eq:AnsatzDiscretization} and \eqref{eq:AnsatzDiscretizationLagrangeMult}
into the contact virtual work~$\WeakForm{\lm} = \int_{\CoupBound}\teo\lmb \, \left(\variation\indexedSlave{\teo\disp} - \variation\indexedMaster{\teo\disp}\right) \dx \Bound$
yields
\begin{align}
\label{eq:DiscreteWeakFormInterfaceFlux}
\begin{split}
\WeakForm{\lm} \approx \WeakForm{\lm,\Disc} =
\sum_{\indNode=1}^{\nNodesSlaveWithLM}\sum_{\indNodeSlave=1}^{\nNodesSlave}
\trans{\indexedNode{\mao\lmb}{\indNode}}
\overbrace{
\left[\int_{\indexedSlave{\CoupBoundDisc}}
\indexedNode{\ShapeFuncLM}{\indNode}
\indexedDom{\indexedNode{\ShapeFunc}{\indNodeSlave}}{1}
\dx\Bound\right]}^{\MoD\left[\indNode,\indNodeSlave\right]}
\variation\indexedDom{\indexedNode{\mao\sol}{\indNodeSlave}}{1}
- \sum_{\indNode=1}^{\nNodesSlaveWithLM}\sum_{\indNodeMaster=1}^{\nNodesMaster}
\trans{\indexedNode{\mao\lmb}{\indNode}}
\underbrace{
\left[\int_{\indexedSlave{\CoupBoundDisc}}
\indexedNode{\ShapeFuncLM}{\indNode}
\left(\indexedDom{\indexedNode{\ShapeFunc}{\indNodeMaster}}{2}\circ\disc{\MapMaToSl}\right)
\dx\Bound\right]}_{\MoM\left[\indNode,\indNodeMaster\right]}
\variation\indexedDom{\indexedNode{\mao\sol}{\indNodeMaster}}{2}.
\end{split}
\end{align}
The mortar matrices~$\MoD$ and~$\MoM$ associated with the {\slave} and {\master} side of the coupling interface are then assembled
from the nodal blocks~$\MoD\left[\indNode,\indNodeSlave\right]$ and~$\MoM\left[\indNode,\indNodeMaster\right]$
defined in~\eqref{eq:DiscreteWeakFormInterfaceFlux}, respectively.
In general, both~$\MoD$ and~$\MoM$ are rectangular matrices.
If $\nNodesSlaveWithLM = \nNodesSlave$
(which is common practice except for a few cases, {\eg} higher-order FEM~\cite{Lamichhane2005a,Lamichhane2007a,Popp2012b}),
$\MoD$ becomes square.
Furthermore, if~$\indexedNode{\ShapeFuncLM}{\indNode}$ are chosen as so-called dual shape functions
that satisfy a biorthogonality relationship with the standard shape functions~$\indexedNode{\ShapeFunc}{\indNode}$,
then $\MoD$ becomes a diagonal matrix and, thus, easy and computationally cheap to invert~\cite{Flemisch2007a,Lamichhane2005a,Lamichhane2007a,Popp2012b,Scott1990a,Wohlmuth2000a,Wohlmuth2001a,Wohlmuth2012a}.

We stress that both summands in~\eqref{eq:DiscreteWeakFormInterfaceFlux} contain integrals
over the {\slave} side~$\indexedDom{\CoupBoundDisc}{1}$ of the discrete coupling surface,
where the discretization is indicated by the additional subscript~$\sub{\AnyQuantity}{\Disc}$.
A suitable discrete mapping~${\disc{\MapMaToSl}:\indexedMaster{\CoupBoundDisc}\rightarrow\indexedSlave{\CoupBoundDisc}}$
from the {\master} side to the {\slave} side of the coupling interface is required,
because the discrete coupling surfaces~$\indexedMaster{\CoupBoundDisc}$ and~$\indexedSlave{\CoupBoundDisc}$ do not coincide anymore in general,
especially when considering {\nonmatching} meshes on curved interfaces.
These projections are usually based on a continuous field of normal vectors defined on the {\slave} side~$\indexedSlave{\CoupBoundDisc}$,
{\cf}~\cite{Popp2009a,Yang2005a}.

We note that the mortar matrices~$\MoD$ and~$\MoM$ also occur in the discrete representation of the Hertz-Signorini-Moreau conditions,
{\cf} \cite{Popp2010a} for example.

\subsection{Evaluation of mortar integrals}
\label{sec:EvaluationMortarIntegrals}

In general,
the evaluation of both~$\MoD\left[\indNode,\indNodeSlave\right]$ and~$\MoM\left[\indNode,\indNodeMaster\right]$ in~\eqref{eq:DiscreteWeakFormInterfaceFlux}
requires information from both the discrete {\slave} interface~$\indexedSlave{\CoupBoundDisc}$
and the discrete {\master} interface~$\indexedMaster{\CoupBoundDisc}$.
Firstly, this inevitably involves the discrete mapping~$\disc{\MapMaToSl}$
to project finite element nodes and quadrature points between {\slave} and {\master} sides.
In practice,
mortar integration is often performed on a piecewise flat geometrical approximation of the {\slave} surface~$\indexedSlave{\CoupBoundDisc}$
as proposed in~\cite{Puso2004c}.
For further details and an in-depth mathematical analysis, see~\cite{Dickopf2009a,Puso2004a,Puso2004b}.
Secondly, the {\slave}-sided integration domain~$\indexedSlave{\CoupBoundDisc}$ has to be split into so-called mortar segments,
such that both~$\indexedDom{\indexedNode{\ShapeFuncLM}{\indNode}}{1}$ and~$\indexedDom{\indexedNode{\ShapeFunc}{\indNodeMaster}}{2}$
are $C^1$-continuous on these segments,
as kinks in the function to be integrated would deteriorate the achievable accuracy of the numerical quadrature.
These mortar segments are arbitrarily shaped polygons,
which will then be decomposed into triangles to perform quadrature.
While the evaluation of~$\MoD\left[\indNode,\indNodeSlave\right]$ involves quantities
solely defined on the {\slave} interface~$\indexedSlave{\CoupBoundDisc}$,
the evaluation of~$\MoM\left[\indNode,\indNodeMaster\right]$ requires to integrate the product of
{\master} side shape functions~$\indexedDom{\indexedNode{\ShapeFunc}{\indNodeMaster}}{2}$
and {\slave} side shape functions~$\indexedDom{\indexedNode{\ShapeFuncLM}{\indNode}}{1}$
over the discrete {\slave} interface~$\indexedSlave{\CoupBoundDisc}$.

\Algref{alg:MortarSegmentation} outlines the necessary steps to perform segmentation
and numerical quadrature for one pair of {\slave} and {\master} elements.
\begin{algorithm}
\SetKwProg{myproc}{Procedure}{}{}
\SetKwFunction{proc}{MortarSegmentation}
\For{all {\slave} elements}{
\For{all associated {\master} elements in the vicinity of the current {\slave} element}{
Project {\master} elements onto {\slave} side\\
Find mesh intersection of {\slave} and {\master} elements via a clipping algorithm, see {\eg}~\cite{Foley1997a}.\\
Divide clip polygon into triangular integration cells.\\
Perform quadrature to compute entries of~$\MoD\left[\indNode,\indNodeSlave\right]$ and~$\MoM\left[\indNode,\indNodeMaster\right]$
according to~\eqref{eq:DiscreteWeakFormInterfaceFlux}.
}
}
\caption{Segment-based mortar integration for three-dimensional problems}
\label{alg:MortarSegmentation}
\end{algorithm}
We refer to \cite{Puso2004c} for a detailed description of all steps outlined in \algref{alg:MortarSegmentation}.
Although segment-based quadrature as described in \algref{alg:MortarSegmentation} undoubtedly delivers the highest achievable accuracy
for the numerical integration of~$\MoD\left[\indNode,\indNodeSlave\right]$ and~$\MoM\left[\indNode,\indNodeMaster\right]$ in three dimensions,
it comes at high computational expenses
related to mesh projection and intersection, subsequent triangulation as well as numerical quadrature.
In practice and also in the present work,
both mortar operators~$\MoD\left[\indNode,\indNodeSlave\right]$ and~$\MoM\left[\indNode,\indNodeMaster\right]$ are usually evaluated using segment-based integration
in order to guarantee conservation of linear momentum~\cite{Puso2004c}.
More efficient but possibly less accurate integration algorithms have been discussed in~\cite{Brivadis2015a,Farah2015a,Maday2002a,Wilking2017a}.

Having today's parallel computing architectures with distributed memory in mind,
the evaluation of~\eqref{eq:DiscreteWeakFormInterfaceFlux} brings along two major implications on the software and algorithm design:
\begin{enumerate}
\item The evaluation of the integrands in~\eqref{eq:DiscreteWeakFormInterfaceFlux}
requires information from both the {\slave} and {\master} side.
{\Slave} data is readlily available locally on each parallel process.
The implementation has to enable access also to {\master} side data,
that might be owned by another process or is stored on a different compute node.
\item The computational cost and time is mostly associated with numerical integration over the {\slave} side of the interface.
Parallelization can reduce the computational time by distributing the integration domain,
{\ie} the {\slave} interface, over multiple parallel processes.
\end{enumerate}
Therefore, we deduce the following requirements:
\begin{itemize}
\setlength\itemsep{-0.2em}
\item[{\reqMaster}:] Enable access to all required {\slave} and {\master} data during evaluation of mortar integrals
while keeping the memory demand and parallel communication low.
\item[{\reqSlave}:] Use parallel resources efficiently for numerical integration over the {\slave} side of the mortar interface,
also targeting parallel scalabity.
\end{itemize}
We will elaborate on these implications in \secsref{sec:DistributedMemory}{sec:LoadBalancing}
and outline various approaches to satisfy both requirements~{\reqMaster} and~{\reqSlave} in the context of parallel computing.

\section{Storing data of the contact interface on a parallel machine}
\label{sec:DistributedMemory}

When executing the FEM solver on a parallel machine,
data needs to be distributed among the different MPI ranks or compute nodes.
Now, we first summarize the basics of overlapping domain decomposition to distribute chunks of the discretization to individual processes.
Then, we discuss the implications on access to the relevant interface data during contact evaluation,
before we present and discuss several strategies to ensure access to the necessary data without excessive data redundancy.
Overall, this section is devoted to strategies in order to satisfy our basic requirement~{\reqMaster}.

\subsection{Overlapping domain decomposition}


We base our considerations on the existence of an FEM solver that can be executed on parallel computers
with a multitude of CPUs and/or compute nodes using a distributed memory architecture.
In our case, this FEM solver is our {\inhouse} code {\baci}~\cite{BaciURL}.
For optimal parallel treatment, the code base utilizes
\define{overlapping domain decomposition (DD)} techniques~\cite{Dolean2015a,Quarteroni2005a,Smith2008a,Toselli2005a}.
Using~$\nproc$ to denote the number of available parallel processes,
the computational domain~$\dom$ is divided into $\nproc$ subdomains~$\subdomainID{\dom}{\indSubdomain}$,
$\indSubdomain\in\{0,1,\hdots,\nsubdomain-1\}$.
A one-to-one mapping of subdomains to processes is employed, such that~$\nproc = \nsubdomain$.

An exemplary overlapping DD into four subdomains distributed to processes~$\indProc\in\{0,1,2,3\}$ is shown in \figref{fig:BasicsOverlappingDD}.
\begin{figure}
\begin{centering}


\begin{tikzpicture}[scale=0.8,node distance=0.2cm,font=\small]

\def\h{0.7} 
\def\hScaling{0.85} 
\def\len{4.2} 
\def\r{0.07} 

\def\twoH{1.4}
\def\threeH{2.1}
\def\fourH{2.8}
\def\fiveH{3.5}

\def\dist{0.75}
\def\myspace{1.0}

\begin{scope}[shift={(2*\myspace,0)}]

\draw [color=white,fill=red!20] (0,0) -- (3.5*\h,0) -- (3.5*\h,2.5*\h) -- (0,2.5*\h) -- cycle;
\draw [color=white,fill=green!20!black!10] (3.5*\h,0) -- (\len,0) -- (\len,3.5*\h) -- (3.5*\h,3.5*\h) -- cycle;
\draw [color=white,fill=blue!60!black!10] (3.5*\h,3.5*\h) -- (\len,3.5*\h) -- (\len,\len) -- (2.5*\h,\len) -- (2.5*\h,4.5*\h) -- (3.5*\h,4.5*\h) -- cycle;
\draw [color=white,fill=yellow!20] (0,2.5*\h) -- (3.5*\h,2.5*\h) -- (3.5*\h,4.5*\h) -- (2.5*\h,4.5*\h) -- (2.5*\h,\len) -- (0,\len) -- cycle;

\input{./fig/basics_overlapping_dd/mesh.tex}
    
\draw [very thick, dashed, color=orange] (3.5*\h,-0.5*\h) -- (3.5*\h,2.5*\h) -- (-0.5*\h,2.5*\h);
\draw [very thick, dashed, color=orange] (3.5*\h,2.5*\h) -- (3.5*\h,4.5*\h) -- (2.5*\h,4.5*\h) -- (2.5*\h,6.5*\h);
\draw [very thick, dashed, color=orange] (3.5*\h,3.5*\h) -- (6.5*\h,3.5*\h);
    
\node [above right] at (0,\len+0.5*\h) {Mesh and subdomains:};
\node [above left] at (0,0) {$\subdomainID{\dom}{0}$};
\node [above right] at (\len,0) {$\subdomainID{\dom}{1}$};
\node [below right] at (\len,\len) {$\subdomainID{\dom}{2}$};
\node [below left] at (0,\len) {$\subdomainID{\dom}{3}$};
\end{scope}

\begin{scope}[scale=1.0,shift={(0,-1.3*\len)}]
\draw [fill=black!30] (0,0) -- (4*\h,0) -- (4*\h,2*\h) -- (4*\h,3*\h) -- (0,3*\h) -- cycle;
\draw [pattern=north east lines, pattern color=blue] (0,2*\h) -- (3*\h,2*\h) -- (3*\h,0) -- (4*\h,0) -- (4*\h,3*\h) -- (0,3*\h) -- cycle;

\input{./fig/basics_overlapping_dd/mesh.tex}
\foreach \x in {0,\h,...,\threeH}
  \foreach \y in {0,\h,...,\twoH}
    \draw [fill=black] (\x,\y) circle (\r);
\foreach \y in {0,\h,...,\threeH}
  \draw [fill=white] (4*\h,\y) circle (\r);    
\foreach \x in {0,\h,...,\fourH}
  \draw [fill=white] (\x,3*\h) circle (\r);    
 
\node [above right] at (0,\len) {proc $0$:};

\end{scope}

\begin{scope}[scale=1.0, shift={(\len+1*\myspace,-1.3*\len)}]
\draw [fill=black!30] (4*\h,0) -- (\len,0) -- (\len,4*\h) -- (3*\h,4*\h) -- (3*\h,3*\h) -- (4*\h,3*\h) -- cycle;
\draw [pattern=north east lines, pattern color=blue] (3*\h,0) -- (4*\h,0) -- (4*\h,3*\h) -- (\len,3*\h) -- (\len,4*\h) -- (3*\h,4*\h) -- cycle;

\input{./fig/basics_overlapping_dd/mesh.tex}
\foreach \x in {\fourH,\fiveH,...,\len}
  \foreach \y in {0,\h,...,\threeH}
    \draw [fill=black] (\x,\y) circle (\r);
\foreach \y in {0,\h,...,\fourH}
  \draw [fill=white] (3*\h,\y) circle (\r);    
\foreach \x in {\threeH,\fourH,...,\len}
  \draw [fill=white] (\x,4*\h) circle (\r);   
  
\node [above right] at (0,\len) {proc $1$:};
\end{scope}

\begin{scope}[scale=1.0, shift={(2*\len+2*\myspace,-1.3*\len)}]
\draw [fill=black!30] (\len,4*\h) -- (\len,\len) -- (2*\h,\len) -- (2*\h,5*\h) -- (3*\h,5*\h) -- (3*\h,4*\h) -- cycle;
\draw [pattern=north east lines, pattern color=blue] (3*\h,3*\h) -- (\len,3*\h) -- (\len,4*\h) -- (4*\h,4*\h) -- (4*\h,5*\h) -- (3*\h,5*\h) -- (3*\h,\len) -- (2*\h,\len) -- (2*\h,4*\h) -- (3*\h,4*\h) -- cycle;

\input{./fig/basics_overlapping_dd/mesh.tex}
\foreach \x in {\fourH,\fiveH,...,\len}
  \draw [fill=black] (\x,4*\h) circle (\r);
\foreach \x in {\threeH,\fourH,...,\len}
  \foreach \y in {\fiveH,\len}
    \draw [fill=black] (\x,\y) circle (\r);    
\foreach \x in {\threeH,\fourH,...,\len}
  \draw [fill=white] (\x,3*\h) circle (\r);    
\foreach \x in {\twoH,\threeH}
  \draw [fill=white] (\x,4*\h) circle (\r);    
\foreach \y in {\fiveH,\len}
  \draw [fill=white] (2*\h,\y) circle (\r);   
  
\node [above right] at (0,\len) {proc $2$:};
\end{scope}

\begin{scope}[scale=1.0, shift={(3*\len+3*\myspace,-1.3*\len)}]
\draw [fill=black!30] (0,3*\h) -- (3*\h,3*\h) -- (3*\h,5*\h) -- (2*\h,5*\h) -- (2*\h,\len) -- (0,\len) -- cycle;
\draw [pattern=north east lines, pattern color=blue] (0,2*\h) -- (4*\h,2*\h) -- (4*\h,5*\h) -- (3*\h,5*\h) -- (3*\h,\len) -- (2*\h,\len) -- (2*\h,4*\h) -- (3*\h,4*\h) -- (3*\h,3*\h) -- (0,3*\h) -- cycle;

\input{./fig/basics_overlapping_dd/mesh.tex}
\foreach \x in {0,\h,...,\threeH}
  \foreach \y in {\threeH,\fourH}
    \draw [fill=black] (\x,\y) circle (\r);
\foreach \x in {0,\h,...,\twoH}
  \foreach \y in {\fiveH,\len}
    \draw [fill=black] (\x,\y) circle (\r);
\foreach \x in {0,\h,...,\fourH}
  \draw [fill=white] (\x,2*\h) circle (\r);    
\foreach \y in {\threeH,\fourH,...,\fiveH}
  \draw [fill=white] (4*\h,\y) circle (\r);    
\foreach \y in {\fiveH,\len}
  \draw [fill=white] (3*\h,\y) circle (\r);    
  
\node [above right] at (0,\len) {proc $3$:};
\end{scope}

\begin{scope}[scale=0.8,shift={(\len+6*\myspace,\len)}]

\draw [very thick, dashed, color=orange] (-0.5*\h,0) -- (0.5*\h,0);
\node [right=\h] at (0,0) {Subdomain boundaries};

\node at (0,-\dist) {$\subdomainID{\dom}{\indSubdomain}$};
\node [right=\h] at (0,-\dist) {Subdomain~$\indSubdomain$ on process~$\indProc=\indSubdomain$};

\draw [fill=black] (0,-2*\dist) circle (\r);
\node [right=\h] at (0,-2*\dist) {Owned nodes};

\draw [fill=white] (0,-3*\dist) circle (\r);
\node [right=\h] at (0,-3*\dist) {Ghosted nodes};

\draw [fill=black!30] (-1/2*\hScaling*\h,-4*\dist-1/2*\hScaling*\h) rectangle (1/2*\hScaling*\h,-4*\dist+1/2*\hScaling*\h);
\node [right=\h] at (0,-4*\dist) {Owned elements};

\draw [pattern=north east lines, pattern color=blue] (-1/2*\hScaling*\h,-5*\dist-1/2*\hScaling*\h) rectangle (1/2*\hScaling*\h,-5*\dist+1/2*\hScaling*\h);
\node [right=\h] at (0,-5*\dist) {Elements integrated by multiple processes};

\end{scope}
  
\end{tikzpicture}
\caption{Exemplary overlapping domain decomposition and parallel assembly
involving four subdomains~$\subdomainID{\dom}{\indSubdomain}, \indSubdomain\in\{0,1,2,3\}$
assigned to four parallel processes~$\indProc\in\{0,1,2,3\}$.
Since each process can only assembly into unknowns of owned nodes,
elements spanning across the subdomain boundaries need to be evaluated by multiple processes.
This requires ghosting of nodes and elements,
which entails parallel communication among multiple processes.
}
\label{fig:BasicsOverlappingDD}
\end{centering}
\end{figure}
While each node in the finite element discretization is uniquely assigned to a subdomain~$\subdomainID{\dom}{\indSubdomain}$,
elements might span subdomain boundaries.
We stress that processes can only access data of nodes that they own themselves.
This has implications on finite element evaluation and assembly:
A process~$\indProc$ can only assemble into those entries of the global residual vector and those rows of the global Jacobian matrix
that are associated with nodes in~$\procID{\dom}{\indProc}$.
Hence, elements that span across subdomain boundaries will be evaluated by all processes
that own at least one of this element's nodes
such that each process can assemble quantities associated with its own nodes.\footnote{As an alternative, linear algebra data structures, that are specialized for FEM computations, are available. They allow to assemble into off-process rows. Naturally, communication among parallel processes is required.}
This requires communication of data prior to the evaluation,
{\ie} data of off-process nodes needs to be communicated.
This is often referred to as \define{ghosting}.
Ideally, subdomains exhibit a small surface-to-volume ratio to minimize the amount of data subject to ghosting.

In our code base {\baci}, we employ the hypergraph partitioning package \SoftwarePackage{Zoltan}~\cite{Boman2012a}
with the \SoftwarePackage{ParMETIS} backend
to decompose the computational domain~$\dom$ into $\nproc$ subdomains~$\subdomainID{\dom}{\indSubdomain}$.
Parallel data structures and parallel linear algebra is enabled through the {\trilinos}\footnote{https://trilinos.github.io} packages
\SoftwarePackage{Epetra}, \SoftwarePackage{Tpetra}, and \SoftwarePackage{Xpetra}.
Iterative solvers for sparse systems of linear equations are taken
form the {\trilinos} packages \SoftwarePackage{AztecOO}~\cite{Heroux2007a} and \SoftwarePackage{Belos}~\cite{Bavier2012a}
with scalable {\multilevel} preconditioners from \SoftwarePackage{ML}~\cite{gee2006a}
and \SoftwarePackage{MueLu}~\cite{BergerVergiat2019a}.

\subsection{Implications of distributed memory on the contact search and evaluation}

\begin{figure}
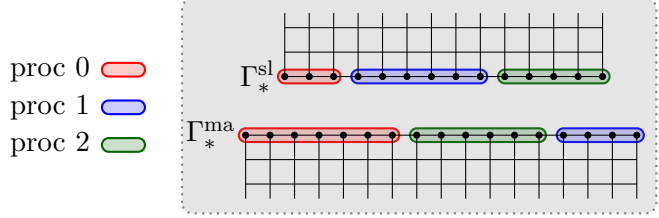

\begin{center}

\begin{tikzpicture}

\input{./fig/interface_storage/common_defines.tex}

\begin{scope}[scale=1.3,rotate=0] 

  \begin{scope}[xshift=0,yshift=0]    
  
    \input{./fig/interface_storage/slave_subdomains.tex}
    
    \input{./fig/interface_storage/master_defines.tex}
    \draw[thick,\colorProcOne,fill=\colorProcOne,fill opacity=0.2] \convexpath{mA,mB}{\padding};
    \draw[thick,\colorProcThree,fill=\colorProcThree,fill opacity=0.2] \convexpath{mC,mD}{\padding};
    \draw[thick,\colorProcTwo,fill=\colorProcTwo,fill opacity=0.2] \convexpath{mE,mF}{\padding};
  
    \input{./fig/interface_storage/meshes.tex}
    
    
  \end{scope} 
  
  \begin{scope}[xshift=-2.0cm,yshift=-0.1cm]
    \input{./fig/interface_storage/par_dist_color_legend.tex}
  \end{scope} 

\end{scope} 

\end{tikzpicture}
\caption{Without particular measures, DDs of {\master} and {\slave} side of the interface distribute each interface side to some processes.
Geometrically close portions of the {\master} and {\slave} interface are not guaranteed to reside on the same process.
Without further measures, each process~$\indProc$ can only identify possibly contacting pairs of {\slave}/{\master} elements
from the subset of {\master} elements owned by process~$\indProc$, {\ie} {\master} elements that reside in the set~$\subdomainID{\dom}{\indProc}\cap\indexedMaster{\CoupBound}$.
This can and needs to be alleviated by extending the ghosting of the {\master} side of the interface.
(For simplicity of visualization, coloring of ownership omits ghosted elements stemming from the overlapping interface DD.)}
\label{fig:InterfaceDDBasicLayout}
\end{center}
\end{figure}

Without loss of generality and for ease of presentation,
we assume that the entire discretization of a two-body contact problem has undergone an overlapping DD
and that each subdomain~$\indSubdomain\in\{0, \hdots, \nsubdomain - 1\}$ has been assigned to a process~$\indProc\in\{0, \hdots, \nproc - 1\}$.
For the purpose of illustration, we will discuss the case of $\nproc = 3$ subdomains and further assume
that every process owns a part of the {\master} and of the {\slave} interface as illustrated in \figref{fig:InterfaceDDBasicLayout}.
Please note that our considerations also hold,
if some processes only own a part of either the {\slave} or the {\master} side of the interface discretization
or even if some processes do not own any portion of the contact interface at all.

When process~$\indProc$ is performing contact search and evaluation on its share of the {\slave} interface,
it needs access to data from the geometrically close {\master} side of the interface.
In a parallel computing environment,
the required data from the {\master} side of the interface does not necessarily reside on that same process~$\indProc$.
Still, access has to be enabled in order to
\begin{itemize}
\setlength\itemsep{-0.2em}
\item identify pairs of {\slave}/{\master} elements, that potentially are in active contact.
This step is usually referred to as ``contact search''.
\item  evaluate the second integrand in~\eqref{eq:DiscreteWeakFormInterfaceFlux},
where the shape functions~$\indexedDom{\indexedNode{\ShapeFunc}{\indNodeMaster}}{2}$
defined on the {\master} side need to be evaluated and projected onto the {\slave} side.
\end{itemize}
If the required data of the {\master} side resides on a different parallel process~$\altIndProc$ than the current {\slave}-sided process~$\indProc$,
this data has to be communicated or ``ghosted'' ({\cf} \figref{fig:BasicsOverlappingDD})
from process~$\altIndProc$ to process~$\indProc$ in order to be known by process~$\indProc$.
Therefore, the ghosting of the {\master} interface discretization has to be extended.
Since such an extension will impact the inter-processes communication demand as well as the on-process memory demand,
we will introduce models for communication and memory demands in \secref{sec:ModelingStorageDemand}.
More importantly, we will discuss various approaches for extending the ghosting of the {\master} interface discretization
in \secsref{sec:MasterRedundantStorage}{sec:MasterDistributedStorage},
where we will also discuss the impact of these ghosting extension strategies on the memory demand.

\subsection{Models for communication and memory demand}
\label{sec:ModelingStorageDemand}

Starting from an overlapping DD and distributed storage of both interface discretizations,
data needs to be communicated among processes to facilitate the mortar evaluation.
We will use~$\costCommunication$ to denote the amount of data to be sent over the interconnect of all compute nodes and processes.
Since data related to the {\slave} side of the interface discretization just remains on its process~$\indProc$,
$\procID{\indexedSlave{\costCommunication}}{\indProc} = 0$.
As has already been indicated in \figref{fig:InterfaceDDBasicLayout},
process~$\indProc$ owning the portion~$\subdomainID{\indexedSlave{\Bound}}{\indSubdomain}$ of the {\slave} side of the mortar interface
requires the {\master} side's data from those processes owning the geometrically close {\master} elements.
Hence, usually $\procID{\indexedMaster{\costCommunication}}{\indProc} > 0$,
especially if a situation as depicted in \figref{fig:InterfaceDDBasicLayout} occurs.
Although an explicit expression to compute~$\procID{\indexedMaster{\costCommunication}}{\indProc}$ cannot be given,
as it highly depends on the software implementation at hand,
it for sure is related to the number of nodes~$\nnode$ and elements~$\nEle$ to be communicated.
We denote this relation by
\begin{align}
\label{eq:CostCommunicationRelationToMasterDiscretization}
\procID{\indexedMaster{\costCommunication}}{\indProc} \propto \costCommunicationMeasure(\nnode, \nEle)
\end{align}
with~$\costCommunicationMeasure(\nnode, \nEle)$ referring to an implementation-specific measure describing the cost of parallel communication.
The total amount of data to be communicated to process~$\indProc$ sums up to
\begin{align}
\label{eq:ModelingCommunication}
\costCommunication = \sum_{\indProc}^{\nproc - 1} \procID{\indexedMaster{\costCommunication}}{\indProc}.
\end{align}
Obviously, $\costCommunication$ increases with an increasing number of subdomains.
More importantly, however, it is impacted by the individual contributions~$\procID{\indexedMaster{\costCommunication}}{\indProc}$.
Especially when the number of subdomains, that are required to solve a given problem, is fixed,
reducing $\procID{\indexedMaster{\costCommunication}}{\indProc}$ is key to reduce the overall cost of communication.
Naturally, $\costCommunication = 0$ if $\nproc = 1$.

From the domain decomposition of the underlying bulk field,
the memory demand~$\procID{\bulk{\costStorage}}{\indProc}$ per process~$\indProc$ is given.
For the mortar interface discretizations,
we use~$\procID{\indexedSlave{\costStorage}}{\indProc}$ to denote the memory demand of the {\slave} interface portion~$\subdomainID{\indexedSlave{\Bound}}{\indSubdomain}$ on process~$\indProc$.
Furthermore, $\procID{\indexedMaster{\costStorage}}{\indProc}$ refers to the memory demand of the {\master} interface portion~$\subdomainID{\indexedMaster{\Bound}}{\indSubdomain}$ on process~$\indProc$.
Then, the total memory demand~$\procID{\costStorage}{\indProc}$ on process~$\indProc$ is given as
\begin{align}
\label{eq:ModelingMemoryDemand}
\procID{\costStorage}{\indProc} = \procID{\bulk{\costStorage}}{\indProc} + \procID{\indexedSlave{\costStorage}}{\indProc} + \procID{\indexedMaster{\costStorage}}{\indProc}
\quad
\forall \indProc\in\left\{0,\hdots,\nproc-1\right\}.
\end{align}
Note that~$\procID{\costStorage}{\indProc}$ includes the amount of memory required for owned nodes/elements as well as for ghost nodes/elements
originating from the overlapping DD with an element overlap of~$1$.
We stress that~$\procID{\bulk{\costStorage}}{\indProc}$ is fully determined by the overlapping DD of the underlying bulk fields
and that~$\procID{\indexedSlave{\costStorage}}{\indProc}$ is only governed by the overlapping DD of the {\slave} interface discretization,
that might arise from any of the schemes proposed in \secref{sec:LoadBalancing} later.
At this point, only the {\master} interface's contribution~$\procID{\indexedMaster{\costStorage}}{\indProc}$ can be controlled
by choosing a specific ghosting extension strategy.

%
%

\subsection{Redundant storage: the straightforward case}
\label{sec:MasterRedundantStorage}

\begin{figure}
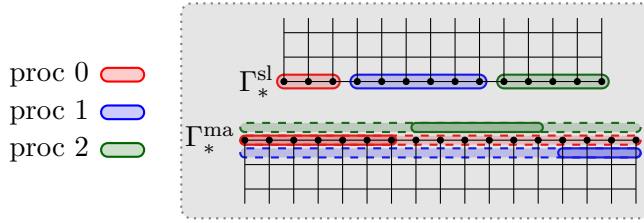

\begin{center}

\begin{tikzpicture}

\input{./fig/interface_storage/common_defines.tex}

\begin{scope}[scale=1.3,rotate=0] 

  \input{./fig/interface_storage/slave_subdomains.tex}
  
  \node (mLeftProcOne) at (0.0,0.0) {};
  \node (mMidProcOne) at (6*\h,0.0) {};
  \node (mRightProcOne) at (\masterLength,0.0) {};
  \draw[thick,\colorProcOne,fill=\colorProcOne,fill opacity=0.2] \convexpath{mLeftProcOne,mMidProcOne}{2/3*\padding};
  \draw[thick,dashed,\colorProcOne,fill=\colorProcOne,fill opacity=0.2] \convexpath{mLeftProcOne,mRightProcOne}{2/3*\padding};

  \node (mLeftProcTwo) at (0.0,-\verticalShift) {};
  \node (mMidProcTwo) at (13*\h,-\verticalShift) {};
  \node (mRightProcTwo) at (\masterLength,-\verticalShift) {};
  \draw[thick,\colorProcTwo,fill=\colorProcTwo,fill opacity=0.2] \convexpath{mMidProcTwo,mRightProcTwo}{2/3*\padding};
  \draw[thick,dashed,\colorProcTwo,fill=\colorProcTwo,fill opacity=0.2] \convexpath{mLeftProcTwo,mRightProcTwo}{2/3*\padding};
  
  \node (mLeftProcThree) at (0.0,\verticalShift) {}; 
  \node (mMidLeftProcThree) at (7*\h,\verticalShift) {};
  \node (mMidRightProcThree) at (12*\h,\verticalShift) {};
  \node (mRightProcThree) at (\masterLength,\verticalShift) {};
  \draw[thick,\colorProcThree,fill=\colorProcThree,fill opacity=0.2] \convexpath{mMidLeftProcThree,mMidRightProcThree}{2/3*\padding};
  \draw[thick,dashed,\colorProcThree,fill=\colorProcThree,fill opacity=0.2] \convexpath{mLeftProcThree,mRightProcThree}{2/3*\padding};

  
  \input{./fig/interface_storage/meshes.tex}
  

  \begin{scope}[xshift=-2.0cm,yshift=-0.1cm]
    \input{./fig/interface_storage/par_dist_color_legend.tex}
  \end{scope} 
  
\end{scope} 
\end{tikzpicture}
\caption{Fully redundant storage of the {\master} interface discretization: Solid lines indicate data that is owned by a particular process due to the initial DD.
Dashed lines indicate data that is available through the extended ghosting.
With fully redundant storage of the {\master} discretization on each process,
each process~$\indProc$ can immediately identify all pairs of {\slave}/{\master} elements, that are possibly in active contact.}
\label{fig:MasterRedundantStorage}
\end{center}
\end{figure}

The probably most straightforward remedy for the issue of undetected {\master}/{\slave} pairs described in \figref{fig:InterfaceDDBasicLayout}
is to fully extend the {\master} side's ghosting to all processes,
{\ie} to store the entire {\master} side of the interface redundantly on every process~$\indProc$.
This scenario of distributed storage of the {\slave} interface discretization, but redundant storage of the {\master} interface discretization
is illustrated in \figref{fig:MasterRedundantStorage} for an exemplary number of three processes.
The {\slave} interface~$\indexedSlave{\CoupBound}$ is decomposed into three subdomains and distributed to the processes \procName{0}, \procName{1}, and \procName{2},
indicated by coloring.
The {\master} interface~$\indexedMaster{\CoupBound}$ starts out from its initial DD (colored boxes with solid lines) as already seen in \figref{fig:InterfaceDDBasicLayout}.
Then, its ghosting is extended over the entire {\master} interface~$\indexedMaster{\CoupBound}$ (colored boxes with dashed lines),
such that~$\indexedMaster{\CoupBound}$ is now stored redundantly on all three processes.

The redundant storage of the {\master} side of the interface just requires a one-time setup and communication cost at the beginning of the simulation
in order to extend the ghosting of {\master} data to the entire {\master} interface,
but then enables access to every bit of {\master} interface data from every process~$\indProc\in\{0,1,\hdots,\nproc-1\}$
without further communication among parallel processes.
After the ghosting has been extended following the idea of fully redundant storage,
all algorithmic steps, {\eg} the contact search or the evaluation of~\eqref{eq:DiscreteWeakFormInterfaceFlux},
can be performed immediately without further communication.

In terms of the communication cost~$\costCommunication$,
however, this approach is rather expensive:
since the entire {\master} discretization needs to be communicated to every {\slave} processor~$\indProc$,
the total communication cost can be estimated via~\eqref{eq:CostCommunicationRelationToMasterDiscretization}
and~\eqref{eq:ModelingCommunication} as
\begin{align}
\label{eq:ModelCostCommunicationRedundantMaster}
\costCommunication \approx \nproc \costCommunicationMeasure(\nNodesMaster, \nEleMaster),
\end{align}
where \emph{all} nodes and elements of the {\master} side of the interface discretization enter the cost estimate.
The model~\eqref{eq:ModelCostCommunicationRedundantMaster} suffers only from a slight over-estimation,
since a part of the {\master} surface might already be located on the target process and, thus, does not need to be communicated.
Yet, this over-estimation becomes smaller for an increasing number of subdomains.

Since the entire {\master} discretization has to be stored on each process along with a portion of the {\slave} discretization,
the memory demand of this approach can grow quite excessively when going to large {\master} interface discretizations.
The maximum problem size, for which this strategy still works, cannot be given theoretically.
It strongly depends on several key factors, for example the exact specifications of the computing hardware or
intricate details of the software implementation.
Considering the memory model~\eqref{eq:ModelingMemoryDemand},
the per-process {\master} contribution~$\procID{\indexedMaster{\costStorage}}{\indProc}$
has to be replaced by the memory consumption~$\indexedMaster{\costStorage}$ of the entire {\master} interface
since each process stores the entire {\master} discretization.
Since~$\indexedMaster{\costStorage}$ grows with mesh refinement,
the total storage demand~$\procID{\costStorage}{\indProc}$ on process~$\indProc$ is not bounded.
This limits the applicability of redundant storage to small and medium sized interface discretizations,
depending on the hardware at hand.

Besides the possibly unbounded memory demand,
fully redundant ghosting of the {\master} side also comes with a {\runtime} cost:
when process~$\indProc$ loops over all of its nodes/elements of the {\master} discretization,
then it actually loops over \emph{all} nodes/elements of the entire {\master} discretization,
although most of the nodes/elements are irrelevant on process~$\indProc$
as they are not located in the geometric vicinity of process~$\indProc$'s {\slave} nodes/elements.
Naturally, the code is not aware of any concept of vicinity prior to the contact search,
so this cost cannot be avoided with this approach.



\subsection{Distributed storage: going to large problems}
\label{sec:MasterDistributedStorage}\label{sec:Binning}

As soon as the memory demand~$\procID{\costStorage}{\indProc}$ exceeds the available memory on a computing node,
redundant storage as described in \secref{sec:MasterRedundantStorage} should not be applied anymore
to avoid performance degradation due to memory swapping.
Following~\eqref{eq:ModelingMemoryDemand}, the total storage demand~$\procID{\costStorage}{\indProc}$ per process
can be reduced by reducing the storage demand of the {\master} interface.
In particular, when storing also the {\master} interface discretization in a distributed fashion,
its storage demand per process can be reduced to~$\procID{\indexedMaster{\costStorage}}{\indProc} < \indexedMaster{\costStorage}$ for $\indProc\in\left\{0,\hdots,\nproc-1\right\}, \nproc \geq 2$.
Similarly, when the growth in {\runtime} for loops over {\master} nodes/elements becomes prohibitive,
reducing the portion of the {\master} interface stored on each process~$\indProc$ is expected to speed up simulations.
Still, each portion of the {\slave} interface needs to have access to those parts of the {\master} interface
that reside in its geometric proximity ({\cf} \figref{fig:InterfaceDDBasicLayout}).
In turn, measuring geometric proximity requires access to all pairs of {\slave} and {\master} nodes.

This situation can be remedied by different algorithmic modifications:
Within a token-based evaluation strategy, {\eg} inspired by Round-Robin (RR) scheduling~\cite{Brucker2007a},
the parallel decomposition and distribution of the {\slave} interface is fixed.
On the {\master} side, just the decomposition into subdomains is fixed,
while the subdomain-to-process mapping is shifted by one process per RR iteration
until every process has owned each {\master} interface subdomain once.
Since an RR loop requires $\nproc$ iterations for a complete evaluation of all {\slave} elements,
its {\runtime} cost is high and has even proven to be prohibitive in large-scale applications,
which we have also observed in our own experiments.

As an alternative,
the incorporation of the notion of \emph{proximity} already into the extension of the {\master} side's ghosting offers a promising solution.
Hence, we resort to pre-computing ghosting data based on a \emph{geometrically motivated binning approach},
where we exploit the fact that the contact search needs to identify all {\master} elements in the proximity of a given {\slave} element.
This idea is inspired by \cite{Plimpton1995a},
where a similar parallel algorithm is used for the spatial decomposition of atoms in short-range molecular dynamics simulations.

\begin{figure}
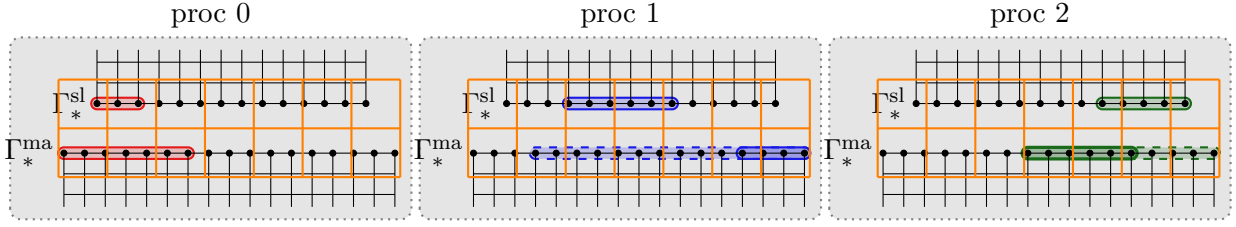

\begin{center}

\begin{tikzpicture}

\input{./fig/interface_storage/common_defines.tex}

\def\colorBins{orange} 
\def\binSize{0.59cm} 
\def\xshiftBins{0.0cm} 
\def\numXBins{7} 

\begin{scope}[scale=1.1] 
\begin{scope}[scale=1.0,rotate=0,xshift=-\boxWidth-\dist] 

  \node [above, align=center] at (\boxLeft+\boxWidth/2,\boxBottom+\boxHeight) {proc 0};

  \node (sA) at (\xOffset,\gap) {};
  \node (sB) at (\xOffset+0.5cm,\gap) {};
  \draw[thick,\colorProcOne,fill=\colorProcOne,fill opacity=0.2] \convexpath{sA,sB}{\padding};

  \node (mOwnedLeftProcOne) at (0.0,0.0) {};
  \node (mOwnedRightProcOne) at (6*\h,0.0) {};
  \draw[thick,\colorProcOne,fill=\colorProcOne,fill opacity=0.2] \convexpath{mOwnedLeftProcOne,mOwnedRightProcOne}{\padding};

  \input{./fig/interface_storage/meshes.tex}

  \input{./fig/interface_storage/bins.tex}


\end{scope} 

\begin{scope}[scale=1.0,rotate=0] 

  \node [above, align=center] at (\boxLeft+\boxWidth/2,\boxBottom+\boxHeight) {proc 1};

  \node (sC) at (\xOffset+0.75cm,\gap) {};
  \node (sD) at (\xOffset+2.0cm,\gap) {};
  \draw[thick,\colorProcTwo,fill=\colorProcTwo,fill opacity=0.2] \convexpath{sC,sD}{\padding};

  \node (mOwnedLeftProcTwo) at (13*\h,0) {};
  \node (mOwnedRightProcTwo) at (16*\h,0) {};
  \draw[thick,\colorProcTwo,fill=\colorProcTwo,fill opacity=0.2] \convexpath{mOwnedLeftProcTwo,mOwnedRightProcTwo}{\padding};
  \node (mGhostedLeftProcTwo) at (3*\h,0) {};
  \node (mGhostedRightProcTwo) at (16*\h,0) {}; 
  \draw[thick,dashed,\colorProcTwo,fill=\colorProcTwo,fill opacity=0.2] \convexpath{mGhostedLeftProcTwo,mGhostedRightProcTwo}{\padding};

  \input{./fig/interface_storage/meshes.tex}

  \input{./fig/interface_storage/bins.tex}

\end{scope} 

\begin{scope}[scale=1.0,rotate=0,xshift=\boxWidth+\dist] 

  \node [above, align=center] at (\boxLeft+\boxWidth/2,\boxBottom+\boxHeight) {proc 2};

  \node (sE) at (\xOffset+2.25cm,\gap) {};
  \node (sF) at (\xOffset+3.25cm,\gap) {};
  \draw[thick,\colorProcThree,fill=\colorProcThree,fill opacity=0.2] \convexpath{sE,sF}{\padding};

  \node (mOwnedLeftProcThree) at (7*\h,0) {};
  \node (mOwnedProcThree) at (12*\h,0) {};
  \draw[very thick,\colorProcThree,fill=\colorProcThree,fill opacity=0.2] \convexpath{mOwnedLeftProcThree,mOwnedProcThree}{\padding};
  \node (mGhostedLeftProcThree) at (7*\h,0) {}; 
  \node (mGhostedProcThree) at (16*\h,0) {};
  \draw[thick,dashed,\colorProcThree,fill=\colorProcThree,fill opacity=0.2] \convexpath{mGhostedLeftProcThree,mGhostedProcThree}{\padding};

  \input{./fig/interface_storage/meshes.tex}

  \input{./fig/interface_storage/bins.tex}


\end{scope} 
\end{scope} 

\end{tikzpicture}
\caption{Extended ghosting of the {\master} interface using a binning scheme ---
We exemplarily show three parallel processes and depict each one in its own sketch for the sake of presentation.
Bins are sketched in solid orange lines.
On~$\indexedSlave{\CoupBound}$, mesh entities (such as nodes and elements) are owned by the respective process anyway.
On~$\indexedMaster{\CoupBound}$, solid lines indicate data that is owned by this process,
while dashed lines indicate data that has been ghosted via binning.}
\label{fig:InterfaceDDBinning}
\end{center}
\end{figure}

In the context of mortar methods,
we will first construct an axis aligned bounding box around the mortar interface,
{\ie} a cuboid box that is oriented along the Cartesian axes and encloses all nodes of the mortar interface.
Then, this bounding box will be covered with a set of Cartesian bins
that are independent of the finite element meshes of the contacting bodies ({\cf} \figref{fig:InterfaceDDBinning}).
Since the contacting bodies are moving relative to the background bins,
{\slave} nodes or elements can migrate between individual bins over time.
In order to not loose track of individual nodes or elements due to this motion,
the minimal bin size~$\minBinSize$ is chosen as
\begin{align*}
\minBinSize = \max_{\nEleSlave} \indexedSlave{\hEle} + 2\cdot\Dt\cdot\indexedCoupling{\average{\dot{\mao\sol}}}
\end{align*}
with~$\max_{\nEleSlave} \indexedSlave{\hEle}$ being the largest element edge of the {\slave} discretization,
$\Dt$ representing the time step size,
$\indexedCoupling{\dot{\mao\sol}}$ denoting the vector of nodal interface velocities
and~$\average{\AnyQuantity}$ referring to the mean value of~$\AnyQuantity$, respectively.
If the interface velocity is not available in static problems,
it can be replaced via a finite difference approximation {\wrt} to the previous load step.
Analogously, the axis aligned bounding box embracing all mortar nodes
is expanded by~$\minBinSize$ in each direction.
Then, the actual bin size~$\binSize$ and number of bins per direction is computed
based on the dimensions of the expanded axis aligned bounding box and the minimal bin size~$\minBinSize$.
We then apply \algref{alg:Binning} to compute process-specific lists~$\procID{\indexedMaster{\listGhostEles}}{\indProc}$ of {\master} elements to be ghosted for each process~$\indProc$.
\begin{algorithm}
\SetKwProg{myproc}{Procedure}{}{}
\SetKwFunction{proc}{Binning}
Sort all $\nEleSlave$ {\slave} elements into bins\\
Sort all $\nEleMaster$ {\master} elements into bins\\
\For{each process~$\indProc$}{
Find set of bins~$\listOfBinsEnclosingSubdomain{\indSubdomain}$ enclosing $\indProc$'s {\slave} subdomain~$\subdomainID{\indexedSlave{\Bound}}{\indSubdomain}$\\
\For{each bin~$\indBin$ in~$\listOfBinsEnclosingSubdomain{\indSubdomain}$}{
Find set of bins~$\listOfNeighboringBins{\indBin}$ neighboring the current bin~$\indBin$\\
Collect all {\master} elements of~$\listOfNeighboringBins{\indBin}$ in~$\procID{\indexedMaster{\listGhostEles}}{\indProc}$
}
}
\caption{Geometrically motivated binning to pre-compute ghosting of the {\master} interface}
\label{alg:Binning}
\end{algorithm}

\Figref{fig:InterfaceDDBinning} illustrates the binning approach detailed in \algref{alg:Binning} for three processes.
For \procName{0}, no further ghosting is required in this example,
since all required {\master} elements already reside in the neighboring bins of the set of bins~$\listOfBinsEnclosingSubdomain{0}$
enclosing all {\slave} elements of~$\indexedSlave{\CoupBoundSubdomain{0}}$.
In contrast, the {\master} elements owned by \procName{1} are not contained in~$\listOfBinsEnclosingSubdomain{1}$
and do not participate to the evaluation of mortar terms in~$\indexedSlave{\CoupBoundSubdomain{1}}$.
The {\master} elements of interest, {\ie} those in the neighboring bins of~$\listOfBinsEnclosingSubdomain{1}$, need to be ghosted,
which leaves out the {\master} elements in the left most bin covering~$\indexedMaster{\CoupBound}$.
Finally, \procName{2} already owns some of the required elements and only needs to ghost some additional elements.

The communication cost~$\procID{\indexedMaster{\costCommunication}}{\indProc}$ for each processor~$\indProc$
now depends on the number of nodes/elements in the current bin~$\indBin$ and its neighboring bins.
Due to the Cartesian character of bins,
each bin has 8 or 26 neighbors in 2D or 3D, respectively.
Based on a constant bin size and assuming uniform mesh sizes,
the cost measure~$\costCommunicationMeasure$ per subdomain introduced in~\eqref{eq:CostCommunicationRelationToMasterDiscretization}
is now evaluated with $8\times\procID{\nNodesMaster}{\indProc}$ or~$26\times\procID{\nNodesMaster}{\indProc}$ nodes
and $8\times\procID{\nEleMaster}{\indProc}$ or~$26\times\procID{\nEleMaster}{\indProc}$ elements
for 2D and 3D problems, respectively.
With an increasing number of subdomains and under the assumption of uniform meshes,
the total cost for communication is then bounded by
\begin{align}
\label{eq:CommunicationCostBinning}
\costCommunication \leq
\begin{cases}
8\cdot\procID{\indexedMaster{\costCommunication}}{\indProc} & \text{for 2D}\\
26\cdot\procID{\indexedMaster{\costCommunication}}{\indProc} & \text{for 3D}
\end{cases}
\end{align}
which is a significant reduction for large core counts compared to~\eqref{eq:ModelCostCommunicationRedundantMaster}.
The scalar factors in~\eqref{eq:CommunicationCostBinning} originate from the number of neighboring bins in 2D and 3D, respectively.

Regarding memory demand as estimated via~\eqref{eq:ModelingMemoryDemand},
the {\master} side's demand~$\procID{\indexedMaster{\costStorage}}{\indProc}$
now comprises of all {\master} elements stored on process~$\indProc$
plus all {\master} elements in neighboring bins.
Assuming bin sizes similar to the size of subdomains~$\subdomainID{\dom}{\indSubdomain}$ as well as evenly sized {\master} elements,
the {\master} side's storage demand is bounded by~$5\times\procID{\indexedMaster{\costStorage}}{\indProc}$
or~$9\times\procID{\indexedMaster{\costStorage}}{\indProc}$ for 2D and 3D problems, respectively.
While the number of bins and, thus, the effort to sort {\master} elements into bins increases with a smaller characteristic bin size~$\binSize$,
the storage demands for each process~$\indProc$ diminishes even more.

\subsection{Intermediate discussion of ghosting strategies}

So far, we have concerned ourselves with strategies to satisfy the requirement~{\reqMaster}.
Before addressing~{\reqSlave} in \secref{sec:LoadBalancing},
we briefly discuss some properties of the presented strategies for the ghosting of the {\master} interface.

While the fully redundant ghosting presented in \secref{sec:MasterRedundantStorage} appears as straightforward,
easy to implement, and only needs to be done once at the beginning of the simulation,
its runtime cost for communication as well as its memory demand can become prohibitive
when going to large problems.
The RR approach, in turn, alleviates the issue of excessive growth of memory demand.
Yet, the number of necessary RR iterations equals the number of processes~$\nproc$,
rendering this approach impractical for~$\nproc \gg 1$ (especially as it has to be applied in every time/load step).
Although the binning approach proposed in \secref{sec:Binning} needs to be applied in every time/load step,
it appears as the only approach without impractical restrictions when going to large problem sizes:
Through the choice of the number and size of the bins,
the amount of data to be ghosted can be controlled, such that only those {\master} elements will be ghosted,
that are likely to be required during contact search and evaluation.
In sum, the applicability of the binning approach is neither affected by the number of parallel processes
nor greatly impacts the parallel communication or total memory demand.

We will later supplement our assessment with detailed numerical experiments in \secref{sec:NumExTwoCubesContactWeakScaling},
but want to anticipate the main finding here:
For the largest problems with $25M$ mesh nodes and 25k interface nodes,
the process with the largest ghosting demand asks for the redundant ghosting of $25921$ nodes,
while binning reduces this number to $1212$ nodes,
which amounts to a reduction of more than $20\times$.
On average across all MPI ranks,
these numbers can be improved through load balancing
which will be introduced in \secref{sec:LoadBalancing}.

\section{Balancing the work load among multiple parallel processes}
\label{sec:LoadBalancing}

Now, we discuss strategies for an optimal distribution of the work load to multiple parallel processes.
These strategies are intended to satisfy the requirement~{\reqSlave} from \secref{sec:EvaluationMortarIntegrals}.
We assume that requirement~{\reqMaster} has already been satisfied by any of the methods described in \secref{sec:DistributedMemory}
and, thus, all data is accessible whenever needed.

In \secsrangeref{sec:ScalabilityConcepts}{sec:InterfaceDD}, we first present some general considerations applicable to all type of mortar interface problems,
before we move to the specific scenario of dynamically evolving contact problems in \secref{sec:InterfaceDDDynamic}.

\subsection{The concepts of strong and weak scalability}
\label{sec:ScalabilityConcepts}

When assessing the performance of a parallel code and/or algorithm,
an important question is
whether adding more computational resources will actually {\speedup} the algorithm's performance at the proper rate.
Two concepts are commonly followed and investigated:
\begin{itemize}
\item For a fixed problem size,
\define{strong scalability} is given,
if the computational time diminishes at the same rate as the used hardware resources grow.
The \define{strong scaling limit} is reached,
when increasing the hardware resources does not lead to a further reduction of computational time.
See \cite{Amdahl1967a}.
\item \define{Weak scalability} expects a constant computational time
when increasing the problem size and the parallel resources at the same rate,
{\ie} when the work load per process is kept constant.
See \cite{Gustafson1988a}.
\end{itemize}
As it is well established in many research and application codes (and also in our code base {\baci}~\cite{BaciURL}),
weak scalability of the finite element evaluation of the pure bulk field ({\ie} volume element evaluation)
without the presence of any mortar interface can be achieved under uniform mesh refinement.

\subsection{Curse of dimensionality}
\label{sec:ClashOfDimensionality}

In surface-coupled problems with~$\ndim$ spatial dimensions,
the coupling surface is always a~$\ndim-1$ dimensional geometric entity.
Originally described in \cite{Bellmann1957a},
this \define{curse of dimensionality} between the bulk and the interface discretization
becomes problematic under uniform mesh refinement.
Denoting the characteristic mesh size with~$\hEle$,
the number of unknowns in the bulk discretization grows at~$\HOT{\hEle}{\ndim}$
while the surface discretization of the coupling interface exhibits a growth rate of~$\HOT{\hEle}{\ndim-1}$ only.

This becomes evident in practice
when a first and simple DD of the interface discretizations is now obtained
by aligning the interface subdomains of the {\slave} and {\master} side
with the subdomains of the underlying bulk discretizations.
Although this approach is straightforward to implement and also avoids off-process assembly, thus reducing parallel communication,
it does not result in an optimal parallel distribution for the evaluation of the mortar coupling terms.
Since computing the interface contributions,
{\ie} the mortar segmentation process, integration and assembly of the mortar matrices~$\MoD$ and~$\MoM$
to only name the most important tasks,
is all done on the {\slave} interface discretization,
all numerical tasks might be performed by very few parallel processes only,
while others idle.

\begin{figure}
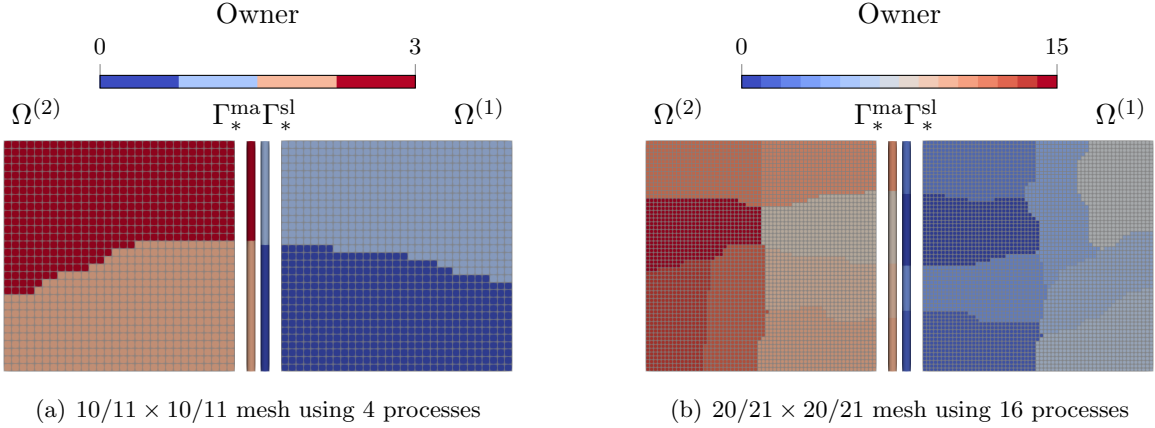

\centering
\hfill
\subfigure[$10/11\times 10/11$ mesh using $4$ processes]{\label{fig:DimClashNoRedistSmall}\input{fig/meshtying2D_two_blocks/case_vol_01_owner.tex}}
\hfill
\subfigure[$20/21\times 20/21$ mesh using $16$ processes]{\label{fig:DimClashNoRedistLarge}\input{fig/meshtying2D_two_blocks/case_vol_02_owner.tex}}
~\hfill
\caption{Perfect alignment of bulk and interface domain decomposition:
Subdomains in the {\master} side of the interface (left vertical strip) coincide with the {\master} side's bulk discretization (left square),
while the {\slave} side's interface subdomains (right vertical strip) are aligned with the {\slave} side's bulk discretization (right square).}
\label{fig:DimClashNoRedist}
\end{figure}

For simplicity of visualization, this is illustrated using a two-dimensional {\meshtying} problem in \figref{fig:DimClashNoRedist},
where the domain decomposition of the mortar interface's {\slave} and {\master} side is fully aligned with the underlying bulk discretizations.
Considering a coarse discretization distributed to four parallel processes as shown in \figref{fig:DimClashNoRedistSmall},
the {\slave} interface is divided into two subdomains and the {\master} interface is owned by two processes only as well.
Consequently, there are two processes, that do not own a share of the {\slave} interface,
and two other processes not owning any node of the {\master} interface.
In \figref{fig:DimClashNoRedistLarge},
the mesh has been refined by a factor of two in each direction, and $16$ processes have been used
such that the load per process remains constant in the bulk discretization.
While the bulk discretization is now split into $16$ subdomains,
the {\slave} and {\master} interface are shared only among four processes each.
In sum, only $4$ processes tackle the expensive evaluation of mortar terms on the {\slave} side of the interface,
while $12$ processes are completely left out.
Even in these small and only two-dimensional problem, it becomes evident
that the alignment of interface subdomains with bulk subdomains potentially leaves a huge fraction of all processes idle during interface evaluation.
While it is true that using more processes improves the parallelization of the bulk discretization,
it does not necessarily contribute to a good and scalable parallelization of the interface computations.
We stress that this issue is even more pronounced for the three-dimensional case and larger numbers of parallel processes.

\subsection{Improving the domain decomposition of interface discretizations}
\label{sec:InterfaceDD}

\begin{figure}
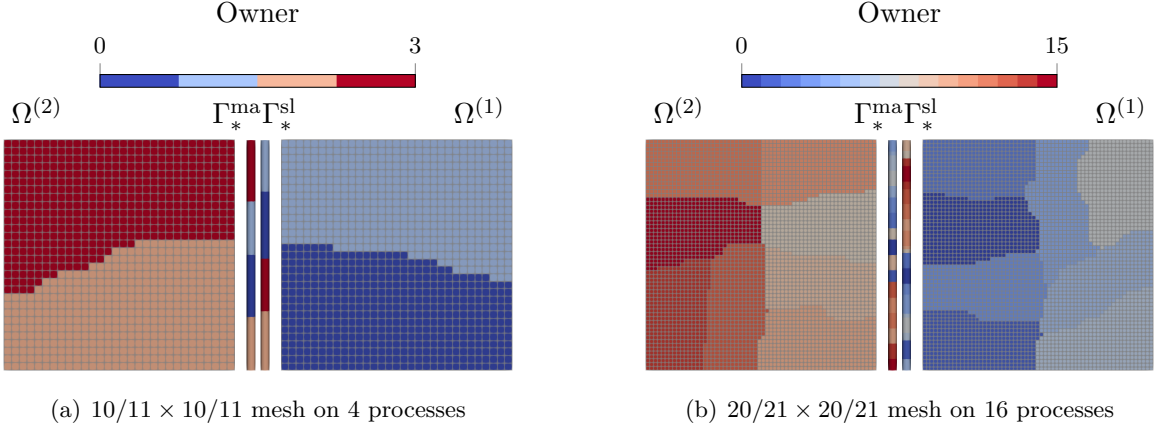

\centering
\hfill
\subfigure[$10/11\times 10/11$ mesh on $4$ processes]{\label{fig:DimClashRedistToAllSmall}\input{fig/meshtying2D_two_blocks/case_vol_03_owner.tex}}
\hfill
\subfigure[$20/21\times 20/21$ mesh on $16$ processes]{\label{fig:DimClashRedistToAllLarge}\input{fig/meshtying2D_two_blocks/case_vol_04_owner.tex}}
~\hfill
\caption{Independent interface domain decomposition using all~$\nproc$ parallel processes: }
\label{fig:InterfaceDDToAll}
\end{figure}

To overcome the curse of dimensionality and to satisfy~{\reqSlave}, we allow the {\slave} and {\master} side of the interface
to be decomposed into subdomains independently from the underlying bulk discretizations
in order to achieve optimal parallel scalability of the computational tasks associated with both
the integration and assembly in the bulk domains~$\procID{\dom}{1}$ and~$\procID{\dom}{2}$
as well as integration and assembly on the interfaces~$\indexedSlave{\Bound}$ and~$\indexedMaster{\Bound}$.
In a first and straightforward approach, one can divide both interfaces~$\indexedSlave{\Bound}$ and~$\indexedMaster{\Bound}$ into $\nproc$ subdomains,
such that each parallel process handles a portion of the interface as illustrated in \figref{fig:InterfaceDDToAll}.
This is particularly important for the {\slave} side
which needs to perform all computations related to the integration of the mortar terms in~\eqref{eq:DiscreteWeakFormInterfaceFlux}.

For the coarse and the fine mesh, both~$\indexedSlave{\Bound}$ and~$\indexedMaster{\Bound}$ are distributed to~$4$ and~$16$ parallel processes, respectively.
A clear advantage of this strategy is that all parallel processes participate in the interface treatment, so idling is mostly avoided.
However, the fine mesh already indicates that the interface subdomains may become very small, {\ie} they consist only of a few elements.
Recalling the curse of dimensionality outlined in \secref{sec:ClashOfDimensionality},
this will become an issue at large scale where the bulk field is divided into~$\nproc$ subdomains of reasonable size,
while the subdomain size of the interface decreases when refining the mesh and adding parallel processes at the same rate.
Having many but very small interface subdomains does not leave any process idle,
but also yields interface subdomains with a large surface-to-volume ratio which indicates an increasing communication overhead.
In sum, this strategy distributes the computational work of the interface evaluation more evenly to all processes
than just adopting the interface subdomains from the underlying bulk discretization.
The numerical experiments in \secref{sec:NumEx} confirm this statement.

Conceptually, there is still room for further optimizations, in particular related to parallel communication among processes,
{\eg} by setting a lower bound~$\nEleMin$ on the number of elements per interface subdomain to reduce the communication overhead.
Such an approach needs to compromise between the amount of parallel communication among processes and the number of idling processes.
In this work, we have refrained from exploring this research direction,
since the distribution of the interface to all parallel processes already delivers satisfying scaling behavior for many practical applications.

\subsection{Interface domain decompositions for dynamically evolving interfaces}
\label{sec:InterfaceDDDynamic}

In many applications, the interface configuration evolves over time,
{\eg} as in contact problems with large sliding or contact of rolling bodies.
In such cases, the interface DD can come out of balance,
resulting in some processes to do significantly more work than others, which possibly idle.
Then, a rebalancing can become necessary
to distribute the computational work evenly to all participating processes.

In each time step, we track the time spent in the evaluation of all mortar terms for each processor
as well as the number of {\slave} elements per processor.
We then estimate the imbalance among all processes by
\begin{align}
\imbalanceRatioTime = \frac{\max_\indProc \left(\tEvaluateOnProc{\indProc}\right)}{\min_\indProc \left(\tEvaluateOnProc{\indProc}\right)}, \quad
\imbalanceRatioEles = \frac{\max_\indProc \left(\procID{\nEleSlave}{\indProc}\right)}{\min_\indProc \left(\procID{\nEleSlave}{\indProc}\right)}
\end{align}
with~$\imbalanceRatioTime$ and~$\imbalanceRatioEles$ denoting the imbalance in contact evaluation time and number of {\slave} elements per processor, respectively.
The theoretical optimum of a perfect balancing of the mortar-related workload is given for $\imbalanceRatioTime = 1$ and $\imbalanceRatioEles = 1$, respectively,
{\ie} when all processes spend exactly the same time in mortar evaluation and when all processes own the exact same number of {\slave} elements.
If in any time step these imbalance estimates exceed user-given thresholds~$\imbalanceThresholdTime \geq 1$ and~$\imbalanceRatioEles \geq 1$
for contact evaluation time and number of {\slave} elements per processor, respectively,
{\ie} if
\begin{align}
\label{eq:TriggerRebalancing}
\imbalanceRatioTime \geq \imbalanceThresholdTime
\quad\vee\quad
\imbalanceRatioEles \geq \imbalanceRatioEles,
\end{align}
then we re-compute the interface DD to obtain a DD with better load balancing.
As this load balancing procedure is triggered dynamically by the current state of the simulation,
we refer to it as \define{dynamic load balancing}.

Naturally, $\imbalanceThresholdTime = 1$ will trigger rebalancing in every time step,
such that each time step can rely on the best possible interface DD.
In practice, the cost for rebalancing needs to be taken into account,
such that practical computations require $\imbalanceThresholdTime > 1$.
We will study the impact of the actual choice of~$\imbalanceThresholdTime$ on the run time in \secref{sec:NumExRollingCyl}.

The main difference between the two imbalance measures~$\imbalanceRatioTime$ and~$\imbalanceRatioEles$ is
that~$\imbalanceRatioEles$ does not account for the time to evaluate a given {\master}/{\slave} pair,
while~$\imbalanceRatioTime$ relies on actual wall clock timings.
Thus, situations with~$\imbalanceRatioEles \gg 1$, but $\imbalanceRatioTime$ fairly close to $1$ can occur,
if the contact search identifies a huge number of pairs of {\master} and {\slave} elements as close to each other,
but the subsequent mortar evaluation cannot find a valid projection and, thus,
most of the computational work to evaluate~\eqref{eq:DiscreteWeakFormInterfaceFlux} is skipped for such pairs of elements.
In sum, the time-based trigger~$\imbalanceRatioTime$ is expected to be more effective
to avoid idling processes in practical simulations.

\subsection{Implication on finite element assembly and communication patterns}

Although the {\slave} side's interface discretization might exhibit its independent DD to improve scalability,
all system quantities, {\eg} the Jacobian matrix~$\mat\jacobian$ and the residual vector~$\mao\res$,
are distributed among parallel processes following the DD of the underlying bulk discretizations.
After evaluation of the mortar element matrices defined in~\eqref{eq:DiscreteWeakFormInterfaceFlux}
within a mortar element in interface subdomain~$\subdomainID{\Bound}{\altIndSubdomain}$, $\altIndSubdomain\in\{0, \hdots, \nsubdomain-1\}$,
on process~$\altIndProc\in\{0,\hdots,\nproc-1\}$,
a contribution to~$\mat\jacobian$ and~$\mao\res$ associated with node~$\indNode$ in~$\subdomainID{\dom}{\indSubdomain}$, $\indSubdomain\in\{0, \hdots, \nsubdomain-1\}$,
owned by process~$\indProc\in\{0,\hdots,\nproc-1\}$
can only be assembled by process~$\indProc$.
Hence, if $\indProc\neq\altIndProc$, communication is required to send data from process~$\altIndProc$ to process~$\indProc$
in order to assemble into global system quantities.
From the perspective of the evaluating process, this is referred to as \define{off-process assembly}.
Communication can only be avoided if and only if~$\indProc=\altIndProc$.

It is true that off-process assembly increases the amount of communication and, thus, puts a cost burden onto the entire algorithm.
Although this is not desirable, it is usually the much cheaper price to pay
than to just stick to the one-to-one matching of interface and underlying bulk DDs.
The {\speedup} of the cost-intensive evaluation of mortar terms through an independent DD of the interface discretizations
easily amortizes the additional cost of communication related to off-process assembly.
We will study timings of the mortar evaluation and off-process assembly in detail
in the numerical experiments in \secref{sec:NumEx}.

\section{Numerical experiments}
\label{sec:NumEx}

We first study parallel redistribution and scalability in a simple two-block contact example
in \secref{sec:NumExTwoCubesContact}
before moving on to dynamic contact problems in \secref{sec:NumExRollingCyl}.

All computations are done with our {\inhouse} {\multiphysics} research code {\baci}~\cite{BaciURL}.
All scaling studies have been run on our {\inhouse} cluster
(20 nodes with 2x Intel Xeon Gold 5118 (Skylake-SP) 12 core CPUs, 196 GB RAM per node, Mellanox Infiniband Interconnect).

\subsection{Contact of two cubes}
\label{sec:NumExTwoCubesContact}

For a first assessment of the scalability of the contact evaluation,
we consider a simple two-block contact problem
with a small block (dimensions $0.8 \times 0.8 \times 0.8$)
and a slightly bigger block (dimensions $1.0 \times 1.0 \times 1.0$),
where contact will occur between two flat surfaces of the blocks.
To reduce the complexity of the contact problem and to exclude nonlinearities due to changes in the contact active set,
the faces opposite to the contact interface are fixed with Dirichlet boundary conditions,
while the blocks initially penetrate each other at the contact interface by $0.001$.
The smaller block acts as the {\slave} side and its entire contact area is already initialized as ``active''.
Application of the contact algorithms will then result in a slight compression of both blocks, such that the initial penetration vanishes.
This problem setup allows to distill the computational effort spent on the redistribution, ghosting, and contact evaluation.
In fact, for the parallel scaling studies, we only evaluate all contact terms, but then do not even solve the contact problem
to allow for an even more concise focus on the scaling behavior of the contact evaluation.

Both blocks use a Neo--Hooke material with Young's modulus~$E=10$ and Poisson's ratio $\nu=0.3$.
Denoting the mesh refinement factor with~$\meshfac$,
both blocks are discretized with $5\meshfac$ linear hexahedral elements along their edges.

\begin{figure}
\begin{center}
\input{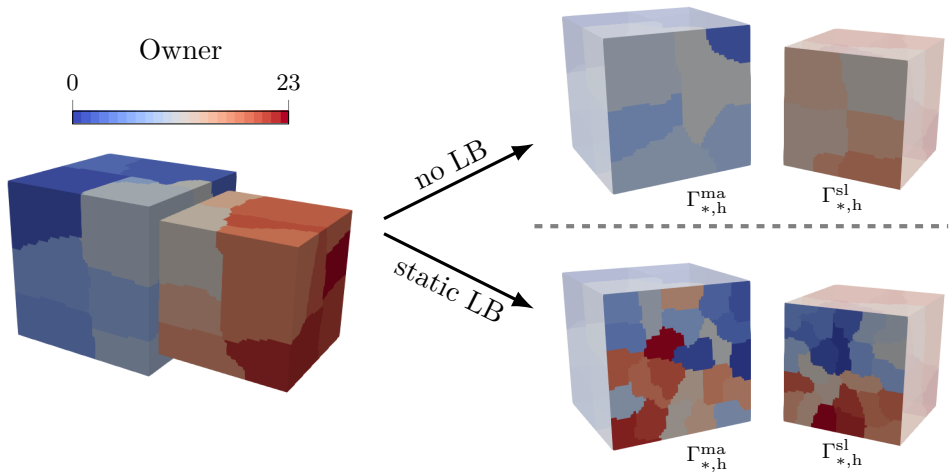}
\caption{Two cubes in contact: colors represent the owning MPI rank of a volume/interface subdomain.}
\label{fig:TwoCubesOwnerPlot}
\end{center}
\end{figure}

As an exemplary visualization,
\figref{fig:TwoCubesOwnerPlot} illustrates the assignment of subdomains to MPI ranks for a simulation with 24 MPI ranks.
Since the discretization of both blocks uses the same number of elements per block,
the volume DD exhibits 12 subdomains for each block.
Without load balancing,
the interface DD evidently matches the underlying volume DD
({\cf} the top right picture in \figref{fig:TwoCubesOwnerPlot}).
In particular,
the {\slave} side of the interface is shared by only 6 (out of 24) processes,
such that the remaining 18 processes idle during the expensive mortar evaluation.
While the DD of the solid volume is not affected by the interface load balancing,
the interface DD now yields 24 subdomains for both sides of the interface
({\cf} the bottom right picture in \figref{fig:TwoCubesOwnerPlot}).
This allows to share the computational work for the mortar evaluation among \emph{all} processes.

\subsubsection{Weak scaling}
\label{sec:NumExTwoCubesContactWeakScaling}

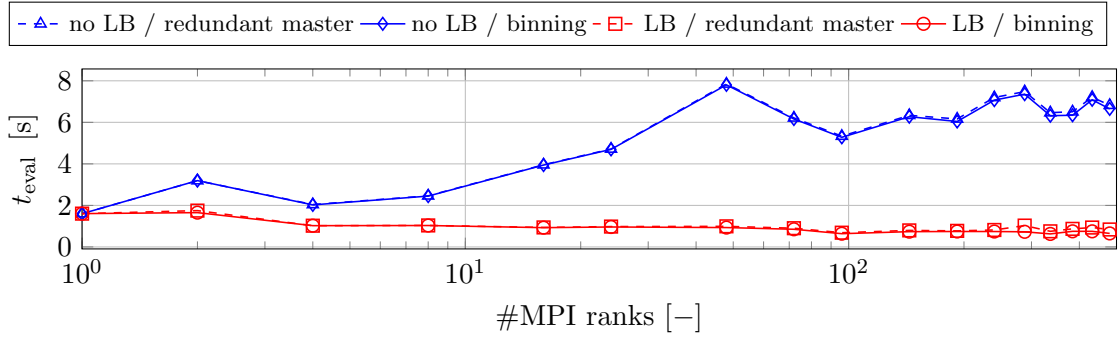
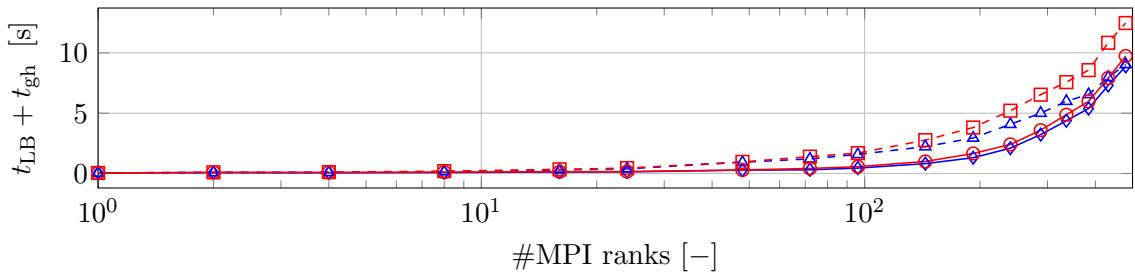
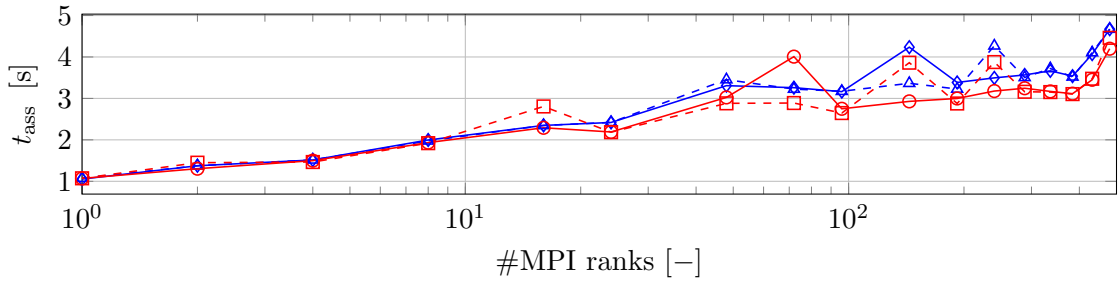
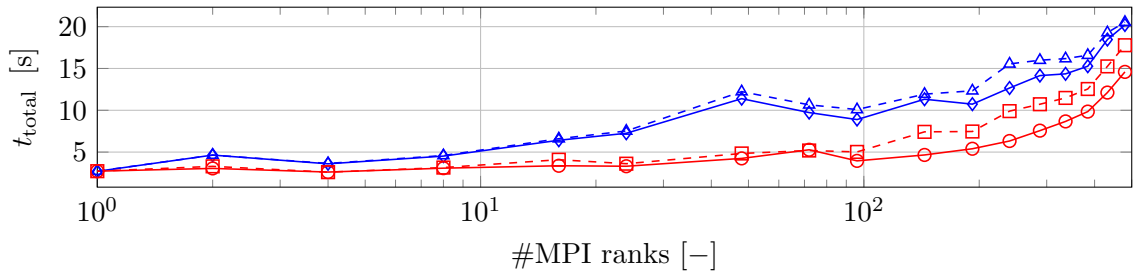
\begin{figure}

\tikzstyle{noneRedundant}=[semithick, dashed, blue, mark=triangle, mark options={scale=1.2, solid, semithick}]
\tikzstyle{noneBinning}=[semithick, solid, blue, mark=diamond, mark options={scale=1.2, solid, semithick}]
\tikzstyle{staticRedundant}=[semithick, dashed, red, mark=square, mark options={scale=1.2, solid, semithick}]
\tikzstyle{staticBinning}=[semithick, solid, red, mark=o, mark options={scale=1.2, solid, semithick}]

\begin{center}

\begin{tikzpicture}
\begin{axis}[%
  hide axis,
  xmin=0,
  xmax=1,
  ymin=0,
  ymax=1,
  legend cell align=left,
  legend style={font=\footnotesize, at={(0.5,1.2)}, anchor=center},
  legend columns = 4
  ]

\addlegendimage{noneRedundant}
\addlegendentry{no LB / redundant {\master}}
\addlegendimage{noneBinning}
\addlegendentry{no LB / binning}

\addlegendimage{staticRedundant}
\addlegendentry{LB / redundant {\master}}
\addlegendimage{staticBinning}
\addlegendentry{LB / binning}

\end{axis}
\end{tikzpicture}

\subfigure[Contact evaluation time]{
\label{fig:ContactTwoCubesWeakScalingContactEvaluationTime}
\begin{tikzpicture}

\pgfplotstableread{data/cubes_weak_scaling/cubes_weak_None_redundant_master_scaling.data}\noneRedundant
\pgfplotstableread{data/cubes_weak_scaling/cubes_weak_None_binning_scaling.data}\noneBinning
\pgfplotstableread{data/cubes_weak_scaling/cubes_weak_Static_redundant_master_scaling.data}\staticRedundant
\pgfplotstableread{data/cubes_weak_scaling/cubes_weak_Static_binning_scaling.data}\staticBinning

\begin{axis}[scale only axis,
  xmin=1,
  xmax=500,
  xmode=log,
  width=0.8\textwidth, height=0.1\textheight,
  ylabel={$\tEvaluate~\left[\second\right]$},
  xlabel={\#MPI ranks~$\left[-\right]$},
  grid=major]

\addplot[noneRedundant] table[x = num_procs, y = contact_time] from \noneRedundant;
\addplot[noneBinning] table[x = num_procs, y = contact_time] from \noneBinning;

\addplot[staticRedundant] table[x = num_procs, y = contact_time] from \staticRedundant;
\addplot[staticBinning] table[x = num_procs, y = contact_time] from \staticBinning;

\end{axis}
\end{tikzpicture}
} 
\subfigure[Time for redistribution and ghosting]{
\label{fig:ContactTwoCubesWeakScalingRedistGhostTime}
\begin{tikzpicture}

\pgfplotstableread{data/cubes_weak_scaling/cubes_weak_None_redundant_master_scaling.data}\noneRedundant
\pgfplotstableread{data/cubes_weak_scaling/cubes_weak_None_binning_scaling.data}\noneBinning
\pgfplotstableread{data/cubes_weak_scaling/cubes_weak_Static_redundant_master_scaling.data}\staticRedundant
\pgfplotstableread{data/cubes_weak_scaling/cubes_weak_Static_binning_scaling.data}\staticBinning

\begin{axis}[scale only axis,
  xmin=1,
  xmax=500,
  xmode=log,
  width=0.8\textwidth, height=0.1\textheight,
  ylabel={$\tRedistribute + \tGhosting~\left[\second\right]$},
  xlabel={\#MPI ranks~$\left[-\right]$},
  grid=major]

\addplot[noneRedundant] table[x = num_procs, y = redist_ghost_time] from \noneRedundant;
\addplot[noneBinning] table[x = num_procs, y = redist_ghost_time] from \noneBinning;

\addplot[staticRedundant] table[x = num_procs, y = redist_ghost_time] from \staticRedundant;
\addplot[staticBinning] table[x = num_procs, y = redist_ghost_time] from \staticBinning;

\end{axis}
\end{tikzpicture}
} 
\subfigure[Time for assembly of contact terms in linear system]{
\label{fig:ContactTwoCubesWeakScalingAssemblyTime}
\begin{tikzpicture}

\pgfplotstableread{data/cubes_weak_scaling/cubes_weak_None_redundant_master_scaling.data}\noneRedundant
\pgfplotstableread{data/cubes_weak_scaling/cubes_weak_None_binning_scaling.data}\noneBinning
\pgfplotstableread{data/cubes_weak_scaling/cubes_weak_Static_redundant_master_scaling.data}\staticRedundant
\pgfplotstableread{data/cubes_weak_scaling/cubes_weak_Static_binning_scaling.data}\staticBinning

\begin{axis}[scale only axis,
  xmin=1,
  xmax=500,
  xmode=log,
  width=0.8\textwidth, height=0.1\textheight,
  ylabel={$\tAssemble~\left[\second\right]$},
  xlabel={\#MPI ranks~$\left[-\right]$},
  grid=major]

\addplot[noneRedundant] table[x = num_procs, y = modify_time] from \noneRedundant;
\addplot[noneBinning] table[x = num_procs, y = modify_time] from \noneBinning;

\addplot[staticRedundant] table[x = num_procs, y = modify_time] from \staticRedundant;
\addplot[staticBinning] table[x = num_procs, y = modify_time] from \staticBinning;

\end{axis}
\end{tikzpicture}
} 

\subfigure[Total time for contact evaluation]{
\label{fig:ContactTwoCubesWeakScalingTotalContactTime}
\begin{tikzpicture}

\pgfplotstableread{data/cubes_weak_scaling/cubes_weak_None_redundant_master_scaling.data}\noneRedundant
\pgfplotstableread{data/cubes_weak_scaling/cubes_weak_None_binning_scaling.data}\noneBinning
\pgfplotstableread{data/cubes_weak_scaling/cubes_weak_Static_redundant_master_scaling.data}\staticRedundant
\pgfplotstableread{data/cubes_weak_scaling/cubes_weak_Static_binning_scaling.data}\staticBinning

\begin{axis}[scale only axis,
  xmin=1,
  xmax=500,
  xmode=log,
  width=0.8\textwidth, height=0.1\textheight,
  ylabel={$\tContactTotal~\left[\second\right]$},
  xlabel={\#MPI ranks~$\left[-\right]$},
  grid=major]

\addplot[noneRedundant] table[x = num_procs, y = total_eval_assemble_time] from \noneRedundant;
\addplot[noneBinning] table[x = num_procs, y = total_eval_assemble_time] from \noneBinning;

\addplot[staticRedundant] table[x = num_procs, y = total_eval_assemble_time] from \staticRedundant;
\addplot[staticBinning] table[x = num_procs, y = total_eval_assemble_time] from \staticBinning;

\end{axis}
\end{tikzpicture}
} 

\caption{Weak scaling of contact time for two cubes}
\label{fig:ContactTwoCubesWeakScaling}
\end{center}
\end{figure}

We perform a weak scaling study.
The smallest problem using 1 MPI rank consists of 55,566 displacement unknowns,
while 441/400 nodes/elements reside on the {\slave} side of the contact interface.
The largest problem using 480 MPI ranks contains 25,039,686 displacement unknowns,
with 25,921/25,600 nodes/elements located on the {\slave} side of the contact interface.
We target a load of $\approx$50k displacement unknowns per MPI rank under weak scaling conditions.
Timing results are shown in \figref{fig:ContactTwoCubesWeakScaling}.
With load balancing, the pure contact evaluation time remains constant under weak scaling conditions
as shown in \figref{fig:ContactTwoCubesWeakScalingContactEvaluationTime}
and as expected for finite element evaluations.
Manifesting the curse of dimensionality described in \secref{sec:ClashOfDimensionality} though,
the case without load balancing does not equally benefit from adding hardware resources
since most of the additional processes do not participate in the mortar evaluation.
While the choice of load balancing does not impact the serial case ($\nproc = 1$) of course,
the contact evaluation without load balancing requires twice as much time on 2, 4, and 8 MPI ranks than with load balancing,
since the processes owning a piece of the {\master} side of the interface do not contribute to the contact evaluation.
For an increasing number of MPI ranks,
this gap increases.

Regarding the time spent in redistribution and extending the interface ghosting, $\tRedistribute+\tGhosting$,
an increase with an increasing number of MPI ranks is expected,
as the size of the MPI communicator grows and, thus, mandates increased communication.
Obviously, this time component is rather independent of the parallel distribution,
but is largely impacted by the ghosting strategy:
Since the redundant ghosting of the {\master} side requires to communicate all interface nodes and elements of the {\master} side to all MPI ranks,
the timings for redundant ghosting exceed the time for ghosting via the geometrically motivated binning approach,
where the amount of data to be communicated among processes is reduced based on geometric information.

It becomes evident from \figref{fig:ContactTwoCubesWeakScalingAssemblyTime},
that the time for assembling of all contact terms into the global linear system is only slightly impacted by load balancing,
while the impact of the ghosting strategy appears to be negligible.

Finally, we assess the total cost of contact evaluation
which is the most relevant target quantity for practical applications.
It is given by the total time~$\tContactTotal = \tRedistribute+\tGhosting+\tEvaluate+\tAssemble$
for (possibly) redistributing, ghosting, evaluation, and assembly of the contact interface
and is shown in \figref{fig:ContactTwoCubesWeakScalingTotalContactTime}.
Again, ghosting via binning (``binning'') results in a lower total time~$\tContactTotal$
than the redundant ghosting of the {\master} side (``redundant {\master}'').
Moreover, load balancing (``LB'') allows all MPI ranks to participate in the evaluation of the contact terms,
yielding a faster total contact time than without load balancing (``no LB'').
Dominated by the contact evaluation time~$\tEvaluate$,
the case without load balancing does not scale beyond 8 MPI ranks,
while load balancing shows good weak scalability up to 200 MPI ranks.
Overall, our proposed strategy of load balancing in combination with binning delivers the fastest contact evaluation for all mesh sizes
and also features the smallest increase in total contact time when increasing the problem size.

\begin{figure}

\tikzstyle{noneRedundantColMax}=[semithick, dashed, blue, mark=square, mark options={scale=1.2, solid, semithick}]
\tikzstyle{noneBinningColMax}=[semithick, solid, blue, mark=square, mark options={scale=1.2, solid, semithick}]
\tikzstyle{noneBinningRowMax}=[very thick, dotted, blue, mark=square, mark options={scale=1.2, solid, semithick}]

\tikzstyle{staticRedundantColMax}=[semithick, dashed, red, mark=o, mark options={scale=1.2, solid, semithick}]
\tikzstyle{staticBinningColMax}=[semithick, solid, red, mark=o, mark options={scale=1.2, solid, semithick}]
\tikzstyle{staticBinningRowMax}=[very thick, dotted, red, mark=o, mark options={scale=1.2, solid, semithick}]

\begin{center}

\begin{tikzpicture}
\begin{axis}[%
  hide axis,
  xmin=0,
  xmax=1,
  ymin=0,
  ymax=1,
  legend cell align=left,
  legend style={font=\footnotesize, at={(0.5,1.2)}, anchor=center},
  legend columns = 3
  ]

\addlegendimage{noneRedundantColMax}
\addlegendentry{no LB / redundant {\master} / max. \#ghosted}
\addlegendimage{noneBinningColMax}
\addlegendentry{no LB / binning / max. \#ghosted}
\addlegendimage{noneBinningRowMax}
\addlegendentry{no LB / max. \#owned}

\addlegendimage{staticRedundantColMax}
\addlegendentry{LB / redundant {\master} / max. \#ghosted}
\addlegendimage{staticBinningColMax}
\addlegendentry{LB / binning / max. \#ghosted}
\addlegendimage{staticBinningRowMax}
\addlegendentry{LB / max. \#owned}

\end{axis}
\end{tikzpicture}

%
%
%
%
%
%
%

\begin{tikzpicture}

\pgfplotstableread{data/cubes_weak_scaling/cubes_weak_None_redundant_master_ghosting.data}\noneRedundant
\pgfplotstableread{data/cubes_weak_scaling/cubes_weak_None_binning_ghosting.data}\noneBinning
\pgfplotstableread{data/cubes_weak_scaling/cubes_weak_Static_redundant_master_ghosting.data}\staticRedundant
\pgfplotstableread{data/cubes_weak_scaling/cubes_weak_Static_binning_ghosting.data}\staticBinning

\begin{axis}[scale only axis,
  xmin=0,
  xmax=480,
  ymode=log,
  width=0.8\textwidth, height=0.1\textheight,
  ylabel={\#elements~$\left[-\right]$},
  xlabel={\#MPI ranks~$\left[-\right]$},
  grid=major]

\addplot[noneRedundantColMax] table[x = num_procs, y = m_eles_col_max] from \noneRedundant;
\addplot[noneBinningColMax] table[x = num_procs, y = m_eles_col_max] from \noneBinning;

\addplot[staticRedundantColMax] table[x = num_procs, y = m_eles_col_max] from \staticRedundant;
\addplot[staticBinningColMax] table[x = num_procs, y = m_eles_col_max] from \staticBinning;

\addplot[noneBinningRowMax] table[x = num_procs, y = m_eles_row_max] from \noneBinning;

\addplot[staticBinningRowMax] table[x = num_procs, y = m_eles_row_max] from \staticBinning;

\end{axis}
\end{tikzpicture}

\caption{Number of owned and ghosted elements of the interface's {\master} side under weak scaling conditions}
\label{fig:ContactTwoCubesWeakScalingInterfaceGhosting}
\end{center}
\end{figure}
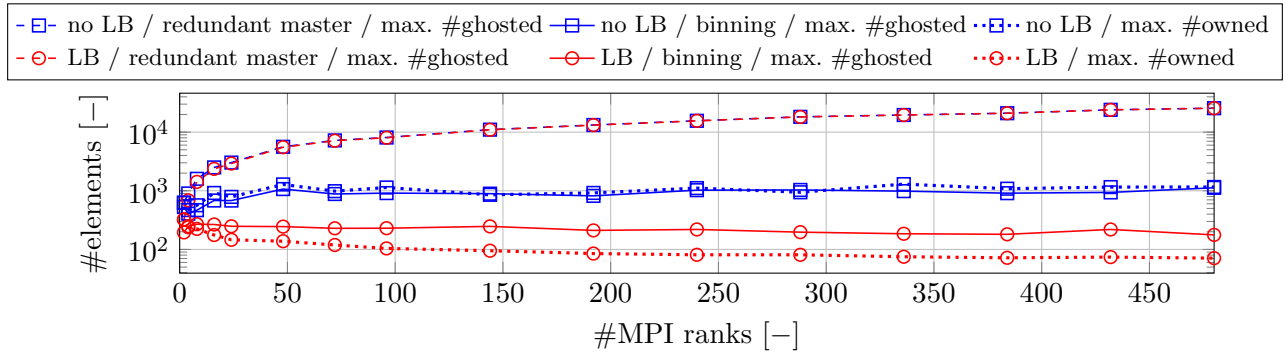

\Figref{fig:ContactTwoCubesWeakScalingInterfaceGhosting} illustrates the impact of both the load balancing and the ghosting strategy
on the number of owned and ghosted {\master} side elements
by reporting the maximum number of elements per MPI rank among all processes.
Naturally, load balancing, where all processes hold a portion of the contact interface,
leads to a lower number of owned entities per processor than no load balancing,
where the interface is stored only by a subset of all processes:
Depicted by the dotted lines,
the number of owned elements per MPI rank is smaller in case of load balancing than without load balancing,
in particular by a factor of $100$ for more than $48$ MPI ranks in this example.
The influence of the ghosting strategy is shown with dashed and solid lines:
While binning (solid lines) just adds a small number of nodes or elements to be ghosted from other processes,
fully redundant ghosting (dashed lines) drastically increases the number of ghosted elements.
For large examples, this increase can exceed two orders of magnitude.
We observe that the number of owned elements is consistently smaller than the number of ghosted elements when using load balancing,
while this is not the case without load balancing.
This peculiarity is just an artifact of the visualization,
since the MPI rank with the maximum number of owned entities is not necessarily the same as the one with the maximum number of ghosted entities.
Using the ratio of ghosted elements to owned elements to indicate the additional overhead in memory and parallel communication
due to the distributed memory paradigm,
we make the following key observation:
ghosting via binning as proposed in \secref{sec:Binning} is much more efficient in terms of memory and parallel communication than fully redundant ghosting.
Please note that the respective diagram for owned and ghosted nodes of the interface's {\master} side essentially looks the same
and, thus, is not shown for the conciseness of the presentation.

\subsubsection{Strong scaling}
\label{sec:NumExTwoCubesContactStrongScaling}

To assess the strong scaling behavior under different load balancing and ghosting strategies,
we study three different meshes and problem sizes detailed in \tabref{tab:TwoCubesMeshesStrongScaling}.
\begin{table}
\begin{center}
\begin{tabular}{c|rr|rr}
\multirow{2}{*}{Mesh ID} & \multicolumn{2}{c}{Bulk domain} & \multicolumn{2}{|c}{{\Slave} interface} \\
~ & number of nodes & number of DoFs & number of nodes & number of elements\\
\hline
2M & 715,822 & 2,147,466 & 5,041 & 4,900\\
5M & 1,769,472 & 5,308,416 & 9,216 & 9,025\\
10M & 3,543,122 & 10,629,366 & 14,641 & 14,400
\end{tabular}
\caption{Three meshes and problem sizes for strong scaling of a two-block contact example}
\label{tab:TwoCubesMeshesStrongScaling}
\end{center}
\end{table}
The strong scaling behavior is reported in \figref{fig:ContactTwoCubesStrongScaling}.
\begin{figure}

\tikzstyle{noneRedundantTwoM}=[thick, dashed, blue, mark=triangle, mark options={scale=0.9, solid, semithick}]
\tikzstyle{noneRedundantFiveM}=[thick, dashed, blue, mark=square, mark options={scale=0.9, solid, semithick}]
\tikzstyle{noneRedundantTenM}=[thick, dashed, blue, mark=o, mark options={scale=0.9, solid, semithick}]

\tikzstyle{noneBinningTwoM}=[thick, solid, blue, mark=triangle, mark options={scale=0.9, solid, semithick}]
\tikzstyle{noneBinningFiveM}=[thick, solid, blue, mark=square, mark options={scale=0.9, solid, semithick}]
\tikzstyle{noneBinningTenM}=[thick, solid, blue, mark=o, mark options={scale=0.9, solid, semithick}]

\tikzstyle{staticRedundantTwoM}=[thick, dashed, red, mark=triangle, mark options={scale=0.9, solid, semithick}]
\tikzstyle{staticRedundantFiveM}=[thick, dashed, red, mark=square, mark options={scale=0.9, solid, semithick}]
\tikzstyle{staticRedundantTenM}=[thick, dashed, red, mark=o, mark options={scale=0.9, solid, semithick}]

\tikzstyle{staticBinningTwoM}=[thick, solid, red, mark=triangle, mark options={scale=0.9, solid, semithick}]
\tikzstyle{staticBinningFiveM}=[thick, solid, red, mark=square, mark options={scale=0.9, solid, semithick}]
\tikzstyle{staticBinningTenM}=[thick, solid, red, mark=o, mark options={scale=0.9, solid, semithick}]

\tikzstyle{perfectScaling}=[thick, dashed, gray]

\begin{center}

\begin{tikzpicture}
\begin{axis}[%
  hide axis,
  xmin=0,
  xmax=1,
  ymin=0,
  ymax=1,
  legend cell align=left,
  legend style={font=\footnotesize, at={(0.5,1.2)}, anchor=center},
  legend columns = 4
  ]

\addlegendimage{noneRedundantTwoM}
\addlegendentry{2M / no LB / redundant {\master}}
\addlegendimage{noneBinningTwoM}
\addlegendentry{2M / no LB / binning}

\addlegendimage{staticRedundantTwoM}
\addlegendentry{2M / LB / redundant {\master}}
\addlegendimage{staticBinningTwoM}
\addlegendentry{2M / LB / binning}

\addlegendimage{noneRedundantFiveM}
\addlegendentry{5M / no LB / redundant {\master}}
\addlegendimage{noneBinningFiveM}
\addlegendentry{5M / no LB / binning}

\addlegendimage{staticRedundantFiveM}
\addlegendentry{5M / LB / redundant {\master}}
\addlegendimage{staticBinningFiveM}
\addlegendentry{5M / LB / binning}

\addlegendimage{noneRedundantTenM}
\addlegendentry{10M / no LB / redundant {\master}}
\addlegendimage{noneBinningTenM}
\addlegendentry{10M / no LB / binning}

\addlegendimage{staticRedundantTenM}
\addlegendentry{10M / LB / redundant {\master}}
\addlegendimage{staticBinningTenM}
\addlegendentry{10M / LB / binning}

%
%
%

\addlegendimage{perfectScaling}
\addlegendentry{perfect scaling}

\end{axis}
\end{tikzpicture}

\subfigure[Contact evaluation time]{
\label{fig:ContactTwoCubesStrongScalingContactEvaluationTime}
\begin{tikzpicture}

\pgfplotstableread{data/cubes_strong_scaling/cubes_strong_2M_None_redundant_master_scaling.data}\noneRedundantTwoM
\pgfplotstableread{data/cubes_strong_scaling/cubes_strong_2M_None_binning_scaling.data}\noneBinningTwoM
\pgfplotstableread{data/cubes_strong_scaling/cubes_strong_2M_Static_redundant_master_scaling.data}\staticRedundantTwoM
\pgfplotstableread{data/cubes_strong_scaling/cubes_strong_2M_Static_binning_scaling.data}\staticBinningTwoM

\pgfplotstableread{data/cubes_strong_scaling/cubes_strong_5M_None_redundant_master_scaling.data}\noneRedundantFiveM
\pgfplotstableread{data/cubes_strong_scaling/cubes_strong_5M_None_binning_scaling.data}\noneBinningFiveM
\pgfplotstableread{data/cubes_strong_scaling/cubes_strong_5M_Static_redundant_master_scaling.data}\staticRedundantFiveM
\pgfplotstableread{data/cubes_strong_scaling/cubes_strong_5M_Static_binning_scaling.data}\staticBinningFiveM

\pgfplotstableread{data/cubes_strong_scaling/cubes_strong_10M_None_redundant_master_scaling.data}\noneRedundantTenM
\pgfplotstableread{data/cubes_strong_scaling/cubes_strong_10M_None_binning_scaling.data}\noneBinningTenM
\pgfplotstableread{data/cubes_strong_scaling/cubes_strong_10M_Static_redundant_master_scaling.data}\staticRedundantTenM
\pgfplotstableread{data/cubes_strong_scaling/cubes_strong_10M_Static_binning_scaling.data}\staticBinningTenM

%

\begin{axis}[scale only axis,
  xmin=1,
  xmax=500,
  xmode=log,
  ymode=log,
  width=0.3\textwidth, height=0.25\textwidth,
  ylabel={$\tEvaluate~\left[\second\right]$},
  xlabel={\#MPI ranks~$\left[-\right]$},
  grid=major]

\addplot[noneRedundantTwoM] table[x = num_procs, y = contact_time] from \noneRedundantTwoM;
\addplot[noneBinningTwoM] table[x = num_procs, y = contact_time] from \noneBinningTwoM;
\addplot[staticRedundantTwoM] table[x = num_procs, y = contact_time] from \staticRedundantTwoM;
\addplot[staticBinningTwoM] table[x = num_procs, y = contact_time] from \staticBinningTwoM;

\addplot[noneRedundantFiveM] table[x = num_procs, y = contact_time] from \noneRedundantFiveM;
\addplot[noneBinningFiveM] table[x = num_procs, y = contact_time] from \noneBinningFiveM;
\addplot[staticRedundantFiveM] table[x = num_procs, y = contact_time] from \staticRedundantFiveM;
\addplot[staticBinningFiveM] table[x = num_procs, y = contact_time] from \staticBinningFiveM;

\addplot[noneRedundantTenM] table[x = num_procs, y = contact_time] from \noneRedundantTenM;
\addplot[noneBinningTenM] table[x = num_procs, y = contact_time] from \noneBinningTenM;
\addplot[staticRedundantTenM] table[x = num_procs, y = contact_time] from \staticRedundantTenM;
\addplot[staticBinningTenM] table[x = num_procs, y = contact_time] from \staticBinningTenM;

%

\addplot[perfectScaling] coordinates {(1,18) (100, 0.18)};

\end{axis}
\end{tikzpicture}
} 
\subfigure[Time for redistribution and ghosting]{
\label{fig:ContactTwoCubesStrongScalingRedistGhostTime}
\begin{tikzpicture}

\pgfplotstableread{data/cubes_strong_scaling/cubes_strong_2M_None_redundant_master_scaling.data}\noneRedundantTwoM
\pgfplotstableread{data/cubes_strong_scaling/cubes_strong_2M_None_binning_scaling.data}\noneBinningTwoM
\pgfplotstableread{data/cubes_strong_scaling/cubes_strong_2M_Static_redundant_master_scaling.data}\staticRedundantTwoM
\pgfplotstableread{data/cubes_strong_scaling/cubes_strong_2M_Static_binning_scaling.data}\staticBinningTwoM

\pgfplotstableread{data/cubes_strong_scaling/cubes_strong_5M_None_redundant_master_scaling.data}\noneRedundantFiveM
\pgfplotstableread{data/cubes_strong_scaling/cubes_strong_5M_None_binning_scaling.data}\noneBinningFiveM
\pgfplotstableread{data/cubes_strong_scaling/cubes_strong_5M_Static_redundant_master_scaling.data}\staticRedundantFiveM
\pgfplotstableread{data/cubes_strong_scaling/cubes_strong_5M_Static_binning_scaling.data}\staticBinningFiveM

\pgfplotstableread{data/cubes_strong_scaling/cubes_strong_10M_None_redundant_master_scaling.data}\noneRedundantTenM
\pgfplotstableread{data/cubes_strong_scaling/cubes_strong_10M_None_binning_scaling.data}\noneBinningTenM
\pgfplotstableread{data/cubes_strong_scaling/cubes_strong_10M_Static_redundant_master_scaling.data}\staticRedundantTenM
\pgfplotstableread{data/cubes_strong_scaling/cubes_strong_10M_Static_binning_scaling.data}\staticBinningTenM

%

\begin{axis}[scale only axis,
  xmin=1,
  xmax=500,
  xmode=log,
  ymode=log,
  width=0.3\textwidth, height=0.25\textwidth,
  ylabel={$\tRedistribute + \tGhosting~\left[\second\right]$},
  xlabel={\#MPI ranks~$\left[-\right]$},
  grid=major]

\addplot[noneRedundantTwoM] table[x = num_procs, y = redist_ghost_time] from \noneRedundantTwoM;
\addplot[noneBinningTwoM] table[x = num_procs, y = redist_ghost_time] from \noneBinningTwoM;
\addplot[staticRedundantTwoM] table[x = num_procs, y = redist_ghost_time] from \staticRedundantTwoM;
\addplot[staticBinningTwoM] table[x = num_procs, y = redist_ghost_time] from \staticBinningTwoM;

\addplot[noneRedundantFiveM] table[x = num_procs, y = redist_ghost_time] from \noneRedundantFiveM;
\addplot[noneBinningFiveM] table[x = num_procs, y = redist_ghost_time] from \noneBinningFiveM;
\addplot[staticRedundantFiveM] table[x = num_procs, y = redist_ghost_time] from \staticRedundantFiveM;
\addplot[staticBinningFiveM] table[x = num_procs, y = redist_ghost_time] from \staticBinningFiveM;

\addplot[noneRedundantTenM] table[x = num_procs, y = redist_ghost_time] from \noneRedundantTenM;
\addplot[noneBinningTenM] table[x = num_procs, y = redist_ghost_time] from \noneBinningTenM;
\addplot[staticRedundantTenM] table[x = num_procs, y = redist_ghost_time] from \staticRedundantTenM;
\addplot[staticBinningTenM] table[x = num_procs, y = redist_ghost_time] from \staticBinningTenM;

%

\addplot[perfectScaling] coordinates {(1,2) (10, 0.2)};

\end{axis}
\end{tikzpicture}
} 

\subfigure[Time for assembly of contact terms into linear system]{
\label{fig:ContactTwoCubesStrongScalingModifyTime}
\begin{tikzpicture}

\pgfplotstableread{data/cubes_strong_scaling/cubes_strong_2M_None_redundant_master_scaling.data}\noneRedundantTwoM
\pgfplotstableread{data/cubes_strong_scaling/cubes_strong_2M_None_binning_scaling.data}\noneBinningTwoM
\pgfplotstableread{data/cubes_strong_scaling/cubes_strong_2M_Static_redundant_master_scaling.data}\staticRedundantTwoM
\pgfplotstableread{data/cubes_strong_scaling/cubes_strong_2M_Static_binning_scaling.data}\staticBinningTwoM

\pgfplotstableread{data/cubes_strong_scaling/cubes_strong_5M_None_redundant_master_scaling.data}\noneRedundantFiveM
\pgfplotstableread{data/cubes_strong_scaling/cubes_strong_5M_None_binning_scaling.data}\noneBinningFiveM
\pgfplotstableread{data/cubes_strong_scaling/cubes_strong_5M_Static_redundant_master_scaling.data}\staticRedundantFiveM
\pgfplotstableread{data/cubes_strong_scaling/cubes_strong_5M_Static_binning_scaling.data}\staticBinningFiveM

\pgfplotstableread{data/cubes_strong_scaling/cubes_strong_10M_None_redundant_master_scaling.data}\noneRedundantTenM
\pgfplotstableread{data/cubes_strong_scaling/cubes_strong_10M_None_binning_scaling.data}\noneBinningTenM
\pgfplotstableread{data/cubes_strong_scaling/cubes_strong_10M_Static_redundant_master_scaling.data}\staticRedundantTenM
\pgfplotstableread{data/cubes_strong_scaling/cubes_strong_10M_Static_binning_scaling.data}\staticBinningTenM

%

\begin{axis}[scale only axis,
  xmin=1,
  xmax=500,
  xmode=log,
  ymode=log,
  width=0.3\textwidth, height=0.25\textwidth,
  ylabel={$\tAssemble~\left[\second\right]$},
  xlabel={\#MPI ranks~$\left[-\right]$},
  grid=major]

\addplot[noneRedundantTwoM] table[x = num_procs, y = modify_time] from \noneRedundantTwoM;
\addplot[noneBinningTwoM] table[x = num_procs, y = modify_time] from \noneBinningTwoM;
\addplot[staticRedundantTwoM] table[x = num_procs, y = modify_time] from \staticRedundantTwoM;
\addplot[staticBinningTwoM] table[x = num_procs, y = modify_time] from \staticBinningTwoM;

\addplot[noneRedundantFiveM] table[x = num_procs, y = modify_time] from \noneRedundantFiveM;
\addplot[noneBinningFiveM] table[x = num_procs, y = modify_time] from \noneBinningFiveM;
\addplot[staticRedundantFiveM] table[x = num_procs, y = modify_time] from \staticRedundantFiveM;
\addplot[staticBinningFiveM] table[x = num_procs, y = modify_time] from \staticBinningFiveM;

\addplot[noneRedundantTenM] table[x = num_procs, y = modify_time] from \noneRedundantTenM;
\addplot[noneBinningTenM] table[x = num_procs, y = modify_time] from \noneBinningTenM;
\addplot[staticRedundantTenM] table[x = num_procs, y = modify_time] from \staticRedundantTenM;
\addplot[staticBinningTenM] table[x = num_procs, y = modify_time] from \staticBinningTenM;

%

\addplot[perfectScaling] coordinates {(1,60) (100, 0.6)};

\end{axis}
\end{tikzpicture}
} 
\subfigure[Total time for contact evaluation]{
\label{fig:ContactTwoCubesStrongScalingTotalContactTime}
\begin{tikzpicture}

\pgfplotstableread{data/cubes_strong_scaling/cubes_strong_2M_None_redundant_master_scaling.data}\noneRedundantTwoM
\pgfplotstableread{data/cubes_strong_scaling/cubes_strong_2M_None_binning_scaling.data}\noneBinningTwoM
\pgfplotstableread{data/cubes_strong_scaling/cubes_strong_2M_Static_redundant_master_scaling.data}\staticRedundantTwoM
\pgfplotstableread{data/cubes_strong_scaling/cubes_strong_2M_Static_binning_scaling.data}\staticBinningTwoM

\pgfplotstableread{data/cubes_strong_scaling/cubes_strong_5M_None_redundant_master_scaling.data}\noneRedundantFiveM
\pgfplotstableread{data/cubes_strong_scaling/cubes_strong_5M_None_binning_scaling.data}\noneBinningFiveM
\pgfplotstableread{data/cubes_strong_scaling/cubes_strong_5M_Static_redundant_master_scaling.data}\staticRedundantFiveM
\pgfplotstableread{data/cubes_strong_scaling/cubes_strong_5M_Static_binning_scaling.data}\staticBinningFiveM

\pgfplotstableread{data/cubes_strong_scaling/cubes_strong_10M_None_redundant_master_scaling.data}\noneRedundantTenM
\pgfplotstableread{data/cubes_strong_scaling/cubes_strong_10M_None_binning_scaling.data}\noneBinningTenM
\pgfplotstableread{data/cubes_strong_scaling/cubes_strong_10M_Static_redundant_master_scaling.data}\staticRedundantTenM
\pgfplotstableread{data/cubes_strong_scaling/cubes_strong_10M_Static_binning_scaling.data}\staticBinningTenM

%

\begin{axis}[scale only axis,
  xmin=1,
  xmax=500,
  xmode=log,
  ymode=log,
  width=0.3\textwidth, height=0.25\textwidth,
  ylabel={$\tContactTotal~\left[\second\right]$},
  xlabel={\#MPI ranks~$\left[-\right]$},
  grid=major]

\addplot[noneRedundantTwoM] table[x = num_procs, y = total_eval_assemble_time] from \noneRedundantTwoM;
\addplot[noneBinningTwoM] table[x = num_procs, y = total_eval_assemble_time] from \noneBinningTwoM;
\addplot[staticRedundantTwoM] table[x = num_procs, y = total_eval_assemble_time] from \staticRedundantTwoM;
\addplot[staticBinningTwoM] table[x = num_procs, y = total_eval_assemble_time] from \staticBinningTwoM;

\addplot[noneRedundantFiveM] table[x = num_procs, y = total_eval_assemble_time] from \noneRedundantFiveM;
\addplot[noneBinningFiveM] table[x = num_procs, y = total_eval_assemble_time] from \noneBinningFiveM;
\addplot[staticRedundantFiveM] table[x = num_procs, y = total_eval_assemble_time] from \staticRedundantFiveM;
\addplot[staticBinningFiveM] table[x = num_procs, y = total_eval_assemble_time] from \staticBinningFiveM;

\addplot[noneRedundantTenM] table[x = num_procs, y = total_eval_assemble_time] from \noneRedundantTenM;
\addplot[noneBinningTenM] table[x = num_procs, y = total_eval_assemble_time] from \noneBinningTenM;
\addplot[staticRedundantTenM] table[x = num_procs, y = total_eval_assemble_time] from \staticRedundantTenM;
\addplot[staticBinningTenM] table[x = num_procs, y = total_eval_assemble_time] from \staticBinningTenM;

%

\addplot[perfectScaling] coordinates {(1,60) (100, 0.6)};

\end{axis}
\end{tikzpicture}
} 

\caption{Two-block contact: strong scaling of contact time}
\label{fig:ContactTwoCubesStrongScaling}
\end{center}
\end{figure}
While meshes 2M and 5M could be run in serial, mesh 10M did not fit into the memory of a single core.
Hence, the graphs for the mesh 10M start at 3 MPI ranks, while 2M and 5M start at 1 MPI rank.

Regarding the pure contact evaluation time~$\tEvaluate$ depicted in \figref{fig:ContactTwoCubesStrongScalingContactEvaluationTime},
the curse of dimensionality as described in \secref{sec:ClashOfDimensionality}
leads to insufficient scaling behavior for the case without load balancing.
For some cases,
the contact evaluation time does not deacrease (or even slightly increase) when adding more processes,
({\cf} the mesh '2M' without load balancing executed on 3, 6, and 12 MPI ranks for example).
In contrast, the proposed load balancing scheme delivers the expected strong scaling behavior across a wide range of MPI ranks,
since \emph{all} MPI ranks participate in the contact evaluation.
Moreover, load balancing results in faster contact evaluation independent of the mesh size and ghosting strategy than no load balancing.
Naturally, strong scaling behavior of the contact evaluation time~$\tEvaluate$ is not affected by the choice of ghosting strategy.

For the combined time~$\tRedistribute+\tGhosting$ for redistribution and ghosting as shown in \figref{fig:ContactTwoCubesStrongScalingRedistGhostTime},
the timings are now dominated by the choice of ghosting strategy.
In particular, fully redundant ghosting of the {\master} interface (dashed lines) requires a consistently larger time across a wide range of MPI ranks.
Ghosting via binning (solid lines) can benefit from additional hardware resources,
until the strong scaling limit is reached and timings are increasing with an increasing number of MPI ranks.
The effect of the load balancing strategy is negligible,
but we note that the extra cost of performing a redistribution leads to slightly higher times with load balancing than without load balancing.

Considering the assembly of all contact terms into the global linear system,
\figref{fig:ContactTwoCubesStrongScalingModifyTime} shows just a small difference with and without load balancing.
Similar to the weak scaling study from \secref{sec:NumExTwoCubesContactWeakScaling},
the ghosting strategy does not impact these timings.
We observe good strong scalability for all studied cases.

Having in mind the overall goal of a fast {\timetosolution},
the total time~$\tContactTotal = \tRedistribute+\tGhosting+\tEvaluate+\tAssemble$
for (possibly) redistributing, ghosting, evaluation, and assembly of the contact interface
is depicted in \figref{fig:ContactTwoCubesStrongScalingTotalContactTime}.
Again, the proposed load balancing strategy results in the best timings and in good, but not perfect scaling behavior.
Stemming from the pure contact evaluation time~$\tEvaluate$ (depicted in \figref{fig:ContactTwoCubesStrongScalingContactEvaluationTime}),
the total time~$\tContactTotal$ without load balancing does not strictly follow the expected strong scaling behavior.
Per definition of~$\tContactTotal$, this diagram combines all characteristics from \figssref{fig:ContactTwoCubesStrongScalingTotalContactTime}{fig:ContactTwoCubesStrongScalingRedistGhostTime}{fig:ContactTwoCubesStrongScalingModifyTime},
namely the better scaling of~$\tEvaluate$ due to load balancing
and the increase in the timing component ~$\tRedistribute+\tGhosting$ for large numbers of MPI ranks due to increased communication and redistribution effort.

In sum, the best scaling behavior is achieved with the proposed approach of load balancing in combination with ghosting via binning.
While both components affect the overall efficiency,
the fastest evaluation times and the best weak and strong scaling behavior can only be achieved
through the \emph{combination} of load balancing with ghosting via binning.
So far, we have limited our analysis to static contact problems without any changes in the contact zone,
where the proposed algorithms demonstrate their beneficial effect on the run time and the weak and strong scaling behavior,
but could not unfold their full potential.
Therefore, we now move to dynamic contact problems,
where the contact zone changes over time and, thus, the load balancing is expected to show an even better effect on the scalability and performance.

\subsection{Rolling cylinder with dynamic contact}
\label{sec:NumExRollingCyl}

This example studies the behavior of parallel algorithms for dynamic contact problems,
{\ie} for uni-lateral contact problems where the contact zone is changing over time.
This will exercise the parallel redistribution of the contact interface discretization to its full extent.

The problem is configured as follows:
An elastic hollow cylinder is pushed onto a deformable block with initially flat surfaces.
After contact has been established, a rotating motion is imposed on the inner surface of the hollow cylinder,
somewhat mimicking a rolling tire.
Both bodies are modeled with a compressible Neo-Hooke material with
Young's modulus~$\youngs=1$, Poisson's ratio~$\poisson=0.3$, and density~$\density=10^{-6}$.

Both bodies are discretized with first-order hexahedral finite elements.
The top surface of the block is chosen as the {\master} side of the contact interface,
while the outer surface of the hollow cylinder takes the role of the {\slave} surface.
For constraint enforcement, a node-based penalty regularization of the mortar approach with a penalty parameter of~$5$ is chosen.
Time integration employs the generalized-$\alpha$ method~\cite{Chung1993a} with spectral radius~$\SpectralRadiusGenAlpha=1.0$.

\begin{figure}
\begin{center}
\subfigure[Evolution of interface DD for different load balancing strategies]{
\begin{tabular}{c|c|c|c}
~ & no LB & static LB & dynamic LB\\
\hline
\rotatebox{90}{Step 0}
& \begin{tikzpicture} \node at (0,0) {\includegraphics[width=0.27\textwidth]{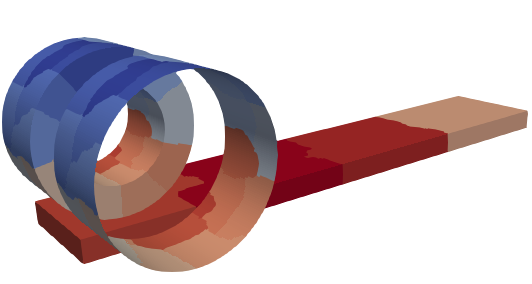}};
\end{tikzpicture}
& \begin{tikzpicture} \node at (0,0) {\includegraphics[width=0.27\textwidth]{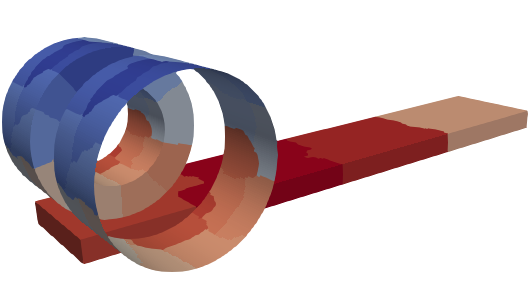}};
\end{tikzpicture}
& \begin{tikzpicture} \node at (0,0) {\includegraphics[width=0.27\textwidth]{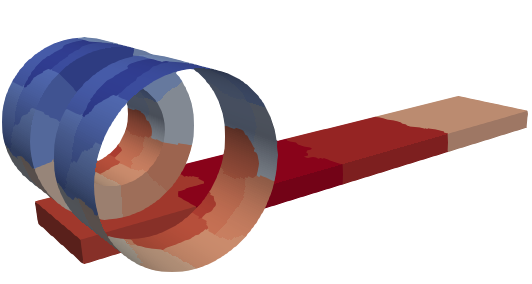}};
\end{tikzpicture}\\
\rotatebox{90}{Step 1}
& \begin{tikzpicture} \node at (0,0) {\includegraphics[width=0.27\textwidth]{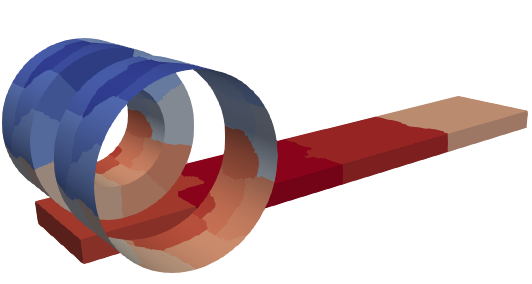}};
\end{tikzpicture}
& \begin{tikzpicture} \node at (0,0) {\includegraphics[width=0.27\textwidth]{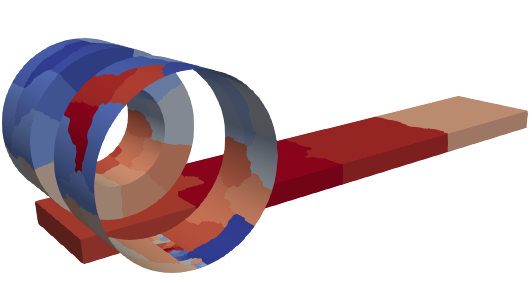}};
\begin{scope}[shift={(-0.75,-0.9)}]
\draw[ultra thick,color=red!50!orange,rotate=-27] (0.0,0.0) ellipse (12pt and 5pt);
\node [right=11pt,color=red!50!orange,circle,inner sep=2pt,fill=white,fill opacity=0.6,text opacity=1] at (0,0) {\footnotesize{I}};
\end{scope}
\end{tikzpicture}
& \begin{tikzpicture} \node at (0,0) {\includegraphics[width=0.27\textwidth]{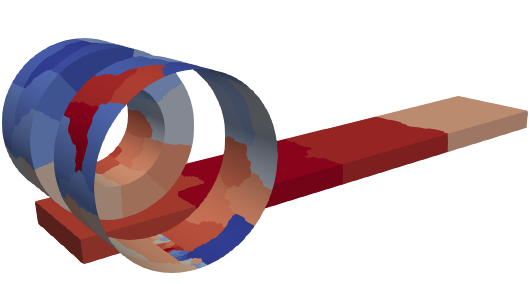}};
\end{tikzpicture}\\
\rotatebox{90}{Step 50}
& \begin{tikzpicture} \node at (0,0) {\includegraphics[width=0.27\textwidth]{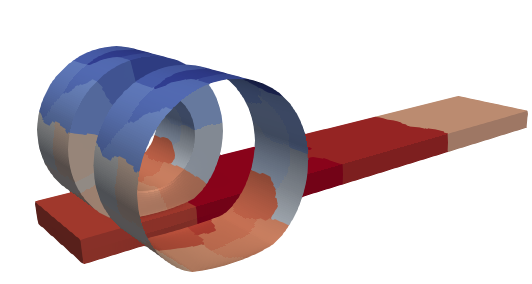}};
\end{tikzpicture}
& \begin{tikzpicture} \node at (0,0) {\includegraphics[width=0.27\textwidth]{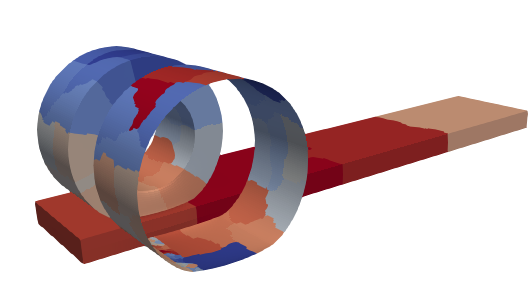}};
\end{tikzpicture}
& \begin{tikzpicture} \node at (0,0) {\includegraphics[width=0.27\textwidth]{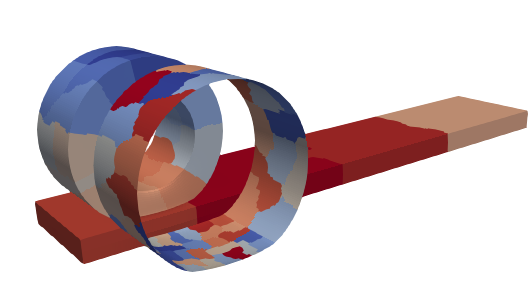}};
\begin{scope}[shift={(-0.5,-0.8)}]
\draw[ultra thick,color=red!50!orange,rotate=-25] (0.0,0.0) ellipse (14pt and 12pt);
\node [right=15pt,color=red!50!orange,circle,inner sep=2pt,fill=white,fill opacity=0.6,text opacity=1] at (0,0) {\footnotesize{II}};
\end{scope}
\end{tikzpicture}\\
\rotatebox{90}{Step 150}
& \begin{tikzpicture} \node at (0,0) {\includegraphics[width=0.27\textwidth]{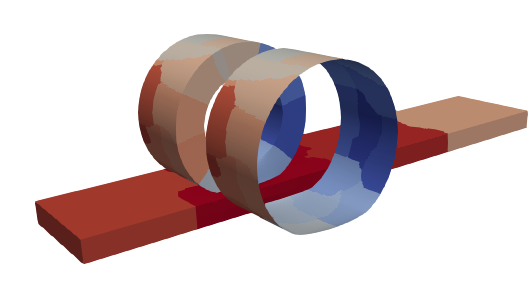}};
\end{tikzpicture}
& \begin{tikzpicture} \node at (0,0) {\includegraphics[width=0.27\textwidth]{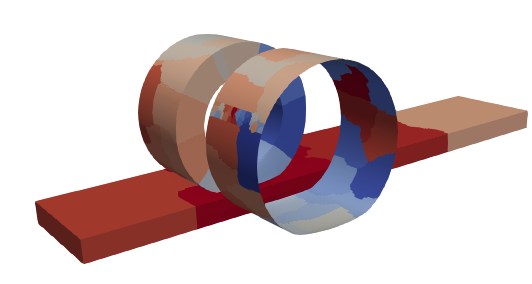}};
\begin{scope}[shift={(-0.26,0.23)}]
\draw[ultra thick,color=red!50!orange,rotate=-22] (0.0,0.0) ellipse (10pt and 5pt);
\node [right=12pt,color=red!50!orange,circle,inner sep=2pt,fill=white,fill opacity=0.6,text opacity=1] at (0,0) {\footnotesize{III}};
\end{scope}
\end{tikzpicture}
& \begin{tikzpicture} \node at (0,0) {\includegraphics[width=0.27\textwidth]{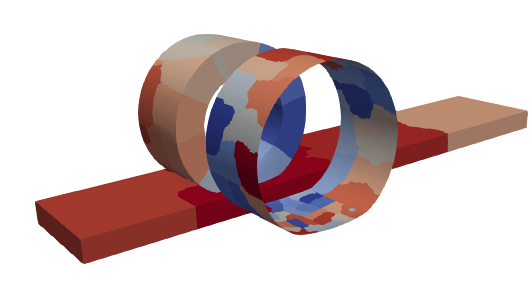}};
\begin{scope}[shift={(0.39,-0.55)}]
\draw[ultra thick,color=red!50!orange,rotate=-22] (0.0,0.0) ellipse (12pt and 10pt);
\node [right=15pt,color=red!50!orange,circle,inner sep=2pt,fill=white,fill opacity=0.6,text opacity=1] at (0,0) {\footnotesize{IV}};
\end{scope}
\end{tikzpicture}
\\
\multicolumn{4}{c}{\input{data/cyl_viz/o/cyl_viz_Static_binning_proximityYes_sizeDDTol-1-03_maxImbalance-1-8_minEleProc-0_mpi024_owner_colorbar.tex}}
\end{tabular}
} 

\subfigure[Close-up of feature I]{
\includegraphics[trim={1cm 0.5cm 9cm 1.5cm},clip,height=0.22\textwidth]{./data/cyl_viz/o/cyl_viz_Static_binning_proximityYes_sizeDDTol-1-03_maxImbalance-1-8_minEleProc-0_mpi024_step001_owner_cylinder.png}
\label{fig:RollingCylFeatI}
} 
\subfigure[Close-up of feature II]{
\includegraphics[trim={3cm 0.5cm 8cm 1.5cm},clip,height=0.22\textwidth]{./data/cyl_viz/o/cyl_viz_Dynamic_binning_proximityYes_sizeDDTol-1-03_maxImbalance-1-8_minEleProc-0_mpi024_step050_owner_cylinder.png}
\label{fig:RollingCylFeatII}
} 
\subfigure[Close-up of feature III]{
\includegraphics[trim={7cm 1.5cm 5cm 1.0cm},clip,height=0.22\textwidth]{./data/cyl_viz/o/cyl_viz_Static_binning_proximityYes_sizeDDTol-1-03_maxImbalance-1-8_minEleProc-0_mpi024_step150_owner_cylinder.png}
\label{fig:RollingCylFeatIII}
} 
\subfigure[Close-up of feature IV]{
\includegraphics[trim={7cm 1.5cm 5cm 1.0cm},clip,height=0.22\textwidth]{./data/cyl_viz/o/cyl_viz_Dynamic_binning_proximityYes_sizeDDTol-1-03_maxImbalance-1-8_minEleProc-0_mpi024_step150_owner_cylinder.png}
\label{fig:RollingCylFeatIV}
} 
\caption{Visualization of volume and interface subdomains for different load balancing strategies in a dynamic contact example.
Interesting features are highlighted with roman numbers I - IV and discussed in the text.}
\label{fig:CylinderDDVisualization}
\end{center}
\end{figure}

\Figref{fig:CylinderDDVisualization} exemplarily compares the volume and interface subdomains for the different load balancing strategies for the case of 24 MPI ranks.
The initial subdomain layout in step~0 is the same for all load balancing strategies.
While the DD of the underlying bodies will not be altered,
we apply the interface load balancing scheme proposed in \secref{sec:InterfaceDD},
which results in different interface DDs for the {\slave} side.
To unclutter the presenation, we only show the evolution of the {\slave} side's interface DDs,
since this is the key ingredient for a scalable mortar evaluation.
Interesting features due to load balancing are highlighted with roman numbers I - IV (see also \figsrangeref{fig:RollingCylFeatI}{fig:RollingCylFeatIV}) and will be discussed below.
In the case of no load balancing (column ``no LB''),
the interface subdomains match the subdomains of the underlying volume DD throughout the entire simulation.
For static LB (column ``static LB''),
an initial interface DD is performed at the beginning of step~1, but it is not updated during the simulation.
Hence, a small strip of {\slave} subdomains is generated during the initial load balancing phase ({\cf} highlight I or \figref{fig:RollingCylFeatI})
and then rotates with the rolling motion of the cylinder as marked by highlight III (see also \figref{fig:RollingCylFeatIII}),
such that it quickly leaves the contact area and, thus, does not contribute to an optimal contact evaluation throughout the entire simulation.
Based on the threshold criterion~\eqref{eq:TriggerRebalancing},
the dynamic load balancing (column ``dyn. LB'') updates the interface DD close to the contact area ({\cf} highlight II or \figref{fig:RollingCylFeatII})
such that the interface DD is nearly optimal in the vicinity of the contact zone and
all processes participate in the evaluation of the contact terms independent of the rolling motion of the cylinder ({\cf} highlight IV or \figref{fig:RollingCylFeatIV}).

\subsubsection{Effect of load balancing on wall clock time and memory consumption}
\label{sec:NumExRollingCylSavings}

We compare the cases of no load balancing, an initial load balancing in the reference configuration,
and the dynamic load balancing proposed in \secref{sec:InterfaceDDDynamic}
on a mesh with 825,600 hexahedral elements consisting of 913,923 nodes and resulting in 2,741,769 displacement unknowns.
We run the simulation on 96 MPI ranks on our {\inhouse} cluster.
For the case of initial and dynamic load balancing,
we limit the relative mismatch in subdomain size of the interface DD
by setting the {\zoltan} parameter \emph{IMBALANCE\_TOL} to  $1.03$~\cite{Boman2012a}.
For dynamic load balancing,
we have tested different thresholds~$\imbalanceThresholdTime \in \{1.01,\allowbreak  1.2,\allowbreak  1.5,\allowbreak  1.8,\allowbreak  2.5,\allowbreak  5.0,\allowbreak  8.0\}$
to trigger rebalancing,
but will only report and discuss selected cases in the following for the sake of presentation,
namely~$\imbalanceThresholdTime \in \{1.01,\allowbreak 1.8,\allowbreak  5.0,\allowbreak  8.0\}$.
To extend the ghosting of the {\master} side's interface discretization,
we rely on the \emph{binning} strategy outlined in \secref{sec:Binning}.
A comparison of the different ghosting strategies is presented in \secref{sec:NumExRollingCylGhostingStrategies}.

We run the simulation for 200 time steps (20 time steps to close the initial gap, then 180 time steps of the rolling motion)
to facilitate a rotation of $180^\circ$,
such that the contact area on the outer cylinder surface substantially moves along the circumferential direction.

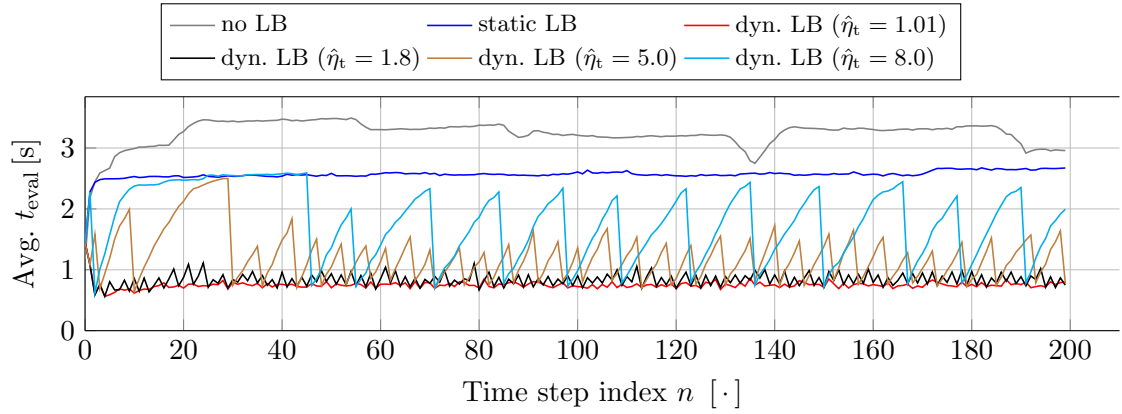
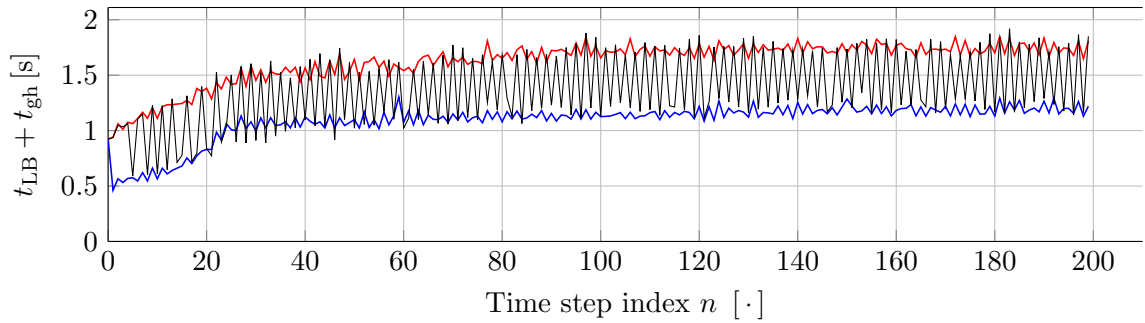
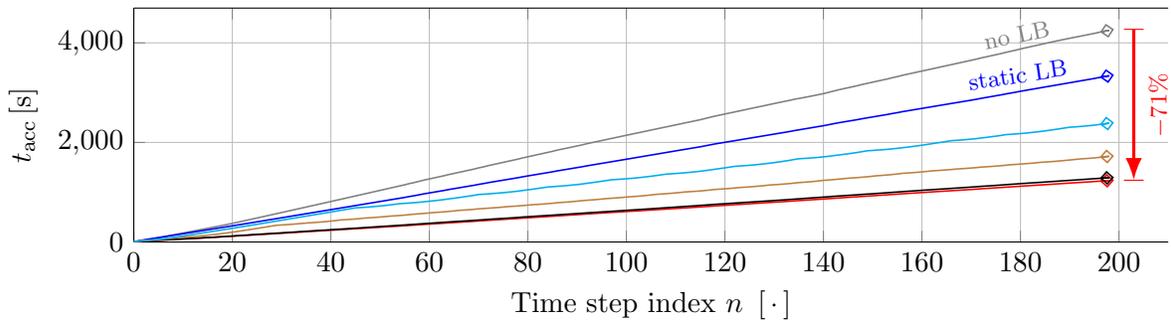
\begin{figure}

\tikzstyle{binningContactEvalTimeNone}=[solid, semithick, color=gray]
\tikzstyle{binningContactEvalTimeStatic}=[solid, semithick, color=blue]
\tikzstyle{binningContactEvalTimeDynamicOne}=[solid, semithick, color=red]
\tikzstyle{binningContactEvalTimeDynamicFour}=[solid, semithick, color=black]
\tikzstyle{binningContactEvalTimeDynamicSix}=[solid, semithick, color=brown]
\tikzstyle{binningContactEvalTimeDynamicSeven}=[solid, semithick, color=cyan]

\begin{center}

\begin{tikzpicture}
\begin{axis}[%
  hide axis,
  xmin=0,
  xmax=1,
  ymin=0,
  ymax=1,
  legend cell align=left,
  legend style={font=\footnotesize, at={(0.5,1.2)}, anchor=center},
  legend columns = 3
  ]

\addlegendimage{binningContactEvalTimeNone}
\addlegendentry{no LB}
\addlegendimage{binningContactEvalTimeStatic}
\addlegendentry{static LB}
\addlegendimage{binningContactEvalTimeDynamicOne}
\addlegendentry{dyn. LB ($\imbalanceThresholdTime = 1.01$)}
\addlegendimage{binningContactEvalTimeDynamicFour}
\addlegendentry{dyn. LB ($\imbalanceThresholdTime = 1.8$)}
\addlegendimage{binningContactEvalTimeDynamicSix}
\addlegendentry{dyn. LB ($\imbalanceThresholdTime = 5.0$)}
\addlegendimage{binningContactEvalTimeDynamicSeven}
\addlegendentry{dyn. LB ($\imbalanceThresholdTime = 8.0$)}

\end{axis}
\end{tikzpicture}

\subfigure[Avgerage contact evaluation time per time step]{
\label{fig:ComparisonLoadBalancingContactTimeAverage}
\begin{tikzpicture}

\pgfplotstableread{data/savings/cyl_savings_None_binning_proximityNo_sizeDDTol-1.03_maxImbalance-0.0_minEleProc-0_mpi096_contact_eval_time.data}\binningContactEvalTimeNone
\pgfplotstableread{data/savings/cyl_savings_Static_binning_proximityYes_sizeDDTol-1.03_maxImbalance-0.0_minEleProc-0_mpi096_contact_eval_time.data}\binningContactEvalTimeStatic

\pgfplotstableread{data/savings/cyl_savings_Dynamic_binning_proximityYes_sizeDDTol-1.03_maxImbalance-1.01_minEleProc-0_mpi096_contact_eval_time.data}\binningContactEvalTimeDynamicOne
\pgfplotstableread{data/savings/cyl_savings_Dynamic_binning_proximityYes_sizeDDTol-1.03_maxImbalance-1.8_minEleProc-0_mpi096_contact_eval_time.data}\binningContactEvalTimeDynamicFour
\pgfplotstableread{data/savings/cyl_savings_Dynamic_binning_proximityYes_sizeDDTol-1.03_maxImbalance-5.0_minEleProc-0_mpi096_contact_eval_time.data}\binningContactEvalTimeDynamicSix
\pgfplotstableread{data/savings/cyl_savings_Dynamic_binning_proximityYes_sizeDDTol-1.03_maxImbalance-8.0_minEleProc-0_mpi096_contact_eval_time.data}\binningContactEvalTimeDynamicSeven

\begin{axis}[scale only axis, axis y line*=left, 
  xmin=0,
  xmax=210,
  ymin=0,
  width=0.8\textwidth, height=0.13\textheight,
  ylabel={Avg. $\tEvaluate\left[\second\right]$},
  xlabel={Time step index~$\indTimeStep~\left[\cdot\right]$},
  grid=major]

\addplot[binningContactEvalTimeNone] table[x = time_step, y = avg_contact_eval_time] from \binningContactEvalTimeNone;
\addplot[binningContactEvalTimeStatic] table[x = time_step, y = avg_contact_eval_time] from \binningContactEvalTimeStatic;
\addplot[binningContactEvalTimeDynamicOne] table[x = time_step, y = avg_contact_eval_time] from \binningContactEvalTimeDynamicOne;
\addplot[binningContactEvalTimeDynamicFour] table[x = time_step, y = avg_contact_eval_time] from \binningContactEvalTimeDynamicFour;
\addplot[binningContactEvalTimeDynamicSix] table[x = time_step, y = avg_contact_eval_time] from \binningContactEvalTimeDynamicSix;
\addplot[binningContactEvalTimeDynamicSeven] table[x = time_step, y = avg_contact_eval_time] from \binningContactEvalTimeDynamicSeven;


\end{axis}
\end{tikzpicture}
} 

\subfigure[Time for ghosting of {\master} interface plus potentially for load balancing (selection of data sets)]{
\label{fig:ComparisonLoadBalancingContactTimeRedistribution}
\begin{tikzpicture}

\pgfplotstableread{data/savings/cyl_savings_None_binning_proximityNo_sizeDDTol-1.03_maxImbalance-0.0_minEleProc-0_mpi096_contact_eval_time.data}\binningContactEvalTimeNone
\pgfplotstableread{data/savings/cyl_savings_Static_binning_proximityYes_sizeDDTol-1.03_maxImbalance-0.0_minEleProc-0_mpi096_contact_eval_time.data}\binningContactEvalTimeStatic

\pgfplotstableread{data/savings/cyl_savings_Dynamic_binning_proximityYes_sizeDDTol-1.03_maxImbalance-1.01_minEleProc-0_mpi096_contact_eval_time.data}\binningContactEvalTimeDynamicOne
\pgfplotstableread{data/savings/cyl_savings_Dynamic_binning_proximityYes_sizeDDTol-1.03_maxImbalance-1.8_minEleProc-0_mpi096_contact_eval_time.data}\binningContactEvalTimeDynamicFour
\pgfplotstableread{data/savings/cyl_savings_Dynamic_binning_proximityYes_sizeDDTol-1.03_maxImbalance-5.0_minEleProc-0_mpi096_contact_eval_time.data}\binningContactEvalTimeDynamicSix
\pgfplotstableread{data/savings/cyl_savings_Dynamic_binning_proximityYes_sizeDDTol-1.03_maxImbalance-8.0_minEleProc-0_mpi096_contact_eval_time.data}\binningContactEvalTimeDynamicSeven

\begin{axis}[scale only axis, axis y line*=left, 
  xmin=0,
  xmax=210,
  ymin=0,
  width=0.8\textwidth, height=0.13\textheight,
  ylabel={$\tRedistribute + \tGhosting\left[\second\right]$},
  xlabel={Time step index~$\indTimeStep~\left[\cdot\right]$},
  grid=major]

\addplot[binningContactEvalTimeStatic] table[x = time_step, y = time_for_redist_ghosting] from \binningContactEvalTimeStatic;
\addplot[binningContactEvalTimeDynamicOne] table[x = time_step, y = time_for_redist_ghosting] from \binningContactEvalTimeDynamicOne;
\addplot[binningContactEvalTimeDynamicFour, thin] table[x = time_step, y = time_for_redist_ghosting] from \binningContactEvalTimeDynamicFour;


\end{axis}
\end{tikzpicture}
} 

\subfigure[Accumulated contact time]{
\label{fig:ComparisonLoadBalancingContactTimeAccumulated}
\begin{tikzpicture}

\pgfplotstableread{data/savings/cyl_savings_None_binning_proximityNo_sizeDDTol-1.03_maxImbalance-0.0_minEleProc-0_mpi096_contact_eval_time.data}\binningContactEvalTimeNone
\pgfplotstableread{data/savings/cyl_savings_Static_binning_proximityYes_sizeDDTol-1.03_maxImbalance-0.0_minEleProc-0_mpi096_contact_eval_time.data}\binningContactEvalTimeStatic

\pgfplotstableread{data/savings/cyl_savings_Dynamic_binning_proximityYes_sizeDDTol-1.03_maxImbalance-1.01_minEleProc-0_mpi096_contact_eval_time.data}\binningContactEvalTimeDynamicOne
\pgfplotstableread{data/savings/cyl_savings_Dynamic_binning_proximityYes_sizeDDTol-1.03_maxImbalance-1.8_minEleProc-0_mpi096_contact_eval_time.data}\binningContactEvalTimeDynamicFour
\pgfplotstableread{data/savings/cyl_savings_Dynamic_binning_proximityYes_sizeDDTol-1.03_maxImbalance-5.0_minEleProc-0_mpi096_contact_eval_time.data}\binningContactEvalTimeDynamicSix
\pgfplotstableread{data/savings/cyl_savings_Dynamic_binning_proximityYes_sizeDDTol-1.03_maxImbalance-8.0_minEleProc-0_mpi096_contact_eval_time.data}\binningContactEvalTimeDynamicSeven

\begin{axis}[scale only axis, axis y line*=left, 
  xmin=0,
  xmax=210,
  ymin=0,
  width=0.8\textwidth, height=0.13\textheight,
  ylabel={$\tAccumulated\left[\second\right]$},
  xlabel={Time step index~$\indTimeStep~\left[\cdot\right]$},
  grid=major]

\addplot[binningContactEvalTimeNone,{-Turned Square[open]}] table[x = time_step, y = accumulated_contact_time] from \binningContactEvalTimeNone;
\addplot[binningContactEvalTimeStatic,{-Turned Square[open]}] table[x = time_step, y = accumulated_contact_time] from \binningContactEvalTimeStatic;
\addplot[binningContactEvalTimeDynamicOne,{-Turned Square[open]}] table[x = time_step, y = accumulated_contact_time] from \binningContactEvalTimeDynamicOne;
\addplot[binningContactEvalTimeDynamicFour,{-Turned Square[open]}] table[x = time_step, y = accumulated_contact_time] from \binningContactEvalTimeDynamicFour;
\addplot[binningContactEvalTimeDynamicSix,{-Turned Square[open]}] table[x = time_step, y = accumulated_contact_time] from \binningContactEvalTimeDynamicSix;
\addplot[binningContactEvalTimeDynamicSeven,{-Turned Square[open]}] table[x = time_step, y = accumulated_contact_time] from \binningContactEvalTimeDynamicSeven;

\node [above,color=gray,rotate=12] at (axis cs:180,3800) {\footnotesize no LB};
\node [above,color=blue,rotate=8.5] at (axis cs:180,3000) {\footnotesize static LB};

\draw [color=red] (axis cs:201,4271.8) -- (axis cs:205,4271.8);
\draw [color=red] (axis cs:201,1240.7) -- (axis cs:205,1240.7);
\draw [color=red,-Latex,very thick] (axis cs:203,4271.8) -- (axis cs:203,1240.7);

\node [below,rotate=90,color=red] at (axis cs:204,2756.2) {\footnotesize $-71\%$};

\end{axis}

\end{tikzpicture}
} 

\caption{Effect of different load balancing strategies on the time spent in the contact evaluation}
\label{fig:ComparisonLoadBalancingContactTime}
\end{center}
\end{figure}

For every time step,
\Figref{fig:ComparisonLoadBalancingContactTimeAverage} reports the average time per nonlinear iteration spent in contact evaluation
(without considering the cost for load balancing).
If no load balancing is performed,
the average contact evaluation time is the largest.
Since the {\slave} side's interface DD is just adopted from the underlying volume discretization,
some processes do never participate in contact evaluation.
Moreover, the number of processes contributing to the contact evaluation changes over time,
so the average contact evaluation time also changes over time steps.
In contrast, static load balancing assures that all parallel processes hold their share of the {\slave} side of the interface,
such that the average contact evaluation time is roughly constant for all time steps (as soon as full contact is established).
Since only a part of all processes contributes to the evaluation of the potentially active part of the {\slave} interface,
the average contact evaluation time is still rather large.
Ultimately, dynamic load balancing triggers a rebalancing based on the current simulation status
to aid a well-balanced distribution of the contact evaluation work to \emph{all} parallel processes.
In \figref{fig:ComparisonLoadBalancingContactTimeAverage},
time steps just after a drop in~$\tEvaluate$ are those, in which a rebalancing has occurred.
Since the effort of mortar evaluation is now distributed to all processes,
the time spent in mortar evaluation drops significantly on average.
Of course, individual time steps with an imbalanced work distribution among processes might take longer,
which will ultimately lead to rebalancing as soon as the rebalancing criterion~\eqref{eq:TriggerRebalancing} is met.
In particular, a very low rebalancing threshold ({\eg} $\imbalanceThresholdTime = 1.01$) requires to rebalance in basically every time step.
Although this results in the overall fastest mortar evaluation,
the additional effort for rebalancing limits the possible {\speedup}.
On the other hand, a loose threshold ({\eg} $\imbalanceThresholdTime \in \{5.0, 8.0\}$)
triggers the rebalancing only a few times over the course of the simulation,
however for some time steps the time spent in the mortar evaluation can grow by a factor of two or even three compared to the ideal case.
In our numerical experiments,
we have found the threshold~$\imbalanceThresholdTime = 1.8$ to deliver a good compromise
between imbalance in per-process workload and the frequency of rebalancing.
Therefore, we will use this threshold value for all further studies.

\Figref{fig:ComparisonLoadBalancingContactTimeRedistribution} reports the time~$\tGhosting$ spent for ghosting of the {\master} side of the contact interface
plus the time~$\tRedistribute$ for rebalancing of the interface DD (if applicable).
For the clarity of the presentation,
we concentrate on three selected cases.
While the time component~$\tGhosting$ for the {\master} side ghosting is rather constant for all three cases,
the time component~$\tRedistribute$ varies:
For static load balancing,
only the first time step requires rebalancing, while all later time steps do not perform load balancing anymore.
Hence, this curve peaks in the first time step and then drops and remains at low values.
For dynamic load balancing with the strict imbalance threshold~$\imbalanceThresholdTime=1.01$,
rebalancing occurs in every time step,
such that this case consistently delivers high values for~$\tRedistribute+\tGhosting$.
Obviously, these two cases can be interpreted as a lower and upper bound as evident from \figref{fig:ComparisonLoadBalancingContactTimeRedistribution}.
The case of dynamic load balancing with~$\imbalanceThresholdTime=1.8$ positions itself in between,
since some time steps require rebalancing, but some do not.

While all cases with dynamic load balancing spend additional time on the redistribution of the interface subdomains,
these additional timings are easily amortized.
To this end,
\figref{fig:ComparisonLoadBalancingContactTimeAccumulated} plots the
time~$\tAccumulated = \sum_{\indContactEvalEvent = 1}^{\numContactEvalEvents} (\tEvaluate+\tRedistribute+\tGhosting)_{\indContactEvalEvent}, \indContactEvalEvent \in \{1,\hdots,\numContactEvalEvents\},$
of all time components related to mortar evaluation over all time steps
accrued over all $\numContactEvalEvents$ contact evaluations of the entire simulation.
The end point markers are intended to highlight also small differences between curves.
Naturally, a strict monotone increase is expected,
while one aims for an as low as possible slope.
Similar to the average contact evaluation time,
static load balancing is beneficial compared to no load balancing at all,
while dynamic load balancing results in the lowest contact evaluation times.
Clearly, the better parallelization due to the dynamic load balancing strategy
pays off the additional cost for occasional rebalancing.
The lower the acceptable imbalance~$\imbalanceThresholdTime$ is,
the lower is the accumulated contact evaluation time~$\tAccumulated$.
Overall, a maximum reduction up to $71\%$ in~$\tAccumulated$ can be achieved through proper dynamic load balancing.
We note that the difference in~$\tAccumulated$ between~$\imbalanceThresholdTime=1.01$ and~$\imbalanceThresholdTime=1.8$ is very small,
indicating that load balancing in every time step does not bring much additional value.

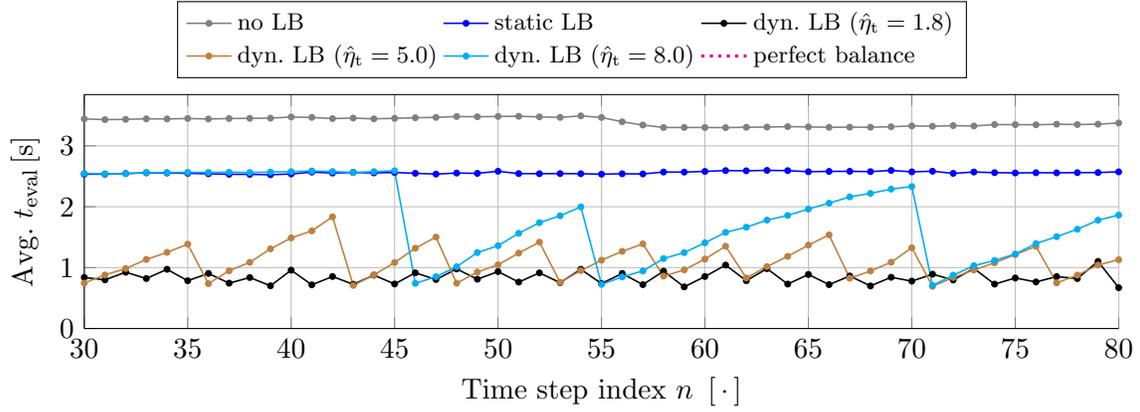
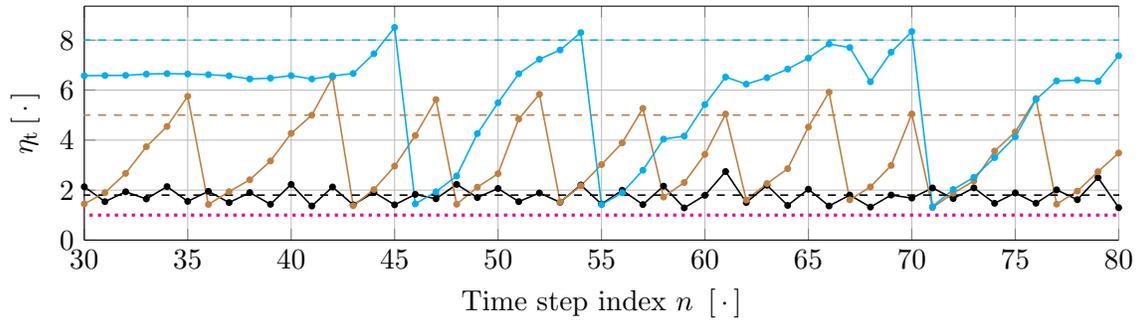
\begin{figure}

\tikzstyle{perfectBalance}=[dotted, very thick, color=magenta]
\tikzstyle{binningContactEvalTimeNone}=[solid, semithick, color=gray, mark=*, mark size=1pt]
\tikzstyle{binningContactEvalTimeStatic}=[solid, semithick, color=blue, mark=*, mark size=1pt]
\tikzstyle{binningContactEvalTimeDynamicOne}=[solid, semithick, color=red, mark=*, mark size=1pt]
\tikzstyle{binningContactEvalTimeDynamicFour}=[solid, semithick, color=black, mark=*, mark size=1pt]
\tikzstyle{binningContactEvalTimeDynamicSix}=[solid, semithick, color=brown, mark=*, mark size=1pt]
\tikzstyle{binningContactEvalTimeDynamicSeven}=[solid, semithick, color=cyan, mark=*, mark size=1pt]

\begin{center}

\begin{tikzpicture}
\begin{axis}[%
  hide axis,
  xmin=0,
  xmax=1,
  ymin=0,
  ymax=1,
  legend cell align=left,
  legend style={font=\footnotesize, at={(0.5,1.2)}, anchor=center},
  legend columns = 3
  ]

\addlegendimage{binningContactEvalTimeNone}
\addlegendentry{no LB}
\addlegendimage{binningContactEvalTimeStatic}
\addlegendentry{static LB}
\addlegendimage{binningContactEvalTimeDynamicFour}
\addlegendentry{dyn. LB ($\imbalanceThresholdTime = 1.8$)}
\addlegendimage{binningContactEvalTimeDynamicSix}
\addlegendentry{dyn. LB ($\imbalanceThresholdTime = 5.0$)}
\addlegendimage{binningContactEvalTimeDynamicSeven}
\addlegendentry{dyn. LB ($\imbalanceThresholdTime = 8.0$)}
\addlegendimage{perfectBalance}
\addlegendentry{perfect balance}

\end{axis}
\end{tikzpicture}

\subfigure[Avgerage contact evaluation time per time step]{
\label{fig:ComparisonLoadBalancingContactTimeAverageCloseUp}
\begin{tikzpicture}

\pgfplotstableread{data/savings/cyl_savings_None_binning_proximityNo_sizeDDTol-1.03_maxImbalance-0.0_minEleProc-0_mpi096_contact_eval_time.data}\binningContactEvalTimeNone
\pgfplotstableread{data/savings/cyl_savings_Static_binning_proximityYes_sizeDDTol-1.03_maxImbalance-0.0_minEleProc-0_mpi096_contact_eval_time.data}\binningContactEvalTimeStatic

\pgfplotstableread{data/savings/cyl_savings_Dynamic_binning_proximityYes_sizeDDTol-1.03_maxImbalance-1.8_minEleProc-0_mpi096_contact_eval_time.data}\binningContactEvalTimeDynamicFour
\pgfplotstableread{data/savings/cyl_savings_Dynamic_binning_proximityYes_sizeDDTol-1.03_maxImbalance-5.0_minEleProc-0_mpi096_contact_eval_time.data}\binningContactEvalTimeDynamicSix
\pgfplotstableread{data/savings/cyl_savings_Dynamic_binning_proximityYes_sizeDDTol-1.03_maxImbalance-8.0_minEleProc-0_mpi096_contact_eval_time.data}\binningContactEvalTimeDynamicSeven

\begin{axis}[scale only axis, axis y line*=left, 
  xmin=30,
  xmax=80,
  ymin=0,
  width=0.8\textwidth, height=0.13\textheight,
  ylabel={Avg. $\tEvaluate\left[\second\right]$},
  xlabel={Time step index~$\indTimeStep~\left[\cdot\right]$},
  grid=major]

\addplot[binningContactEvalTimeNone] table[x = time_step, y = avg_contact_eval_time] from \binningContactEvalTimeNone;
\addplot[binningContactEvalTimeStatic] table[x = time_step, y = avg_contact_eval_time] from \binningContactEvalTimeStatic;
\addplot[binningContactEvalTimeDynamicFour] table[x = time_step, y = avg_contact_eval_time] from \binningContactEvalTimeDynamicFour;
\addplot[binningContactEvalTimeDynamicSix] table[x = time_step, y = avg_contact_eval_time] from \binningContactEvalTimeDynamicSix;
\addplot[binningContactEvalTimeDynamicSeven] table[x = time_step, y = avg_contact_eval_time] from \binningContactEvalTimeDynamicSeven;


\end{axis}
\end{tikzpicture}
} 

\subfigure[Max/min ratio~$\imbalanceRatioTime$ in contact evaluation time across all parallel processes (solid lines)
for different imbalance thresholds~$\imbalanceThresholdTime$ (dashed lines)]{
\label{fig:ComparisonImbalanceThresholdTimeImbalanceCloseUp}
\begin{tikzpicture}

\pgfplotstableread{data/savings/cyl_savings_None_binning_proximityNo_sizeDDTol-1.03_maxImbalance-0.0_minEleProc-0_mpi096_parallel_imbalance.data}\binningContactEvalTimeNone
\pgfplotstableread{data/savings/cyl_savings_Static_binning_proximityYes_sizeDDTol-1.03_maxImbalance-0.0_minEleProc-0_mpi096_parallel_imbalance.data}\binningContactEvalTimeStatic

\pgfplotstableread{data/savings/cyl_savings_Dynamic_binning_proximityYes_sizeDDTol-1.03_maxImbalance-1.8_minEleProc-0_mpi096_parallel_imbalance.data}\binningContactEvalTimeDynamicFour
\pgfplotstableread{data/savings/cyl_savings_Dynamic_binning_proximityYes_sizeDDTol-1.03_maxImbalance-5.0_minEleProc-0_mpi096_parallel_imbalance.data}\binningContactEvalTimeDynamicSix
\pgfplotstableread{data/savings/cyl_savings_Dynamic_binning_proximityYes_sizeDDTol-1.03_maxImbalance-8.0_minEleProc-0_mpi096_parallel_imbalance.data}\binningContactEvalTimeDynamicSeven

\begin{axis}[scale only axis, axis y line*=left, 
  xmin=30,
  xmax=80,
  ymin=0,
  width=0.8\textwidth, height=0.13\textheight,
  ylabel={$\imbalanceRatioTime\left[\cdot\right]$},
  xlabel={Time step index~$\indTimeStep~\left[\cdot\right]$},
  grid=major
  ]

\addplot[binningContactEvalTimeDynamicFour] table[x = time_step, y = imbalance_time] from \binningContactEvalTimeDynamicFour;
\addplot[binningContactEvalTimeDynamicSix] table[x = time_step, y = imbalance_time] from \binningContactEvalTimeDynamicSix;
\addplot[binningContactEvalTimeDynamicSeven] table[x = time_step, y = imbalance_time] from \binningContactEvalTimeDynamicSeven;

\addplot[perfectBalance] coordinates {(0,1.0) (1000, 1.0)};
\addplot[binningContactEvalTimeDynamicFour, dashed] coordinates {(0,1.8) (1000, 1.8)};
\addplot[binningContactEvalTimeDynamicSix, dashed] coordinates {(0,5.0) (1000, 5.0)};
\addplot[binningContactEvalTimeDynamicSeven, dashed] coordinates {(0,8.0) (1000, 8.0)};

\end{axis}
\end{tikzpicture}
} 

\caption{Detailed view of contact evaluation timings and its imbalance}
\label{fig:ComparisonImbalanceThresholdTimeCloseup}
\end{center}
\end{figure}

To demonstrate the effect of the rebalancing trigger~$\imbalanceThresholdTime$ in detail,
\figref{fig:ComparisonImbalanceThresholdTimeCloseup} shows a close-up of the results in \figref{fig:ComparisonLoadBalancingContactTimeAverage}
as well as the evolution of the max/min ratio~$\imbalanceRatioTime$ in contact evaluation time across all parallel processes.
For a clearer visualization, only a subset of the results is plotted.
In \figref{fig:ComparisonLoadBalancingContactTimeAverageCloseUp},
data points after a drop in~$\tEvaluate$ correspond to time steps,
where load balancing has occured since the max/min ratio~$\imbalanceRatioTime$
exceeded the threshold~$\imbalanceThresholdTime$ in the previous time step.
This is in line with~\figref{fig:ComparisonImbalanceThresholdTimeImbalanceCloseUp},
where~$\imbalanceRatioTime$ is plotted over time along with dashed lines to indicate the different thresholds~$\imbalanceThresholdTime$.
We observe that~$\imbalanceRatioTime$ drops close to the perfect balance ({\ie} $\imbalanceRatioTime = 1.0$)
just after it exceeded the threshold level~$\imbalanceThresholdTime$.
In favor of an uncluttered view,
\figref{fig:ComparisonImbalanceThresholdTimeImbalanceCloseUp} shows only results obtained with dynamic load balancing.

So far, we have studied the impact of the load balancing strategy and the imbalance threshold~$\imbalanceThresholdTime$
onto the time spent in the computational treatment of all mortar terms.
In all cases, dynamic load balancing is worth the effort.
Since the present example shows very good behavior for~$\imbalanceThresholdTime=1.8$,
we continue to use this value throughout this example.
We note that the optimal choice of~$\imbalanceThresholdTime$ is problem-dependent.
Yet, we generally recommend to use dynamic load balancing for contact problems with changing contact zones
and select~$\imbalanceThresholdTime$ on a case-by-case basis.

\subsubsection{Strong scaling behavior under dynamic load balancing}
\label{sec:NumExRollingCylStrongScalability}

Now, we study the strong scaling behavior of the contact evaluation time when dynamic load balancing is active.
We therefore study two different problem sizes:
517,185 displacement unknowns referred to as ``$500k$''
and 1,005,993 displacement degrees of freedom denoted by ``$1000k$''.
While we keep the problem sizes fixed,
we solve the problem on an increasing number of MPI ranks 
on our {\inhouse} cluster.

We will compare different load balancing strategies, namely no load balancing (``no LB''),
initial load balancing in the reference configuration (``static LB''),
and dynamic load balancing (with~$\imbalanceThresholdTime=1.8$ (``dyn. LB'') as found useful in \secref{sec:NumExRollingCylSavings}).

\begin{figure}

\tikzstyle{fiveKNone}=[very thick, dotted, blue, mark=triangle, mark options={scale=0.9, solid, thick}]
\tikzstyle{fiveKStatic}=[thick, dashed, blue, mark=square, mark options={scale=0.9, solid, thick}]
\tikzstyle{fiveKDynamic}=[thick, solid, blue, mark=diamond, mark options={scale=1.2, solid, thick}]

\tikzstyle{oneMNone}=[very thick, dotted, red, mark=triangle, mark options={scale=0.9, solid, thick}]
\tikzstyle{oneMStatic}=[thick, dashed, red, mark=square, mark options={scale=0.9, solid, thick}]
\tikzstyle{oneMDynamic}=[thick, solid, red, mark=diamond, mark options={scale=1.2, solid, thick}]

\tikzstyle{perfectScaling}=[thick, dashed, gray]

\begin{center}

\begin{tikzpicture}
\begin{axis}[%
  hide axis,
  xmin=0,
  xmax=1,
  ymin=0,
  ymax=1,
  legend cell align=left,
  legend style={font=\footnotesize, at={(0.5,1.2)}, anchor=center},
  legend columns = 3
  ]

\addlegendimage{fiveKNone}
\addlegendentry{$500k$ / no LB};
\addlegendimage{fiveKStatic}
\addlegendentry{$500k$ / static LB};
\addlegendimage{fiveKDynamic}
\addlegendentry{$500k$ / dyn. LB};

\addlegendimage{oneMNone}
\addlegendentry{$1000k$ / no LB};
\addlegendimage{oneMStatic}
\addlegendentry{$1000k$ / static LB};
\addlegendimage{oneMDynamic}
\addlegendentry{$1000k$ / dyn. LB};

\addlegendimage{perfectScaling}
\addlegendentry{perfect scaling};

\end{axis}
\end{tikzpicture}

\subfigure[Avg. evaluation time per time step]{\label{fig:CylStrongScalingContactEvalTime}
\begin{tikzpicture}

\pgfplotstableread{data/cyl_strong_scaling/cyl_strong_redistNone_500k_scaling.data}\scalingFiveKNone
\pgfplotstableread{data/cyl_strong_scaling/cyl_strong_redistNone_1000k_scaling.data}\scalingOneMNone

\pgfplotstableread{data/cyl_strong_scaling/cyl_strong_redistStatic_500k_scaling.data}\scalingFiveKStatic
\pgfplotstableread{data/cyl_strong_scaling/cyl_strong_redistStatic_1000k_scaling.data}\scalingOneMStatic

\pgfplotstableread{data/cyl_strong_scaling/cyl_strong_500k_scaling.data}\scalingFiveK
\pgfplotstableread{data/cyl_strong_scaling/cyl_strong_1000k_scaling.data}\scalingOneM

\begin{axis}[scale only axis,
  xmode=log,
  ymode=log,
  width=0.35\textwidth, height=0.25\textwidth,
  ylabel={Avg. $\tEvaluate\left[\second\right]$},
  xlabel={\#MPI ranks~$\left[-\right]$},
  grid=major
  ]

\addplot[fiveKNone] table[x = num_procs, y = avg_contact_time] from \scalingFiveKNone;
\addplot[oneMNone] table[x = num_procs, y = avg_contact_time] from \scalingOneMNone;

\addplot[fiveKStatic] table[x = num_procs, y = avg_contact_time] from \scalingFiveKStatic;
\addplot[oneMStatic] table[x = num_procs, y = avg_contact_time] from \scalingOneMStatic;

\addplot[fiveKDynamic] table[x = num_procs, y = avg_contact_time] from \scalingFiveK;
\addplot[oneMDynamic] table[x = num_procs, y = avg_contact_time] from \scalingOneM;

\addplot[perfectScaling] coordinates {(1,9) (10,0.9)};

\end{axis}
\end{tikzpicture}
} 
\subfigure[Avg. redistribution and ghosting time per time step]{\label{fig:CylStrongScalingRedistributeTime}
\begin{tikzpicture}

\pgfplotstableread{data/cyl_strong_scaling/cyl_strong_redistNone_500k_scaling.data}\scalingFiveKNone
\pgfplotstableread{data/cyl_strong_scaling/cyl_strong_redistNone_1000k_scaling.data}\scalingOneMNone

\pgfplotstableread{data/cyl_strong_scaling/cyl_strong_redistStatic_500k_scaling.data}\scalingFiveKStatic
\pgfplotstableread{data/cyl_strong_scaling/cyl_strong_redistStatic_1000k_scaling.data}\scalingOneMStatic

\pgfplotstableread{data/cyl_strong_scaling/cyl_strong_500k_scaling.data}\scalingFiveK
\pgfplotstableread{data/cyl_strong_scaling/cyl_strong_1000k_scaling.data}\scalingOneM

\begin{axis}[scale only axis,
  xmode=log,
  ymode=log,
  width=0.35\textwidth, height=0.25\textwidth,
  ylabel={Avg. $\tRedistribute + \tGhosting\left[\second\right]$},
  xlabel={\#MPI ranks~$\left[-\right]$},
  grid=major
  ]

\addplot[fiveKNone] table[x = num_procs, y = avg_redist_ghost_time] from \scalingFiveKNone;
\addplot[oneMNone] table[x = num_procs, y = avg_redist_ghost_time] from \scalingOneMNone;

\addplot[fiveKStatic] table[x = num_procs, y = avg_redist_ghost_time] from \scalingFiveKStatic;
\addplot[oneMStatic] table[x = num_procs, y = avg_redist_ghost_time] from \scalingOneMStatic;

\addplot[fiveKDynamic] table[x = num_procs, y = avg_redist_ghost_time] from \scalingFiveK;
\addplot[oneMDynamic] table[x = num_procs, y = avg_redist_ghost_time] from \scalingOneM;

\addplot[perfectScaling] coordinates {(1,10) (10,1)};

\end{axis}
\end{tikzpicture}
} 

\subfigure[Avg. total contact time per time step]{\label{fig:CylStrongScalingTotalContactTime}
\begin{tikzpicture}

\pgfplotstableread{data/cyl_strong_scaling/cyl_strong_redistNone_500k_scaling.data}\scalingFiveKNone
\pgfplotstableread{data/cyl_strong_scaling/cyl_strong_redistNone_1000k_scaling.data}\scalingOneMNone

\pgfplotstableread{data/cyl_strong_scaling/cyl_strong_redistStatic_500k_scaling.data}\scalingFiveKStatic
\pgfplotstableread{data/cyl_strong_scaling/cyl_strong_redistStatic_1000k_scaling.data}\scalingOneMStatic

\pgfplotstableread{data/cyl_strong_scaling/cyl_strong_500k_scaling.data}\scalingFiveK
\pgfplotstableread{data/cyl_strong_scaling/cyl_strong_1000k_scaling.data}\scalingOneM

\begin{axis}[scale only axis,
  xmode=log,
  ymode=log,
  width=0.35\textwidth, height=0.25\textwidth,
  ylabel={Avg. $\tEvaluate+\tRedistribute+\tGhosting\left[\second\right]$},
  xlabel={\#MPI ranks~$\left[-\right]$},
  grid=major
  ]

\addplot[fiveKNone] table[x = num_procs, y = avg_total_contact_time] from \scalingFiveKNone;
\addplot[oneMNone] table[x = num_procs, y = avg_total_contact_time] from \scalingOneMNone;

\addplot[fiveKStatic] table[x = num_procs, y = avg_total_contact_time] from \scalingFiveKStatic;
\addplot[oneMStatic] table[x = num_procs, y = avg_total_contact_time] from \scalingOneMStatic;

\addplot[fiveKDynamic] table[x = num_procs, y = avg_total_contact_time] from \scalingFiveK;
\addplot[oneMDynamic] table[x = num_procs, y = avg_total_contact_time] from \scalingOneM;

\addplot[perfectScaling] coordinates {(1,80) (10,8)};

\end{axis}
\end{tikzpicture}
} 
\subfigure[Accumulated contact time]{\label{fig:CylStrongScalingAccumulatedTime}
\begin{tikzpicture}
\pgfplotstableread{data/cyl_strong_scaling/cyl_strong_redistNone_500k_scaling.data}\scalingFiveKNone
\pgfplotstableread{data/cyl_strong_scaling/cyl_strong_redistNone_1000k_scaling.data}\scalingOneMNone

\pgfplotstableread{data/cyl_strong_scaling/cyl_strong_redistStatic_500k_scaling.data}\scalingFiveKStatic
\pgfplotstableread{data/cyl_strong_scaling/cyl_strong_redistStatic_1000k_scaling.data}\scalingOneMStatic

\pgfplotstableread{data/cyl_strong_scaling/cyl_strong_500k_scaling.data}\scalingFiveK
\pgfplotstableread{data/cyl_strong_scaling/cyl_strong_1000k_scaling.data}\scalingOneM

\begin{axis}[scale only axis,
  xmode=log,
  ymode=log,
  width=0.35\textwidth, height=0.25\textwidth,
  ylabel={$\tAccumulated~\left[\second\right]$},
  xlabel={\#MPI ranks~$\left[-\right]$},
  grid=major
  ]

\addplot[fiveKNone] table[x = num_procs, y = accumulated_contact_time] from \scalingFiveKNone;
\addplot[oneMNone] table[x = num_procs, y = accumulated_contact_time] from \scalingOneMNone;

\addplot[fiveKStatic] table[x = num_procs, y = accumulated_contact_time] from \scalingFiveKStatic;
\addplot[oneMStatic] table[x = num_procs, y = accumulated_contact_time] from \scalingOneMStatic;

\addplot[fiveKDynamic] table[x = num_procs, y = accumulated_contact_time] from \scalingFiveK;
\addplot[oneMDynamic] table[x = num_procs, y = accumulated_contact_time] from \scalingOneM;

\addplot[perfectScaling] coordinates {(1,2000) (10,200)};

\end{axis}
\end{tikzpicture}
} 

\caption{Strong scaling of the contact timings under different load balancing strategies}
\label{fig:ContactStrongScaling}
\end{center}
\end{figure}
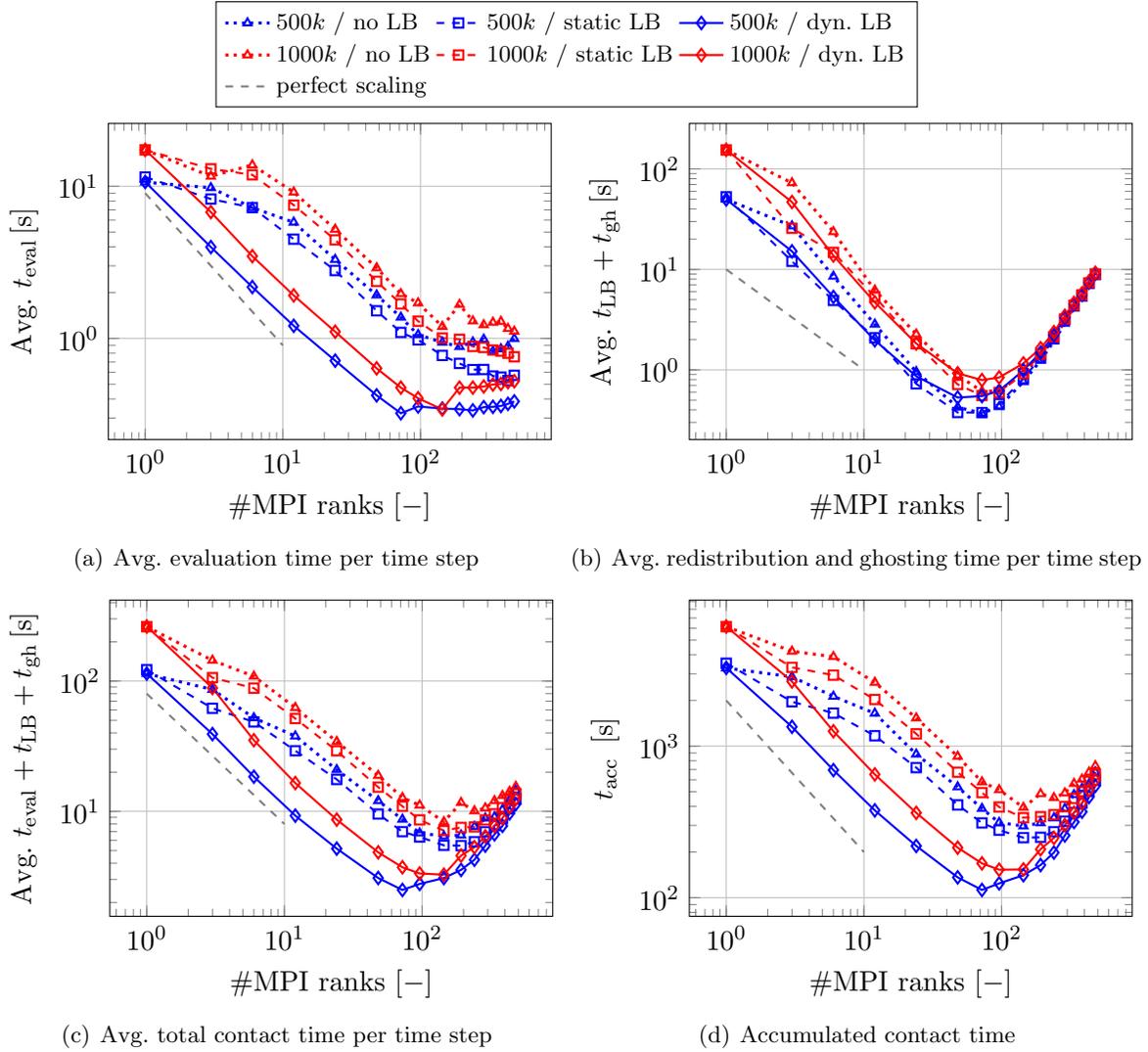

\figref{fig:ContactStrongScaling} shows the strong scaling behavior.
Again, we consider the average contact evaluation time~$\tEvaluate$ per time step,
the time~$\tRedistribute+\tGhosting$ spent in redistribution and ghosting of the interface discretizations,
the average total time~$\tContactTotal = \tEvaluate+\tRedistribute+\tGhosting$ per time step,
and finally the total contact time~$\tAccumulated = \sum_{\indContactEvalEvent = 1}^{\numContactEvalEvents} (\tEvaluate+\tRedistribute+\tGhosting)_{\indContactEvalEvent}, \indContactEvalEvent \in \{1,\hdots,\numContactEvalEvents\},$
accumulated over all $\numContactEvalEvents$ contact evaluations of the entire simulation.
For both problem sizes as well as all quantities of interest,
we observe good strong scaling behavior when using dynamic load balancing:
starting from a small number of MPI ranks,
the time spent on a given task ({\eg.} contact evaluation, redistribution and ghosting, total contact time, accumulated contact time)
is reduced when adding more MPI ranks to tackle the computations,
while the reduction rate is linked to the increase in MPI ranks,
{\ie} delivering perfect strong scaling~\cite{Amdahl1967a}.
As expected, both meshes reached their strong scaling limit at some point,
such that adding more hardware resources does not reduce, but actually increase the execution time,
{\eg} due to a deteriorating computation-to-communication ratio.
Naturally, the strong scaling limit of the large problem ($1000k$) is located at twice the number of MPI ranks
as for the small, half-sized problem ($500k$).
The beneficial effect of dynamic load balancing becomes evident in comparison to ``no LB'' and ``static LB'':
Without any load balancing or just an initial rebalancing of the interface discretizations,
the initial slope in the scaling diagrams is far from optimal.
Once again, this originates from the curse of dimensionality,
since the additional hardware resources do not necessarily participate in the interface evaluation.
For an intermediate number of MPI ranks,
strong scaling is recovered, however absolute timings are much higher than for the same setup with dynamic load balancing.
As already observed in \secref{sec:NumExTwoCubesContactWeakScaling},
static load balancing is consistently a bit faster than no using load balancing at all,
yet it is by far slower than dynamic load balancing.

As demonstrated, the proposed dynamic load balancing scheme is the key factor
to achieve strong scalability of the evaluation of mortar terms in a {\nonlinear} and time-dependent contact simulation.
To the authors' best knowledge,
this constitutes the first time that strong scalability in such a complex setting could be demonstrated.

\subsubsection{Comparison of strategies to extend the {\master} side's ghosting}
\label{sec:NumExRollingCylGhostingStrategies}

While the influence of the load balancing strategy has already been discussed previously,
we now aim to assess the impact of the ghosting strategy on the overall performance of the contact evaluation.
Therefore, we exemplarily consider the mesh from \secref{sec:NumExRollingCylSavings} run on 96 MPI ranks.
Now, we compare the fully redundant storage of the {\master} side of the interface ({\cf} \secref{sec:MasterRedundantStorage})
to the geometrically motivated binning approach ({\cf} \secref{sec:Binning}).
We study again the cases of no, static, and dynamic load balancing.
For the clarity of the presenation,
we only show the case of dynamic load balancing scenario with~$\imbalanceThresholdTime = 1.8$,
but note that other values for~$\imbalanceThresholdTime$ exhibit similar behavior.

\Figref{fig:ComparisonGhostingContactTime} summarizes the wall clock time spent on contact evaluation.
\begin{figure}

\tikzstyle{noneRedundant}=[dotted, thick, gray]
\tikzstyle{staticRedundant}=[dotted, thick, blue]
\tikzstyle{dynamicRedundant}=[dotted, thick, red]

\tikzstyle{noneBinning}=[solid, semithick, gray]
\tikzstyle{staticBinning}=[solid, semithick, blue]
\tikzstyle{dynamicBinning}=[solid, semithick, red]

\begin{center}

\begin{tikzpicture}
\begin{axis}[%
  hide axis,
  xmin=10,
  xmax=50,
  ymin=0,
  ymax=0.4,
  legend cell align=left,
  legend style={font=\footnotesize, at={(0.5,1.2)}, anchor=center},
  legend columns = 3
  ]

\addlegendimage{noneRedundant}
\addlegendentry{no LB / redundant ghosting};
\addlegendimage{staticRedundant}
\addlegendentry{static LB / redundant ghosting};
\addlegendimage{dynamicRedundant}
\addlegendentry{dyn. LB ($\imbalanceThresholdTime = 1.8$) / redundant ghosting};

\addlegendimage{noneBinning}
\addlegendentry{no LB / binning};
\addlegendimage{staticBinning}
\addlegendentry{static LB / binning};
\addlegendimage{dynamicBinning}
\addlegendentry{dyn. LB ($\imbalanceThresholdTime = 1.8$) / binning};

\end{axis}
\end{tikzpicture}

\subfigure[Avgerage contact evaluation time per time step]{
\label{fig:ComparisonGhostingContactTimeAverage}
\begin{tikzpicture}

\pgfplotstableread{data/savings/cyl_savings_None_binning_proximityNo_sizeDDTol-1.03_maxImbalance-0.0_minEleProc-0_mpi096_contact_eval_time.data}\binningContactEvalTimeNone
\pgfplotstableread{data/savings/cyl_savings_Static_binning_proximityYes_sizeDDTol-1.03_maxImbalance-0.0_minEleProc-0_mpi096_contact_eval_time.data}\binningContactEvalTimeStatic
\pgfplotstableread{data/savings/cyl_savings_Dynamic_binning_proximityYes_sizeDDTol-1.03_maxImbalance-1.8_minEleProc-0_mpi096_contact_eval_time.data}\binningContactEvalTimeDynamic

\pgfplotstableread{data/savings/cyl_savings_None_redundant_master_proximityNo_sizeDDTol-1.03_maxImbalance-0.0_minEleProc-0_mpi096_contact_eval_time.data}\redundantMasterContactEvalTimeNone
\pgfplotstableread{data/savings/cyl_savings_Static_redundant_master_proximityYes_sizeDDTol-1.03_maxImbalance-0.0_minEleProc-0_mpi096_contact_eval_time.data}\redundantMasterContactEvalTimeStatic
\pgfplotstableread{data/savings/cyl_savings_Dynamic_redundant_master_proximityYes_sizeDDTol-1.03_maxImbalance-1.8_minEleProc-0_mpi096_contact_eval_time.data}\redundantMasterContactEvalTimeDynamic

\begin{axis}[scale only axis, axis y line*=left,
  xmin=0,
  xmax=210,
  width=0.8\textwidth, height=0.13\textheight,
  ylabel={Avg. $\tEvaluate\left[\second\right]$},
  xlabel={Time step index~$\indTimeStep~\left[\cdot\right]$},
  grid=major]

\addplot[noneBinning] table[x = time_step, y = avg_contact_eval_time] from \binningContactEvalTimeNone;
\addplot[staticBinning] table[x = time_step, y = avg_contact_eval_time] from \binningContactEvalTimeStatic;
\addplot[dynamicBinning] table[x = time_step, y = avg_contact_eval_time] from \binningContactEvalTimeDynamic;

\addplot[noneRedundant] table[x = time_step, y = avg_contact_eval_time] from \redundantMasterContactEvalTimeNone;
\addplot[staticRedundant] table[x = time_step, y = avg_contact_eval_time] from \redundantMasterContactEvalTimeStatic;
\addplot[dynamicRedundant] table[x = time_step, y = avg_contact_eval_time] from \redundantMasterContactEvalTimeDynamic;

\end{axis}
\end{tikzpicture}
} 

\subfigure[Time for ghosting of {\master} interface plus potentially for load balancing]{
\label{fig:ComparisonGhostingContactTimeRedistribution}
\begin{tikzpicture}

\pgfplotstableread{data/savings/cyl_savings_None_binning_proximityNo_sizeDDTol-1.03_maxImbalance-0.0_minEleProc-0_mpi096_contact_eval_time.data}\binningContactEvalTimeNone
\pgfplotstableread{data/savings/cyl_savings_Static_binning_proximityYes_sizeDDTol-1.03_maxImbalance-0.0_minEleProc-0_mpi096_contact_eval_time.data}\binningContactEvalTimeStatic
\pgfplotstableread{data/savings/cyl_savings_Dynamic_binning_proximityYes_sizeDDTol-1.03_maxImbalance-1.8_minEleProc-0_mpi096_contact_eval_time.data}\binningContactEvalTimeDynamic

\pgfplotstableread{data/savings/cyl_savings_None_redundant_master_proximityNo_sizeDDTol-1.03_maxImbalance-0.0_minEleProc-0_mpi096_contact_eval_time.data}\redundantMasterContactEvalTimeNone
\pgfplotstableread{data/savings/cyl_savings_Static_redundant_master_proximityYes_sizeDDTol-1.03_maxImbalance-0.0_minEleProc-0_mpi096_contact_eval_time.data}\redundantMasterContactEvalTimeStatic
\pgfplotstableread{data/savings/cyl_savings_Dynamic_redundant_master_proximityYes_sizeDDTol-1.03_maxImbalance-1.8_minEleProc-0_mpi096_contact_eval_time.data}\redundantMasterContactEvalTimeDynamic

\begin{axis}[scale only axis, axis y line*=left,
  xmin=0,
  xmax=210,
  width=0.8\textwidth, height=0.13\textheight,
  ylabel={$\tRedistribute + \tGhosting\left[\second\right]$},
  xlabel={Time step index~$\indTimeStep~\left[\cdot\right]$},
  grid=major]

\addplot[dynamicBinning] table[x = time_step, y = time_for_redist_ghosting] from \binningContactEvalTimeDynamic;

\addplot[dynamicRedundant] table[x = time_step, y = time_for_redist_ghosting] from \redundantMasterContactEvalTimeDynamic;

\end{axis}
\end{tikzpicture}
} 

\subfigure[Contact time~$\tAccumulated$ accumulated over all time steps]{
\label{fig:ComparisonGhostingContactTimeAccumulated}
\begin{tikzpicture}

\pgfplotstableread{data/savings/cyl_savings_None_binning_proximityNo_sizeDDTol-1.03_maxImbalance-0.0_minEleProc-0_mpi096_contact_eval_time.data}\binningContactEvalTimeNone
\pgfplotstableread{data/savings/cyl_savings_Static_binning_proximityYes_sizeDDTol-1.03_maxImbalance-0.0_minEleProc-0_mpi096_contact_eval_time.data}\binningContactEvalTimeStatic
\pgfplotstableread{data/savings/cyl_savings_Dynamic_binning_proximityYes_sizeDDTol-1.03_maxImbalance-1.8_minEleProc-0_mpi096_contact_eval_time.data}\binningContactEvalTimeDynamic

\pgfplotstableread{data/savings/cyl_savings_None_redundant_master_proximityNo_sizeDDTol-1.03_maxImbalance-0.0_minEleProc-0_mpi096_contact_eval_time.data}\redundantMasterContactEvalTimeNone
\pgfplotstableread{data/savings/cyl_savings_Static_redundant_master_proximityYes_sizeDDTol-1.03_maxImbalance-0.0_minEleProc-0_mpi096_contact_eval_time.data}\redundantMasterContactEvalTimeStatic
\pgfplotstableread{data/savings/cyl_savings_Dynamic_redundant_master_proximityYes_sizeDDTol-1.03_maxImbalance-1.8_minEleProc-0_mpi096_contact_eval_time.data}\redundantMasterContactEvalTimeDynamic

\begin{axis}[scale only axis, axis y line*=left,
  xmin=0,
  xmax=210,
  width=0.8\textwidth, height=0.13\textheight,
  ylabel={$\tAccumulated\left[\second\right]$},
  xlabel={Time step index~$\indTimeStep~\left[\cdot\right]$},
  grid=major]

\addplot[noneBinning] table[x = time_step, y = accumulated_contact_time] from \binningContactEvalTimeNone;
\addplot[staticBinning] table[x = time_step, y = accumulated_contact_time] from \binningContactEvalTimeStatic;
\addplot[dynamicBinning] table[x = time_step, y = accumulated_contact_time] from \binningContactEvalTimeDynamic;

\addplot[noneRedundant] table[x = time_step, y = accumulated_contact_time] from \redundantMasterContactEvalTimeNone;
\addplot[staticRedundant] table[x = time_step, y = accumulated_contact_time] from \redundantMasterContactEvalTimeStatic;
\addplot[dynamicRedundant] table[x = time_step, y = accumulated_contact_time] from \redundantMasterContactEvalTimeDynamic;

\draw [color=red] (axis cs:201,2161.4) -- (axis cs:205,2161.4);
\draw [color=red] (axis cs:201,1297.3) -- (axis cs:205,1297.3);
\draw [color=red,-Latex,very thick] (axis cs:203,2161.4) -- (axis cs:203,1297.3);

\node [below,rotate=90,color=red] at (axis cs:204,1729.3) {\footnotesize $-40\%$};

\end{axis}
\end{tikzpicture}
} 

\caption{Effect of ghosting strategies on the contact timings:
the combination of dynamic load balancing with ghosting via binning consistently delivers the fastest timings for contact evaluation.}
\label{fig:ComparisonGhostingContactTime}
\end{center}
\end{figure}
For the pure contact evaluation time reported in \figref{fig:ComparisonGhostingContactTimeAverage},
the fully redundant ghosting increases the evaluation time for all cases,
since the contact detection needs to account for all master elements,
while ghosting via binning pre-sorts the master elements based on their geometric proximity within neighboring bins.

\Figref{fig:ComparisonGhostingContactTimeRedistribution} depicts the time spent in redistribution and ghosting of interface data.
For the sake of a clear presentation and to really focus on the most relevant case,
we show only the curves for dynamic load balancing.
Evidently, ghosting via binning is faster by a factor of $\approx8-10\times$ than fully redundant ghosting.

\Figref{fig:ComparisonGhostingContactTimeAccumulated} shows the accumulated time for contact evaluation, load balancing, and ghosting,
{\ie} $\tAccumulated = \sum_{\indContactEvalEvent = 1}^{\numContactEvalEvents} (\tEvaluate+\tRedistribute+\tGhosting)_{\indContactEvalEvent}, \indContactEvalEvent \in \{1,\hdots,\numContactEvalEvents\},$
to assess the overall accumulated time spent on all $\numContactEvalEvents$ evaluations of the contact interface over the course of the entire simulation.
For the cases with no and static load balancing,
the ghosting strategy does not impact the overall performance significantly.
For dynamic load balancing though,
the necessity of ghosting after each redistribution makes the difference:
the performance difference between fully redundant ghosting and ghosting via binning as observed in \figref{fig:ComparisonGhostingContactTimeRedistribution}
now accumulates over time,
such that the use of binning results in the overall lowest time spent on contact evaluation.
So, additional savings of $40\%$ of the contact evaluation time can be achieved.
Summing up the study of contact timings,
the best case scenario of dynamic load balancing with ghosting via binning is faster than
\begin{itemize}
\item dynamic load balancing with fully redundant ghosting by a factor of $\approx1.67$,
\item static load balancing by a factor of $\approx2.61$,
\item no load balancing by a factor of $\approx3.30$,
\end{itemize}
which strongly emphasizes the benefits of dynamic load balancing and ghosting via binning in dynamic contact problems.

Finally,
we briefly summarize the impact of the load balancing scheme and the ghosting strategy onto the cost for storage and parallel communication:
If no load balancing is performed (``no LB''),
the maximum number of owned nodes per process is roughly $10\times$ larger than its average across all processes,
since not all processes hold a portion of the interface.
This imbalance is alleviated for static or dynamic load balancing.
Regarding the impact of the ghosting strategy,
ghosting via binning reduces down the number of nodes/elements to be ghosted by a factor of $100\times$ compared to the fully redundant case,
which ultimately also impacts the global memory footprint of the application.

In sum, dynamic contact problems require a good choice of load balancing strategy as well as a suitable ghosting strategy.
In particular, load balancing highly impacts the time spent in contact evaluation.
Despite the additional cost of performing the load balancing operation,
the overall fastest contact evaluation is achieved with dynamic load balancing based on a user-given imbalance threshold~$\imbalanceThresholdTime$.
While we have found~$\imbalanceThresholdTime = 1.8$ to deliver very good results in our numerical studies,
the optimal choice of~$\imbalanceThresholdTime$ can depend on details of the computing hardware, the software implementation, and also the example at hand.
To reduce the amount of parallel communication as well as the memory demand per compute node,
ghosting via binning is by far superior to a fully redundant storage of the master side of the interface discretization.
The overall best performance with respect to both phenomena (run time and communicatio/memory demand) is obtained
through the combination of dynamic load balancing with ghosting via binning.

\section{Concluding remarks}
\label{sec:Conclusion}

Recognizing the tremendous computational effort to evaluate mortar integrals in the context of {\nonmatching} interface discretizations
as they exemplarily arise in contact mechanics,
this paper proposes strategies for efficient storage and parallel computational kernels for mortar interface problems.
Starting from a close look at the computational effort to evaluate mortar integrals,
we have derived two basic requirements for computations on parallel machines with distributed memory architecture:
On the one hand, one needs to enable access to the appropriate interface data to guarantee a correct identification of all {\master}/{\slave} pairs at the mortar interface.
On the other hand, the available parallel hardware needs to be used efficiently,
such that parallel scalability of the mortar evaluation can be achieved.

We have discussed some techniques to guarantee access to all required {\master}/{\slave} pairs during the contact search and mortar evaluation.
While fully redundant ghosting is conceptually easy and straightforward to implement,
it suffers from an elevated memory demands and tremendous communication overhead at large scale,
which ultimately increases the {\timetosolution}.
A geometrically motivated approach using a background grid of Cartesian bins
allows for the efficient identification of nearby {\master} elements,
reduces the per-process memory demand
as well as limits the ghosting data to the nearby {\master} elements.
The binning approach has shown the best timings in large weak and strong scaling studies
and consistently reduces the amount of data to be communicated between parallel processes as well as to be stored within a process.

We have then discussed the curse of dimensionality in overlapping DDs of interface problems,
which requires a special treatment of the interface subdomains.
To this end,
we have proposed to use an interface DD independent from the underlying volume DD
and were able to demonstrate optimal weak and strong scalability of the mortar evaluation time.
To account for dynamic changes in the contact zone,
we have designed a dynamic load balancing scheme for contact problems,
which tracks imbalances among parallel processes and rebalances the computational work
as soon as user-given imbalance thresholds are exceeded.
We have tested the proposed algorithms on a time-dependent {\nonlinear} contact problem undergoing large deformations.
In time measurements on such large-scale examples,
dynamic load balancing outperforms the case of no or only initial load balancing by factors up to $2\!-\!4\times$.
Wall clock time is the lowest,
when only small imbalances are allowed,
although even a large imbalance tolerance delivers faster computations than simulations without any load balancing at all.
For the first time, strong and weak scalability could be shown for time-dependent {\nonlinear} contact problems
undergoing large deformations and dynamically evolving contact zones
through the application of the proposed dynamic load balancing scheme.

In our numerical experiments,
we have studied representative test cases from computational contact mechanics.
We have performed weak and strong scaling studies up to 480 MPI ranks
as well as have assessed the impact of different algorithmic parameters.
From our numerical experiments,
we extract several findings:
\begin{itemize}
\item Ghosting via binning is favorable due to its reduced communication overhead, which also directly reduces the {\timetosolution}.
\item Load balancing is crucial for optimal contact evaluation times.
In particular, systems with a static contact zone benefit from an initial redistribution of the interface,
while contact problems with dynamically evolving contact zones require the proposed dynamic load balancing scheme for optimal performance.
\item For static contact problems, we have found the combination of static load balancing and ghosting via binning to deliver the best results.
\item For dynamic contact problems, we have found the combination of dynamic load balancing and ghosting via binning to deliver the best results.
\end{itemize}
In sum, we recommend to apply static load balancing in combination with ghosting via binning for problems with static contact zones,
while dynamic load balancing in combination with ghosting via binning is preferable for problems with dynamically evolving contact zones.
Following these recommendations,
a fast {\timetosolution} as well as good weak and strong scaling behavior can be achieved.

\section*{Acknowledgements}

This work has been partially funded by the Deutsche Forschungsgemeinschaft (DFG – German Research Foundation)
within the project ``Experimental characterization and numerical simulation of the automated fiber placement (AFP) process for thermoplastic fiber-reinforced plastics''
(project number: 325153381).

\bibliographystyle{abbrv}   
\bibliography{bib_contact_parallel_algos}

\end{document}